\documentclass[11pt]{article}
\usepackage{amsthm}
\usepackage[margin=0.85in]{geometry}

\usepackage{amssymb}
\usepackage{lineno}
\usepackage{tikz-cd}
\usepackage{fancyhdr}
\newcommand{\Qed}{\hfill $\blacksquare$}
\date{}

\begin{document}

\title{Rank and Duality in Representation Theory\thanks{These notes grow out of the 19th Takagi lectures [Howe17-2] delivered by the second named author on July 8-9, 2017, at the the Research Institute for Mathematical Sciences Kyoto University, Japan.}}
\author{Shamgar Gurevich\thanks{S. Gurevich, Department of Mathematics, University of Wisconsin, Madison, WI 53706, USA. e-mail: shamgar@math.wisc.edu.} \, and Roger Howe\thanks{Howe R., Wm. Kenan Jr. Professor of Mathematics, Emeritus, Yale University, New Haven, CT, USA. e-mail: roger.howe@yale.edu.}}
\lhead{Gurevich and Howe}
\rhead{Rank and Duality}
\maketitle

{\it Abstract.} There is both theoretical and numerical evidence that the set of irreducible representations of a reductive group over local or finite fields is naturally partitioned into families according to analytic properties of representations. Examples of such properties are the rate of decay at infinity of ``matrix coefficients" in the local field setting, and the order of magnitude of ``character ratios" in the finite field situation. 

In these notes we describe known results, new results, and conjectures in the theory of ``size" of representations of classical groups over finite fields (when correctly stated, most of them hold also in the local field setting), whose ultimate goal is to classify the above mentioned families of representations and accordingly to estimate the relevant analytic properties of each family.  

Specifically, we treat two main issues: the first is the introduction of a rigorous definition of a notion of size for representations of classical groups, and the second issue is a method to construct and obtain information on each family of representation of a given size. 

In particular, we propose several compatible notions of size that we call {\it $U$-rank, tensor rank and asymptotic rank}, and we develop a method called {\it eta correspondence} to construct the families of representation of each given rank.

Rank suggests a new way to organize the representations of classical groups over finite and local fields - a way in which the building blocks are the "smallest" representations. This is in contrast to Harish-Chandra's philosophy of cusp forms that is the main organizational principle since the 60s, and in it the building blocks are the cuspidal representations which are, in some sense, the "largest". The philosophy of cusp forms is well adapted to establishing the Plancherel formula for reductive groups over local fields, and led to Lusztig's classification of the irreducible representations of such groups over finite fields. However, the understanding of certain analytic properties, such as those mentioned above, seems to require a different approach.

\vskip .2 in

{\bf 0. Introduction}

\vskip .15 in

In these notes we describe a new approach, initiated in [Gurevich-Howe15], for the study of representations of classical groups over finite and local fields, an approach that seems relevant to harmonic analysis. The goal (see examples in [Gurevich-Howe17], [Gurevich-Howe18]) is to develop an effective theory of ``size" for representations, namely to give a precise definition of that notion and develop a method to analyze each family of representations having the same size. 

\vskip .1 in

{\bf 0.1. Motivation.}

\vskip .1 in

The motivation in the finite setting comes from the fact that many problems about finite groups (e.g., properties of random walks, word maps, Cayley graphs, etc.) can be approached using harmonic analysis [Diaconis88], [Diaconis96], [Diaconis-Shahshahani81], [Frobenius1896], [Gluck97], [Guranick-Malle14], [Hildebrand92], [Larsen-Shalev-Tiep11], [Liebeck17], [Liebeck-O'Brien-Shalev-Tiep10], [Liebeck-Shalev05], [Malle14], [Saloff-Coste04], [Shalev07], [Shalev09], [Shalev15]. More precisely, what intervenes in such problems are the {\it character ratios (CRs)}
$$
\frac{\chi_\pi(g)}{\dim(\pi)}, \,\,\, g \in G,   \eqno (0.1.1)
$$ 
where $\pi$ is an irreducible representation (irrep) of the relevant group $G$, and $\chi_\pi(g)=trace(\pi(g)), \, \, g \in G$, is its character. 

In general, it is not feasible to compute the CRs (0.1.1) exactly, but for applications it often suffices to show that the CRs are small for most representations or that certain sums of them exhibit cancellations. 

Looking on numerical character tables of the finite classical groups, one notices that the CRs of the irreducible representations (irreps) tend to form a certain hierarchy according to order of magnitude. Hence, for applications it might be of importance to understand better these natural families of irreps. 

In these notes we introduce several (conjecturally highly compatible)
notions of ``rank" for representations of classical groups. These notions
make the concept of size of a representation rigorous, and they seem, at least in certain cases, to put together representations of similar size of character ratios. 
  
\begin{figure}[!htb]\centering
\includegraphics[height=3.2in, width=4.2in]
{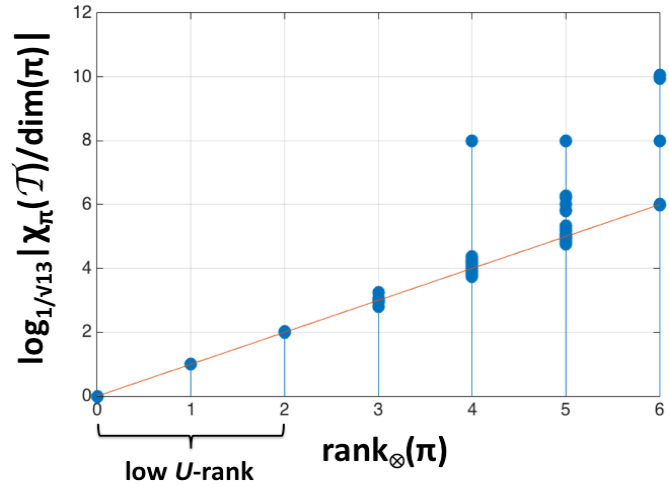}
\caption{$\log _{1/\protect\sqrt{13}}$-scale of CRs at transvection $%
\mathcal{T}$ vs. tensor rank for irreps $\protect\pi $ of $Sp_{6}(\mathbb{F}%
_{13})$.}\label{cr-rank_t-sp6_13}
\end{figure} 

\begin{figure}[!htb]\centering 
\includegraphics[height=2.9in, width=3.9in]
{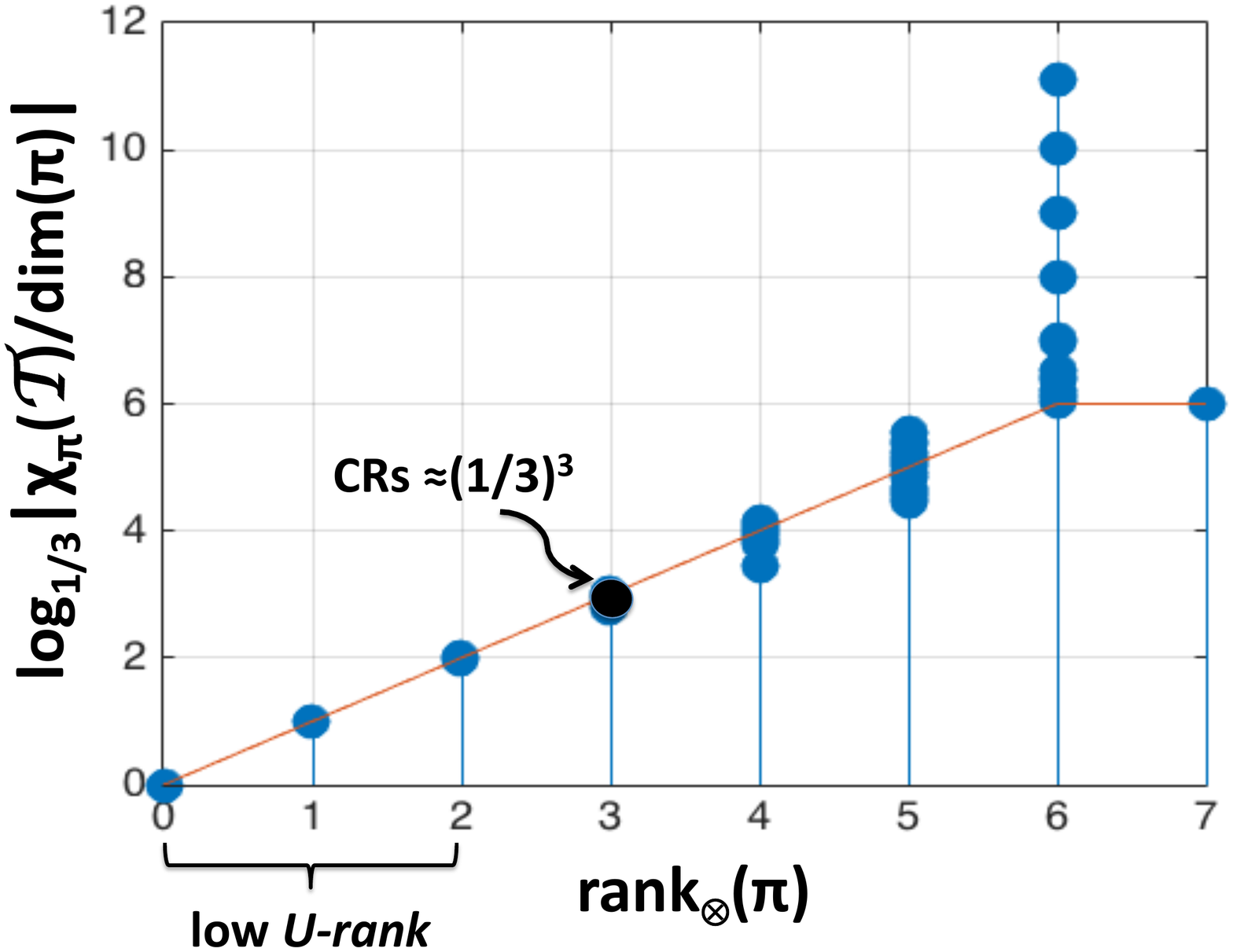}
\caption{$\log _{1/3}$-scale of CRs at transvection $\mathcal{T}$ vs. tensor
rank for irreps $\protect\pi $ of $GL_{7}(\mathbb{F}_{3})$.}
\label{cr-u-vs-rank_t-gl7_3}

\includegraphics[height=2.9in, width=3.9in]
{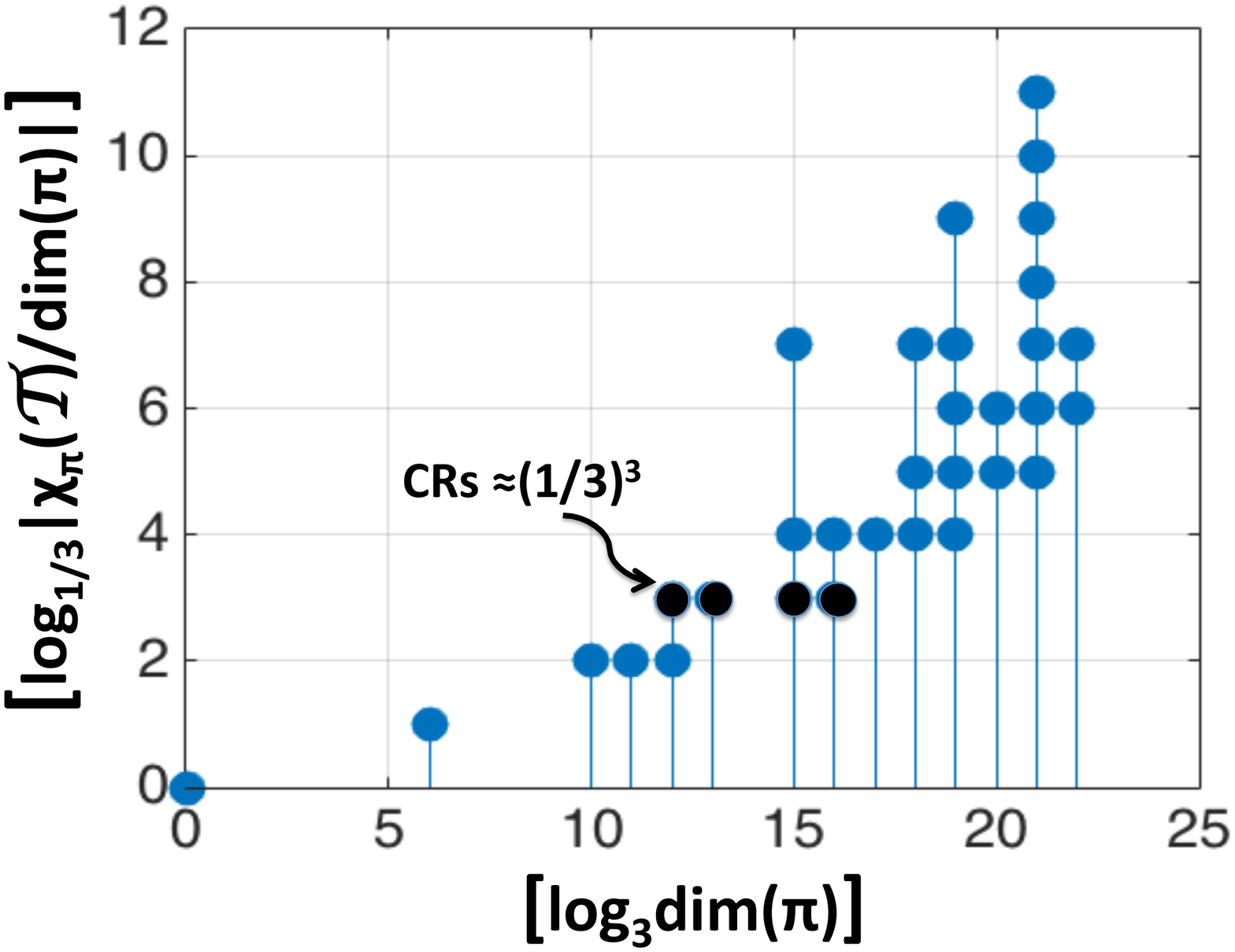}
\caption{$[\log _{1/3}]$-scale of CRs at transvection $\mathcal{T}$ vs. $
[\log _{3}\dim (\protect\pi )]$ for irreps $\protect\pi $ of $GL_{7}(\mathbb{
F}_{3})$.}\label{cr-u-dim-gl7_3}
\end{figure}

Figures \ref{cr-rank_t-sp6_13} and \ref{cr-u-vs-rank_t-gl7_3} illustrate this
numerically. In these figures, the group $G$ is, respectively, $Sp_{6}(
\mathbb{F}_{13})$ and $GL_{7}(\mathbb{F}_{3}),$ the CRs are evaluated on a
transvection $\mathcal{T}$ (i.e., $\mathcal{T}$ is conjugate to a matrix
with $1$'s on the diagonal and only one non-zero off diagonal entry) and $
rank_{\otimes }(\pi )$ denotes the ``tensor rank" of $\pi $ (to be defined later). 

A significant discovery so far is that, although the dimensions of the irreps of a
given rank may vary by quite a lot (they can differ by a large power of $q$), the
CRs of these irreps (for certain elements of interest) are in many cases nearly equal. We will not formally
show this in these notes, but the comparison between Figures 
\ref{cr-u-vs-rank_t-gl7_3} and \ref{cr-u-dim-gl7_3} illustrates\footnote{
We denote by $[x]$ the nearest integer to $x$.} this phenomenon (for more details see [Gurevich-Howe17]).

\vskip .1 in

{\bf 0.2. Earlier developments.}

\vskip .1 in

Before we formulate the main results of these notes, we would like to give
some remarks on the theory of size of representations in a way that places
our work in some historical context.

In the study of representations of algebraic groups over local fields,
considerable attention has been devoted to concepts of size of
representations [Vogan80], [Vogan16], [Wallach82]. For finite groups, since
all irreducible representations are finite dimensional, the dimension of the
space on which a given irreducible representation is realized can serve as a
first proxy for size. However, since most representations of groups over
local fields are infinite dimensional, any notion of size has to appeal to
other considerations. Since many representations are realized on spaces of
functions (or sections of vector bundles) on manifolds or algebraic varieties
(which are often homogeneous spaces for the relevant group) the dimension of
this variety, can provide some measure of size for the representation. For
unipotent groups, the Kirillov orbit method [Kirillov62], [Kirillov04],
[Pukanszky67], assigns a coadjoint orbit in the dual of the Lie algebra to
each representation, and the dimension of this orbit (or half of it, since
these orbits are always even dimensional) can provide an appropriate
measure. Indeed, one of the conclusions of the orbit method for nilpotent
groups is that the representation attached to a given coadjoint orbit can be
realized on sections of a line bundle over a homogeneous space for the
group, of dimension equal to half the dimension of the orbit. Such ideas
extend, in large part, to solvable groups as well [Auslander-Kostant71],
[Bernat-Conze-Duflo-L\'{e}vy-Nahas-Rais-Renouard-Vergne72], [Pukanszky78].

For reductive groups, although the orbit method succeeds in describing a large family of representations, 
including the most ``typical" ones (for example, the ones that appear in the Plancherel formula
[Knapp86], [Vergne83], [Wallach92]), it does not quite cover all
representations, and some of the ones it misses are arguably some of the
most interesting. For real or complex groups $G$, a representation of $G$
gives rise to a ``derived representation" of the Lie algebra 
$\mathfrak{g}$ of $G$, which then further engenders a representation of $U(%
\mathfrak{g})$---the universal enveloping algebra of $\mathfrak{g}$
[Knapp86], [Wallach82]. In this context, the notion of \textit{%
Gelfand-Kirillov dimension} [Knapp86],, [Vogan80], [Vogan16], [Wallach82]
can often be used to assign a number that can be considered a measure of the
size of representations. For representations of groups over non-Archimedean
fields, this concept is not available. Nevertheless, a qualitatively similar
notion of size seems to be relevant for known examples. 

\vskip .1 in

In [Howe82], the notion of \textit{rank} of a representation was introduced
for unitary representations of symplectic groups, and this was extended to
all classical groups in [Li89-1, Li89-2]. The rank of a representation, in
the sense of these papers, is not a very precise measure of size. For
example, representations of a real reductive group can have the same rank in
the sense of [Howe82], but different Gelfand-Kirillov dimensions. However,
it weakly compatible with such notions. Thus, the Gelfand-Kirillov
dimensions of representations of a given rank must vary within a limited
range, and the lower and upper limits of this range increase with the rank.

Unfortunately, the notion of rank discussed in [Howe82] and [Li89-1, Li89-2]
applies only to unitary representations, not to the full class of admissible
representations that form the main preoccupation of the literature of
representation theory for reductive groups. Perhaps because of this
limitation, the idea of rank has found relatively little resonance in that
literature, despite its relevance, as demonstrated in [Howe82],
[Li89-1,Li89-2], for understanding the structure of the unitary dual, as
well as its relationship to ideas developed in the classical theory of
automorphic forms [Gelbart75], [Li92], [Duke-Howe-Li92].

\vskip .1 in

Recently [Gurevich17], [Gurevich-Howe15], [Howe17-1], [Howe17-2], we
realized that there is a way of formulating the notion of rank so that it
makes sense for admissible representations. A major benefit of the new point
of view is that it is relevant for finite fields as well as local fields.
For reasons of ease of exposition, these notes will discuss the new
developments only in the case of finite fields.


\vskip .1 in

{\bf 0.3. Rank of a representation.}

\vskip .1 in

An important invariant attached to any finite group $G$ is its \textit{%
representation ring} [Zelevinsky81] 
$$
R(G)=\mathbb{Z}
\lbrack \widehat{G}], \eqno (0.3.1)
$$
i.e., the ring (aka Grothendieck ring) generated from the set $\widehat{G}$
of isomorphism classes of irreducible representations (irreps) using the
operations of addition and multiplication given, respectively, by direct sum 
$\oplus $ and tensor product $\otimes .$

The main purpose of these notes is to demonstrate that in the case that $G$
is a finite classical group the ring $R(G)$ has a natural rank filtration
that encodes important properties of the representations. In particular,
moving to the associated graded pieces, each member 
$\pi$ of the unitary
dual $\widehat{G}$ gets a well defined non-negative integer that will be its
rank. This integer seems---as illustrated in Figures \ref{cr-rank_t-sp6_13}, %
\ref{cr-u-vs-rank_t-gl7_3}, and \ref{cr-u-dim-gl7_3}, to be intimately
related to certain analytic properties of $\pi$ such as its dimension and size of
character ratios.

\vskip .1 in 

{\bf Remark 0.3.2.}
Another aspect of the new point of view on rank is a connection with the
work of Zelevinsky [Zelevinsky81] on the representations of symmetric groups
and general linear groups over finite fields. Zelevinsky's basic construction, of forming the direct sum over $n$ of the Grothendieck rings of $GL_{n}(\mathbb{F}_{q})$ for a fixed finite field $\mathbb{F}_q$ with $q$ elements, appears naturally [Howe17-1] in connection with the new notion of rank. Also, we have indications for extending Zelevinsky's theory beyond $GL_{n}$.

\vskip .1 in

In fact, in these notes we give several compatible notions of rank (each probably will have its advantages and limitations), and hopefully all of them will serve as valuable analytical tools.

For the sake of introduction, it will be suffice for us to illustrate the main
results and idea using the general linear groups $GL_{n}=GL_{n}(\mathbf{F})$
and the symplectic groups $Sp_{2n}(\mathbf{F})$, where $\mathbf{F}=\mathbb{F}%
_{q}$ is the finite field with $q$ elements, where $q$ (in the $Sp_{2n}$ case) is a power of an odd prime.

\vskip .15 in

{\bf \underline{The $U$-rank}}

\vskip .15 in

The $U$-rank is a notion related to the usual notion of rank of a matrix in linear algebra.

\medskip

{\bf The $GL_{n}$ case:}

\medskip 

Consider the subgroup $U<GL_{n}$ of unipotent upper triangular matrices%
$$
U = \left \{ \left [ \matrix {I_a \ \ \ T \cr 0 \ \ \ I_b} \right ]; \,\, T \in M_{a,b} \right \}, \eqno (0.3.3)
$$
where $a=b=\frac{n}{2}$ if $n$ is even, $a=\frac{n-1}{2}$ and $b=\frac{n+1}{2}$
for $n$ odd, and $I_{a}$ (resp. $I_{b}$) denotes the identity matrix of size 
$a\times a$ (resp. $b\times b$), and $M_{a,b}$ is the space of matrices
of size $a\times b$ over $\mathbf{F}$.

The group $U$ is abelian and can be identified with its group of characters
(aka Pontrjagin dual) $\widehat{U}$, via the map $T\mapsto \chi _{T}$, where 
$\chi _{T}(T^{\prime})=\chi _{0}(trace(T^tT^{\prime }))$, where $T^t$ is the matrix transpose of $T$, and $\chi _{0}$ is
a fixed non-trivial character of $\mathbf{F}$. In particular, every character $\chi \in \widehat{U}$ comes naturally with an integer\footnote{We denotes by $\left\lfloor \,\, \right\rfloor $ the floor function.} $0\leq rank(\chi )\leq \left\lfloor \frac{n}{2}\right\rfloor $, i.e., the rank of the matrix $T$ (0.3.3) associated to it by the identification.

Now, considering a representation $\pi $ of $GL_{n},$ we look at its
restriction $\pi _{|U}$ and define the $U$\textit{-rank }of $\pi $ to be the
integer 
$$
rank_U (\pi) = max \{ rank (\chi); \,\, \chi \, \, \rm{appears \,\, in}  \, \pi_{|U} \}.  \eqno (0.3.4)
$$

In terms of the ring (0.3.1), for $0\leq k\leq \left\lfloor \frac{n}{2}\right\rfloor $ we define 
$$
F_{U}^{k}=F_{U}^{k}(R(GL_{n}))=\{\pi \in R(GL_n); rank_{U}(\pi )\leq k\}.$$
Then the $F_{U}^{k}$ form a filtration of $R(GL_{n})$ that we call the $U$\textit{-rank
filtration}, and of course the members of $F_{U}^{k}\smallsetminus
F_{U}^{k-1}$ are exactly the representations of $U$-rank $k.$

It can be shown (we will come back to more quantitative aspects of this
later) that most of the representations of $GL_{n}$ are of $U$-rank $%
k=\left\lfloor \frac{n}{2}\right\rfloor .$ This motivates us to call the representations of $U$-rank $k<\left\lfloor \frac{n}{2}\right\rfloor ,$ \textit{low }$U$%
\textit{-rank} \textit{representations}.

\medskip 

{\bf The $Sp_{2n}$ case:}

\medskip 

Here, we consider the subgroup $U<Sp_{2n},$ called the \textit{Siegel
unipotent} subgroup, 
$$
U = \left \{ \left [ \matrix {I_n \ \ \ S \cr 0 \ \ \ I_n} \right ]; \,\, S \in S_n^2 \right \}, \eqno (0.3.5)
$$
where $S_{n}^{2}$ denotes the space of $n\times n$ symmetric matrices over $%
\mathbf{F.}$

The group $U\simeq S_{n}^{2}$ is abelian and, exactly as in the $GL_{n}$
case, can be identified with its Pontrjagin dual $\widehat{U}$. This leads
to a well defined notion of rank for every character $\chi \in \widehat{U},$
i.e., the rank of the matrix $S\in S_{n}^{2}$ that corresponds to $\chi,$
in particular, $0\leq rank(\chi )\leq n$. Finally, we have a notion of $U$%
\textit{-rank} for every representation $\pi $ of $Sp_{2n}$ using the
Formula (0.3.4), and again the story of this invariant, as in the
case of $GL_{n},$ can be told in terms of the $U$\textit{-rank filtration}of
the ring $R(Sp_{2n}).$

Again, it can be shown that most of the representations of $Sp_{2n} $ are of maximal $U$-rank $k=n.$ Hence, we will call the representations of $U$-rank $k<n,$ {\it low $U$-rank representations.}

\vskip .15 in

{\bf \underline{The $\otimes$-rank}}

\vskip .15 in

We would like to refine the notion of $U$-rank and filter further the
collection of representations of maximal $U$-rank, i.e., those of $U$-rank $%
k=\left\lfloor \frac{n}{2}\right\rfloor $ in the $GL_{n}$ case, and $k=n$ in
the $Sp_{2n}$ setting. The motivation is to account for the significant
differences (e.g., in dimensions and character ratios) between the members of
that collection. Figures \ref{cr-rank_t-sp6_13} and \ref{cr-u-vs-rank_t-gl7_3} illustrate this. In both pictures the low $U$-rank
representations coincide with these of $\otimes $-rank (to be defined
shortly) $k=0,1,2,$ and a hierarchy in order of magnitude of CRs is clearly
presented well beyond the range of low $U$-rank irreps.

The notion of $\otimes$\textit{-rank} seems (see Figures \ref{cr-rank_t-sp6_13} and \ref{cr-u-vs-rank_t-gl7_3} for illustration) to be an
appropriate refinement. We were led to formulate it by examining the structure of our main construction (the ``eta correspondence") for low $U$-rank irreps.

We will introduce this notion using a specific filtration of the representation ring (0.3.1).

\medskip 

{\bf The $GL_{n}$ case:}

\medskip 

We consider the $U$-rank $k=1$ irreps of $GL_{n}$. They are of smallest
possible dimensions (order of magnitude of $q^{n-1}$) after the one dimensional ones$,$ and can be realized as the non-trivial irreps that appear in the permutation representation $\Pi$ of $GL_{n}$ on $L^{2}(\mathbb{F}_{q}^{n})$, and a twist (i.e., tensor) of these by characters of $GL_{n}$. Now, we introduce the 
$\otimes$\textit{-rank filtration} by letting $F_{\otimes }^{0}=F_{\otimes }^{0}(R(GL_n))$ be the characters of $GL_n$, and for $1 \leq k$ letting $F_{\otimes }^{k}=F_{\otimes }^{k}(R(GL_{n}))$
be the collection of linear combinations of irreps that appear in 
$\ell$-fold, $\ell \leq k$, tensor products of $U$-rank 1 representations.

It is not difficult to show [Howe73-2] that:

\vskip .05 in

{\bf Claim 0.3.6.} The first $k$ such that $F_{\otimes }^{k}=\widehat{GL}_{n}$, is $k=n.$

\vskip .05 in 

Finally, we say that a representation $\pi $ of $GL_{n}$ is of {\it tensor rank} 
$k$, and denote $rank_{\otimes}(\pi )=k,$ if $\pi \in F_{\otimes }^{k}\smallsetminus F_{\otimes }^{k-1}.$

\vskip .15 in 

{\bf The $Sp_{2n}$ case:}

\medskip

Again, we start with the $U$-rank $k=1$ irreps of $Sp_{2n}$. These are the
irreps of smallest possible dimensions $\neq 1.$ There are four of them, with
two of dimension $\frac{q^{n}-1}{2}$, and two of dimension $\frac{q^{n}+1}{2}$.
They can be realized as the components of the two oscillator (aka Weil)
representations of $Sp_{2n}$ [G\'{e}rardin77], [Howe73-1], [Weil64]. Next,
we introduce the $\otimes $\textit{-rank filtration} $F_{\otimes
}^{k}=F_{\otimes }^{k}(R(Sp_{2n}))$ by letting $F_{\otimes }^{k}$ be the
collection of integer linear combinations of $\ell$-fold, $\ell \leq k$,
tensor products of $U$-rank 1 representations.

It can be shown [Adams-Moy93], [Howe73-2], that:

\vskip .05 in

{\bf Claim 0.3.7.} The first $k$ such that $F_{\otimes }^{k}=\widehat{Sp}_{2n}$ is $k=2n.$

\vskip .05 in

Finally, we say that a representation $\pi$ of $Sp_{2n}$ is of 
\textit{tensor rank} $k$, and denote $rank_{\otimes }(\pi )=k$, if 
$\pi \in F_{\otimes }^{k}\smallsetminus F_{\otimes }^{k-1}$.

\vskip .15 in

{\bf \underline{Comparison}}

\vskip .15 in

Tensor rank is defined, in both cases of $GL_{n}$ and $Sp_{2n},$ using $U$-rank $k=1$ irreps. In particular, a simple compatibility property
between $U$-rank and tensor product implies that $rank_{U}\leq rank_{\otimes
}$ for all representations.

Based on theoretical considerations, and numerical evidences, we conjecture
that

\vskip .1 in

{\bf Conjecture 0.3.8 (Agreement of ranks).} 
On the collection of low $U$-rank representations of $GL_{n}$ and $Sp_{2n},$
we have $rank_{U}=rank_{\otimes }.$

\vskip .15 in

{\bf 0.4. Information on representations of a given rank.}

\vskip .1 in

How to get information on the family of representations of a given $U$-rank or $\otimes$-rank? For example information on the character ratios of irreps that are needed for the harmonic analysis problems that we mentioned in Subsection 0.1.

We propose to use the eta correspondence as a method the leads to an answer for the above question. Indeed, it seems [Gurevich-Howe15], [Gurevich-Howe17],[Gurevich17] that this construction might provide the required analytic information about the representations of finite classical groups. Moreover, it does so in a way that for certain purposes complements the information supplied by Harish-Chandra's ``philosophy of cusp forms" [Harish-Chandra70, Howe-Moy86] or more generally by the classification given by Lusztig in [Lusztig84] for the irreps of these groups.

The eta correspondence is a collection of correspondences of representations associated to dual pairs [Gelbart77], [Howe73-2], [Howe77],in the symplectic groups. It is closely related to, but distinct from, the widely used theta correspondence for classical groups over local fields [Gan-Takeda16], [Gelbart75], [Howe77], [Kudla86], [Moeglin-Vigneras-Waldspurger87], [Prasad93].

For classical groups over finite fields, the theta correspondence does not exist and the eta correspondence provides a substitution. In the past, several authors [Adams-Moy93], [Aubert-Michel-Rouquier96], [Srinivasan79] and in particular [Howe73-2], studied certain representation theoretic properties associated with dual pairs in the finite field setting, but no explicit correspondence was exhibited, nor was the concept of rank introduced in this context.

\vskip .1 in

{\bf \underline{The eta correspondence}}

\vskip .1 in

We describe the eta correspondence and how it gives a systematic construction of a large collection (conjecturally all) of the irreps of $GL_{n}$ and $Sp_{2n}$ of a given $U$-rank or $\otimes$-rank $k$.

\vskip .1 in

{\bf The $GL_{n}$ case:}

\vskip .1 in

We start the construction by fixing $n$ and $k$, with $n\geq k,$ and
consider the vector spaces $\mathbf{F}^{n}$ and $\mathbf{F}^{k}$ over $
\mathbf{F}=\mathbb{F}_q.$ We then form the $2nk$-dimensional 
vector space $(\mathbf{F}^{n}\otimes \mathbf{F}^{k})\oplus (\mathbf{F}^{n}\otimes \mathbf{F}^{k})^{\ast }$. This is a direct sum of a space and its dual, so it is in a natural way a symplectic vector space. We denote the associated symplectic group by $Sp_{2nk}$ and note that $GL_{n}$ and $GL_{k}$ can be regarded naturally as a pair of commuting subgroups in $Sp_{2nk}.$ In fact, in this way $GL_{k}$ is the centralizer of $GL_{n}$ in $Sp_{2kn}$ and vice versa,
namely $(GL_{n},GL_{k})$ is a \textit{dual pair} of subgroups in $Sp_{2nk},$ in the sense of [Howe73-2].

The next step of the construction is to take an oscillator representation [G\'{e}rardin77], [Howe73-1], [Weil64] $\omega _{2nk}$ of $%
Sp_{2nk},$ restrict it to $GL_{n}\times GL_{k}$ and look at the
decomposition into isotypic components of irreps of $GL_{k}$ 
$$
\omega _{2nk|GL_{n}\times GL_{k}}\simeq \sum_{\nu \in \widehat{GL}_{k}}\Xi
_{\nu }\otimes \nu,  \eqno (0.4.1)
$$
where the space $\Xi _{\nu }$ is a representation of $GL_{n}.$

We want to say something about the spaces $\Xi_{\nu }$ in (0.4.1). Some of them might be zero and some, but ``not too many", are reducible. Indeed, as was shown in [Howe73-2] the group
algebras
$$
\mathcal{A}_{GL_n} = \mathbb{C} [ \omega _{2nk}(GL_n) ] \,\,\, \rm{and} \,\,\, \mathcal{A}_{GL_k}= \mathbb{C} [ \omega _{2nk}(GL_{k})], \eqno (0.4.2)
$$
are ``approximately" each other's commutant in $End(\omega _{2nk})$ in the
following sense. Denote by $\widehat{GL}_{n}(\omega _{2nk})$ the subset of $\widehat{GL}_{n}$ of irreps that appear in $\omega _{2nk}$. It is not difficult to show that for every $1 \leq k \leq n$, we have a proper containment
$$
\widehat{GL}_{n}(\omega _{2n(k-1)})\subsetneq \widehat{GL}_{n}(\omega _{2nk}), \eqno (0.4.3)
$$ 
and that the sequence stabilizes exactly at the $k=n$ step, i.e., $\widehat{GL}_{n}(\omega _{2n^2})=\widehat{GL}_{n}$. We call the difference between these two sets, i.e., 
$$\widehat{GL}_{n}(\omega_{2nk})\smallsetminus \widehat{GL}_{n}(\omega _{2n(k-1)})
$$
the {\it new} $k$\textit{-spectrum} of $GL_{n},$ and denote it $(\widehat{GL}_{n})_{new,k}$. Finally, we denote by $\left( \omega _{2nk}\right) _{GL_{n}}^{new}\subset
\omega _{2nk}$ the subspace spanned by all isotypic components of members
from $(\widehat{GL}_{n})_{new,k}$. 

The algebras $\mathcal{A}_{GL_{n}}$ and $
\mathcal{A}_{GL_{k}}$ act on $\left( \omega _{2nk}\right) _{GL_{n}}^{new}$
and now we can remove the above ``approximately" phrase, i.e., we have
[Howe73-2]:

\vskip .15 in

{\bf Theorem 0.4.4 (The $\eta$-correspondence for $GL_n$ - non explicit form).}
The algebras $\mathcal{A}_{GL_n}$ and $\mathcal{A}_{GL_k}$ are each
other's commutant acting on $(\omega _{2nk})_{GL_n}^{new}$.
In particular (by the double commutant theorem [Weyl46]), we have a bijection
$$
\nu \longmapsto \eta (\nu )<\Xi _{\nu },  \eqno (0.4.5)
$$
from a subset of $\widehat{GL}_{k}$ onto $(\widehat{GL}_{n})_{new,k}$.

\vskip .15 in

We call (0.4.5) the \textit{eta correspondence}.

\vskip .07 in

It is a good news that, the bijection (0.4.5), and in particular its domain of definition and image, can be made explicit and in a way that leads to various explicit formulas for the $U$-rank and 
$\otimes$-rank $k$ irreps of $GL_n$.      

We start with the case $k<\left\lfloor \frac{n}{2}\right\rfloor$. Here all irreps of $GL_k$ participate in the domain of (0.4.5) while the image of the correspondence consists only of irreps of $U$-rank $k$. To state this precisely, let us introduce some notations.

Denote by $(\widehat{GL}_n)_{U,k}$ the {\it $U$-rank $k$-spectrum} of $GL_n$ consisting of all $U$-rank $k$ irreps, and by $(\omega _{2nk})_{GL_n}^{U,k} \subset \omega _{2nk}$ the subspace spanned by all isotypic components of members from $(\widehat{GL}_{n})_{U,k}$. 

We have,
\vskip .1 in

{\bf Theorem 0.4.6 (Construction of $U$-rank $k$ irreps).} Assume $k<\left\lfloor \frac{n}{2}\right\rfloor$. Then,

\begin{enumerate}
\item The new $k$-spectrum and $U$-rank $k$-spectrum of $\omega _{2nk}$ coincide, i.e., 
$(\omega _{2nk})_{GL_n}^{new} = (\omega _{2nk})_{GL_n}^{U,k}$. 

\item The domain of (0.4.5) is all of $\widehat{GL}_k$. In particular, for every $\nu \in \widehat{GL}_k$ there is a unique $U$-rank $k$ irrep $\eta(\nu)$ in $\Xi _{\nu }$, it appears with multiplicity one, and the other components are of lower $U$-rank.

\end{enumerate}

\vskip .1 in

So, for $k<\left\lfloor \frac{n}{2}\right\rfloor$ the eta correspondence (0.4.5) is an injection from $\widehat{GL}_{k}$ to 
$(\widehat{GL}_{n})_{U,k}$. We have reasons, supported by theoretical and
numerical evidences, to believe that,

\vskip .1 in

{\bf Conjecture 0.4.7 (Exhaustion).} 
Assume $k<\left\lfloor \frac{n}{2}\right\rfloor$. Then, every irrep of $GL_n$ of $U$-rank $k$ appears, maybe after twisting by a character of $GL_n$, in the image of the eta correspondence (0.4.5).

\vskip .1 in

In summary, Theorem 0.4.6 tells us that the for $k<\left\lfloor \frac{n}{2}\right\rfloor$ the new $k$-spectrum can be described using the $U$-rank $k$-spectrum, and (using Conjecture 0.4.7) vice versa. Moreover, it tels us that, in this range, both spectrums can be described by the eta correspondence, and vice versa. 

\vskip .05 in

What about (0.4.5) and the various $k$-spectrums in the range $\left\lfloor \frac{n}{2}\right\rfloor \leq k$?  

\vskip .05 in

Of course, in the above domain the notion of $U$-rank does not seems to be directly effective anymore. But, it follows from the definitions that, the tensor rank of an irrep $\rho$ of $GL_n$ is the minimal $k$ such that $\chi \otimes \rho$ belongs to the new $k$-spectrum $(\widehat{GL}_{n})_{new,k}$, for some character of $GL_n$. So what we want to have, is a good formula for the irreps in the new $k$-spectrum, and an understanding of what irreps of $GL_k$ participate in (0.4.5).

In these notes we obtain (see [Gurevich-Howe17] for more details) the required formulas using Harish-Chandra's ``philosophy of cusp forms" (P-of-CF for short) [Harish-Chandra70, Howe-Moy86] that was put forward by him in the 60s. Moreover, we obtain the following extension of Theorem 0.4.6. Denote by 
$(\widehat{GL}_k)_{new,\geq 2k-n}$ of all irreps of $GL_k$ that appears in the oscillator representation $\omega _{2nl}$ for some $l \geq 2k-n$ and not before.  

\vskip .1 in

{\bf Theorem 0.4.8 (The $\eta$-correspondence for $GL_n$ - explicit form).} Fix $0 \leq k \leq n$. Then, the domain of the eta correspondence (0.4.5) is $(\widehat{GL}_k)_{new,\geq 2k-n}$. In particular, for every $\nu$ in $(\widehat{GL}_k)_{new,\geq 2k-n}$ there is a unique irrep $\eta(\nu)$  in $\Xi _{\nu }$ (0.4.1) which belongs to $(\widehat{GL}_n)_{new,k}$, it appears with multiplicity one, and the other components belong to the lower new spectrum.

\vskip .1 in

In conclusion, we give in these notes a more or less explicit formula for the $\eta$ correspondence (0.4.5), and how it give an explicit identification of the individual members of the set 
$(\widehat{GL}_n)_{\otimes,k}$ of irreps of $GL_n$ of tensor rank $k$.

\vskip .1 in

{\bf The $Sp_{2n}$ case:}

\vskip .1 in

There are some similarities and some differences between the 
$\eta$-correspondence construction in the $Sp_{2n}$ and $GL_n$ cases.

Also in the $Sp_{2n}$ case we start by introducing the dual pairs involved. For the sake of introduction it is enough for us to know [Lam73] that for every $k\neq 0$, up to isometry there are only two non-degenerate symmetric bilinear forms on $\mathbb{F}_q^k$. For now, for simplicity, let us denote either of the two associated orthogonal groups by $O_k$ and 
formulate statements that, unless stated otherwise, work for both of them. 

There is a natural construction so that $(Sp_{2n},O_{k})$ becomes a dual pair in a symplectic group $Sp_{2nk}$. In particular, we can consider the restriction of an oscillator representation $\omega _{2nk}$ of $Sp_{2nk}$ to $Sp_{2n}\times O_k$ and look at the isotypic decomposition with respect to irreps of $O_k$
$$ 
\omega _{2nk|Sp_{2n} \times O_k}\simeq \sum_{\tau \in 
\widehat{O}_{k}}\Omega _{\tau}\otimes \tau . \eqno (0.4.9)
$$
We proceed to use (0.4.9) in order to describe the collection $(%
\widehat{Sp}_{2n})_{U,k}$ of low $U$-rank $k$ irreps of $Sp_{2n}$.

Consider the commuting subalgebras, 
$$
\mathcal{A}_{Sp_{2n}}, \, \mathcal{A}_{O_k}\subset End(\omega _{2nk}),
$$
defined in the same way as those in (0.4.2). Denote by 
$(\omega_{2nk})_{Sp_{2n}}^{U,k}\subset \omega _{2nk}$ the {\it $U$-rank $k$-spectrum}, i.e., the sum of all isotypic components in $\omega_{2nk}$ of members from 
$(\widehat{Sp}_{2n})_{U,k}$. 
In [Gurevich-Howe15] we proved,

\vskip .1 in

{\bf Theorem 0.4.10 (Construction of low $U$-rank irreps).} Suppose $k<n.$ 
The following hold,

\begin{enumerate}
\item For every $\tau \in \widehat{O}_k$ the representation $\Omega _{\tau}$ 
(0.4.9) is of $U$-rank $k$.

\item The centralizer of $\mathcal{A}_{Sp_{2n}}$ acting on 
$(\omega_{2nk})_{Sp_{2n}}^{U,k}$ is $\mathcal{A}_{O_k}$ and vice versa.

\item For each $\tau$ in $\widehat{O}_k$, there is a unique $U$-rank $k$ $Sp_{2n}$ irrep $\eta(\tau)$ in $\Omega_{\tau}$, it appears with multiplicity one, and all other $Sp_{2n}$ constituents of $\Omega_{\tau}$ are of $U$-rank less than $k$. The correspondence
$$
\tau \longmapsto \eta(\tau) \eqno (0.4.11)
$$
is an injection from $\widehat{O}_k$ into $(\widehat{Sp}_{2n})_{U,k}$.
\end{enumerate}

We believe, with theoretical and numerical supports, that the construction of low $U$-rank irreps announced in Theorem 0.4.10 is exhaustive in the following sense. Denote by $O_{k+}$ and $O_{k-}$ the orthogonal groups of the two possible inner product structures on ${\bf F}_q^k$ [Lam73]. It can be shown [Gurevich-Howe15] that the images $\eta(\widehat{O}_{k+})$ and 
$\eta(\widehat{O}_{k-})$ under the mapping (0.4.11) are disjoint, for any $k$, $1 \leq k \leq n$.

We formulate the following conjecture:

\vskip .1 in 

{\bf Conjecture 0.4.12 (Exhaustion).}
For $0\neq k<n$, using both orthogonal groups we obtain all 
low $U$-rank irreps, i.e.,
$$
\eta (\widehat{O}_{k+})\sqcup \eta (\widehat{O}_{k-})=(\widehat{Sp}_{2n})_{U,k}.
$$

\vskip .1 in 

Next, we would like to describe the members of the set 
$(\widehat{Sp}_{2n})_{\otimes ,k}$ of irreps of $\otimes $-rank 
$k$ of $Sp_{2n}.$ In the domain of {\it low $\otimes$-rank} irreps, i.e., 
$\otimes$-rank $k<n$, we have a good description. Namely, it follows from the definition
of $(\widehat{Sp}_{2n})_{\otimes ,k}$ and the $\eta$-correspondence (0.4.11) that,

\vskip .1 in

{\bf Theorem 0.4.13 (Construction of low $\otimes$-rank irreps).}
For $0\neq k<n$, using both orthogonal groups, the $\eta $-correspondences (0.4.11) produce all irreps of $\otimes$-rank $k$, i.e., 
$$
\eta (\widehat{O}_{k+})\sqcup \eta (\widehat{O}_{k-})=(\widehat{Sp}_{2n})_{\otimes ,k}.
$$

\vskip .1 in

We are left to describe the irreps of $\otimes$-rank $k$ for $k \geq n$.
Here our knowledge is still not complete and further research is required to
have the full story as we have for $GL_{n}.$ Let us explain what we do have
and what we still lack.

We do have the extension of the $\eta$-correspondence to construct all the
irreps of a given tensor rank also for $k\geq n$. Indeed, recall that by a
theorem of Witt [Lam05] the collection of orthogonal group over a finite
field is naturally organized in four Witt towers:

\begin{itemize}
\item $O_{m,m},\,\, m\geq 0,\,\,$ - the isometry groups of the direct sums of $m$ copies of the \textit{hyperbolic plane} $P$. This collection is known as the {\it split Witt tower}.

\item $O_{m+2,m},\,\, m\geq 0,\,\,$ - the isometry groups of the direct sums of the
anisotropic plane and $m$ copies of $P$.

\item $O_{m+1,m},\,\, m\geq 0,\,\,$ - the isometry groups of the direct sums of a line equipped with quadratic form of discriminant $1$ and $m$ copies of $P$.

\item $O_{m,m+1},\,\, m\geq 0,\,\,$ - the isometry groups of the direct sum of a line equipped with quadratic form of discriminant $-1$ and $m$ copies of $P$.
\end{itemize}

Now we can describe the eta correspondence in each tower using the notion of
new spectrum as we did in the case of $GL_{n}$ (where we had just one Witt
tower, namely $GL_{k},$ $k\geq 0$). Let us tell how this is done by
examining one tower, that of $O_{m,m}$; the eta correspondence for other towers can be dealt with in a similar manner.

Consider the dual pair $(Sp_{2n},O_{m,m})$ in $Sp_{4nm}$ and an oscillator representation $\omega _{4nm}$ of this group. 
Denote by $\widehat{Sp}_{2n}(\omega _{4nm})$ the subset of $\widehat{Sp}_{2n}$ 
of irreps that appear in $\omega _{4nm}$. It is not difficult to show that 
$\widehat{Sp}_{2n}(\omega _{4n(m-1)})\subset \widehat{Sp}_{2n}(\omega_{4nm})$. 
We call the set 
$\widehat{Sp}_{2n}(\omega_{4nm})\smallsetminus \widehat{Sp}_{2n}(\omega _{4n(m-1)})$ the 
{\it $2m$-new spectrum} and denote it $(\widehat{Sp}_{2n})_{new,2m}$. 
Finally, we denote by $\left( \omega _{4nm}\right)_{Sp_{2n}}^{new}
\subset \omega _{4nm}$ the subspace spanned by all
isotypic components of members from $(\widehat{Sp}_{2n})_{new,2m}$. The
algebras $\mathcal{A}_{Sp_{2n}}$ and $\mathcal{A}_{O_{m,m}}$ act on $\left(
\omega _{4nm}\right) _{Sp_{2n}}^{new}$ and we have [Howe73-2]:

\vskip .15 in

{\bf Theorem 0.4.14 (The $\eta$-correspondence in the split Witt tower).}

The algebras $\mathcal{A}_{Sp_{2n}}$ and $\mathcal{A}_{O_{m,m}}$ are each other's commutant acting on $(\omega_{4nm})_{Sp_{2n}}^{new}$. In particular, we have a bijection,  
$$
\tau \longmapsto \eta (\tau )<\Omega _{\tau },  \eqno (0.4.15)
$$
from a subset of $\widehat{O}_{m,m}$ onto $(\widehat{Sp}_{2n})_{new,2m}$. The representation, $\eta (\tau )$ is the unique member of the $2m$-new spectrum in $\Omega _{\tau }$, and it appeara there with multiplicity one.

\vskip .15 in

We call (0.4.15) the {\it eta correspondence for the tower $O_{m,m}$}. It coincides with (0.4.11) in case $k=2m < n$.

\vskip .07 in

Finally, we can use the eta correspondences in the four Witt towers to realize the irreps of a given tensor rank.

It can be shown [Howe73-2] that every irrep of $Sp_{2n}$ appears in the image of each of the eta correspondences at some step. The first time it appears in a given tower will be called {\it 1st occurrence}. 

\vskip .1 in

{\bf Proposition 0.4.16.}
A representation $\pi \in \widehat{Sp}_{2n}$ is of tensor rank $k$ iff its minimal 1st occurrence among all four Witt towers is $k$.

\vskip .1 in

In summary, Theorem 0.4.14 tells us that we have a construction of all irreps of low tensor rank. On the other hand, looking on Theorem 0.4.14 and Proposition 0.4.16 we learn that we don't know what irreps of the relevant orthogonal groups contribute for the associated tensor rank $k$-spectrum in case $k>n$. The missing information is a subject for future research.

\vskip .15 in

{\bf 0.5. The body of these notes.}

\vskip .1 in

The body of these notes is organized as follows.

In Section 1 we recall parts of the theory of the classical (i.e., general
linear, symplectic, orthogonal, and unitary) groups and their Witt groups.

Section 2 discusses vector spaces of matrices and their Pontrjagin duals. In
particular we recall some basic properties of symmetric, skew-symmetric and
Hermitian matrices.

In Section 3 we look at matrix subgroups (i.e., vector spaces of matrices)
in the specific cases of classical groups.

Section 4 uses the language and notions of linear algebra and harmonic
analysis introduced in Sections 2 and 3, to define the notion of $U$-rank for
representations of the classical groups. Some basic properties of $U$-rank
are discussed as well.

Section 5 is devoted to the construction of the representations of low
$U$-rank of the classical groups. This section is important for the development
of the theory for all classical groups. The construction is an outcome of the 
theories of the Heisenberg and oscillator representations, which we introduce in 
Subsection 5.1. After we have these representations at hand, we move in Subsection
5.2 to discuss the notion of dual pair in the symplectic group $Sp_{2n}$,
and give the important representative examples: The $Sp_{2n}$-$O_{k}$ and $
GL_{n}$-$GL_{k}$ dual pairs. We don't discuss the $U_{n}$-$U_{k}$ dual pair
in these notes, because their theory is very similar to the that of the pair 
$Sp_{2n}$-$O_{k}.$ Finally, in Subsection 5.3 we use the restriction of the
oscillator representation to dual pairs to obtain (see Theorem 5.3.1) the
eta correspondence that, in particular, describes irreps of low $U$-rank of
classical groups. We conjecture (see Conjecture 5.3.3) that our construction
his exhaustive.

In Section 6, we give two extensions of the notion of $U$-rank. The first
one, discussed in Subsection 6.2, comes from considering restrictions of
representations of large classical groups to smaller ``block subgroups" of the same type. We will call this \textit{asymptotic rank}. The second one, presented in 
Subsection 6.4, is the notion of tensor rank discussed also earlier in the introduction. In addition, in Subsection 6.3 we discuss how all these notions of rank are related to filtrations of the associated Grothendieck rings of representations. Finally, in Subsection 6.5 we present a comparison between these three notions of rank, and we formulate the ``agreement conjecture" (see Conjecture 6.5.2, which is in fact
implied by the Exhaustion conjecture 5.3.3) that for low $U$-rank
representations, all these notions of rank agree.

Section 7 gives an explicit construction of the {\it oscillator semigroup}
[Howe73-2]. This structure (motivated by the desire to understand the failures
of dual pairs to generate each other's commutant in the oscillator
representation) is a natural extension of the symplectic group $Sp_{2n}$ to a semigroup 
(that, in an algebro-geometric sense, is a ``compactification" of $Sp_{2n}$). It helps 
in the analysis of the ``new spectrum" discussed earlier in the introduction, namely the extension of the eta correspondence in order to describe higher tensor rank irreps.

In Section 8 the above applications of the oscillator semigroup are given.
In particular, Theorem 8.1.1 describe the commuting algebra of one member
in a dual pair, and Theorem 8.2.1 gives the application for the eta correspondence
(we state this for the symplectic-orthogonal dual pair, but this works in
general [Howe73-2]). 

Section 9 is devoted to an explicit description (see Theorem 9.2.3) of the eta correspondence for $GL_{n}$, and in particular how it gives for this group the explicit description of all irreps of a given tensor rank.

In Section 10 we sketch a proof of Theorem 9.2.3 (complete details appear in [Gurevich-Howe17]). Our argument uses a fairly precise descriptions (in Section 10.2 and 10.3) of the tensor rank of representations, and (in Section 10.4) of the $\eta$ correspondence. For this we invoke Harish-Chandra's ``philosophy of cusp forms" (that we recall in Section 10.1). The corresponding formulas (in Sections 10.2-10.4) are of independent interest. 

Section 11 discusses the conjectural agreement, in the case of $GL_n$, between all  notions of
rank, namely, for low $U$-rank representations we have $U$-rank = asymptotic rank = tensor rank. 
We managed to verify the second and third equalities, and in a more restrictive way also the first one.  

\vskip .15 in

{\bf Acknowledgment}. We thank the organizers of 19th Takagi lectures (RIMS - Kyoto University, Japan, July 8-9, 2017), in particular T. Kobayashi, for the invitation to present this work at the event. 

We thank J. Bernstein for sharing some of his thoughts concerning the organization of representations by small ones. 

We thank J. Cannon of the University of Sydney and S. Goldstein of the University of Wisconsin at Madison for help with the numerical computations using the computer system Magma. 

Parts of these notes were written while we were participating in the program ``Representation Theory of Reductive Groups Over Local Fields and Applications to Automorphic forms", May 3 - June 30, 2017, and visitors at the Weizmann Institute - Israel, and we thank the organizers, especially D. Gourevitch and the math department chair O. Sarig, for the invitation.

\vskip .2 in 

{\bf 1. Classical groups and their Witt semigroups}

\vskip .15 in

We start with some background on classical groups.

\vskip .1 in 

{\bf 1.1. Classical groups over algebraically closed fields.}
\vskip .1 in
The term ``classical groups" refers to a family of algebraic groups. In particular, there
is a family of classical groups attached to each field $\bf F$. (Here we are going to exclude
fields of characteristic 2 from the discussion). If $\bf F$ is algebraically
closed, then its family of classical groups has three series of members: the general linear
groups $GL_n = GL_n (\bf F)$, the orthogonal groups $O_n = O_n (\bf F)$, and the
symplectic groups $Sp_{2n} = Sp_{2n}(\bf F)$. 

The group $GL_n(\bf F)$ is the full group of invertible $n \times n$ matrices with coefficients in $\bf F$. It is isomorphic to the group of linear automorphisms of any vector space over $\bf F$ of dimension $n$.

The group $O_n(\bf F)$ is the subgroup of $GL_n(\bf F)$ consisting of matrices 
$R$ such that $R^tR = I_n$, where $R^t$ denotes the transpose of the matrix $R$, and
$I_n$ is the $n \times n$ identity matrix. Equivalently, it is the subgroup of $GL_n(\bf F)$ 
that preserves the quadratic polynomial 
\vskip .05 in
\centerline {$r^2 = \sum_{j = 1}^n x_j^2,$}
\vskip .05 in
\noindent which, if the base field were $\bf R$, the real numbers, we would call
 the ``standard Euclidean
distance squared". It is a standard fact of linear algebra [Lam73], that if $V$ is any vector space 
over
$\bf F$ of dimension $n$, and $B( \ , \ ) = B$ is a non-degenerate, symmetric bilinear form 
(aka {\it inner product})
on $V$, then the isometry group $O_B$ of $B$, that is, the group of linear automorphisms $R$ of $V$ such that
$$
B(R(\vec v), R(\vec v') = B(\vec v, \vec v')),
$$
for any pair of vectors $\vec v$ and $\vec v'$ in $V$, then $O_B$ is isomorphic to $O_n$.
The group $Sp_{2n}(\bf F)$ is the subgroup of $GL_{2n}(\bf F)$  
that preserves the standard skew symmetric, non-degenerate bilinear form 
(aka {\it symplectic} form)
$$
<\vec w, \vec w'> = \sum_{j = 1}^n x_j y'_j - y_j x'_j, 
$$
where
$$
\vec w = \left [ \matrix {\vec x \cr \vec y} \right ] 
= \left [ \matrix { x_1 \cr x_2 \cr x_2 \cr . \cr . \cr . \cr x_n \cr y_1 \cr y_2 \cr . \cr . \cr . \cr y_n}
\right ], 
$$
is the $2n$-dimensional vector whose first $n$ coordinates are given by a vector
 $\vec x$ in ${\bf F}^n$,
and whose last $n$ coordinates are given by another vector $\vec y$ in ${\bf F}^n$. The vector
$\vec w'$ is defined similarly. The group $Sp_{2n}(\bf F)$ consists of elements $S$ in
$GL_{2n}(\bf F)$ such that 
$$
<S(\vec w), S(\vec w')> = <\vec w, \vec w'>.
$$
It is another standard fact of linear algebra, that the isometry group
 of any symplectic form on a $2n$ dimensional vector space over ${\bf F}$ is isomorphic
 to $Sp_{2n}(\bf F)$.

\vskip .1 in

{\bf 1.2. Rationality issues; split forms and Witt semigroups.}

\vskip .1 in

 If the base field $\bf F$ is not algebraically closed, then a {\it classical group over ${\bf F}$}
  is an algebraic matrix group $G$ acting on a vector space $V$ over ${\bf F}$ such that there
is a field extension $\tilde {\bf F}$ such that the extension of scalars of $G$ acting on
$V\otimes \tilde {\bf F}$ is one of the three types of groups $GL_n(\tilde {\bf F})$,
$O_n(\tilde {\bf F})$, or $Sp_{2n}(\tilde {\bf F})$ described above. The detailed list of classical
groups over $\bf F$ depends on the arithmetic of $\bf F$. For $\bf F  = \bf R$, the real
numbers, the list of classical groups is given in [Weyl46], among other places. For $p$-adic fields, 
the list is similar, except the possible collection of inner products is somewhat more complicated.
\vskip .05 in

\noindent {\bf Orthogonal groups}

\vskip .05 in

When the base field is not algebraically closed, it is still the case that there is only one 
isomorphism class of symplectic form, but there can be many isomorphism classes of 
inner products. The isometry group of any one of these will still qualify as an
 orthogonal group, but the isometry groups of non-equivalent inner products may or may not be
isomorphic. 

We will review the structure of the collection of inner products up to isomorphism. In this 
discussion, each inner product is implicitly defined on some vector space, but the role of
the vector space is fairly passive. Its main salient feature is its dimension, which we will refer
to as the {\it rank} of the inner product under consideration. Somewhat more generally, if we 
are given a symmetric bilinear form $B$ on a vector space $V$ over ${\bf F}$, then we consider
 the {\it radical} $R_B$ of $B$, which is the subspace of vectors $\vec u$ in $V$ that are
$B$-orthogonal to all vectors in $V$, meaning that $B(\vec u, \vec v) = 0$ for all vectors
$\vec v$ in $V$. Then it is easy to show that $B$ factors to an inner product on the 
quotient space $V/R_B$, and the {\it rank} of $B$ is the dimension of $V/R_B$.

There is a natural notion of (orthogonal) direct sum of two inner products. Given inner products $B_1$ and $B_2$, 
defined on vector spaces $V_1$ and $V_2$ over $\bf F$, 
one defines an inner product $B_1 \oplus B_2$ on $V_1 \oplus V_2$ by the rule
$$
(B_1 \oplus B_2) \left ( \left [ \matrix {\vec v_1 \cr \vec v_2} \right ], 
\left [ \matrix {\vec v'_1 \cr \vec v'_2} \right ] \right ) = 
B_1 (\vec v_1, \vec v'_1) + B_2(\vec v_2, \vec v'_2),
$$
with $\vec v_j$ and $\vec v'_j$ in $V_j$.
It is clear from this definition that, if $B'_1$ and $B'_2$ are isomorphic to $B_1$ and
$B_2$ respectively, then $B'_1 \oplus B'_2$ will be isomorphic to $B_1 \oplus B_2$. 
It follows that the collection of isomorphism classes of inner products is a semigroup (the
{\it Witt semigroup}) under
direct sum. Witt's Theorem [Lam73] implies that it is a semigroup {\it with cancellation}. That is,
if $B_1 \oplus B_2 \simeq B_1 \oplus B'_2$, then in fact $B_2 \simeq B'_2$.
Also, the Gram-Schmidt orthogonalization process [Strang76] shows that every inner product 
can be realized as a direct sum of rank one inner products. Moreover, it is
not difficult to see [Lam73] that the collection of rank one inner products is naturally 
parametrized by the group ${\bf F}^{\times}/{\bf F}^{\times 2}$, the multiplicative group of $F$, 
modulo perfect squares.

Among the inner products of dimension 2, there is a distinguished one $B_{1,1}$ 
 known as the {\it hyperbolic plane}. It has a basis $\{ \vec e_1, \vec e_2\}$ such that
$$
B_{1,1} (\vec e_1, \vec e_1) = 0 = B_{1,1} (\vec e_2, \vec e_2) \,\,\,\, {\rm and} \,\,\,\,  B_{1,1}(\vec e_1, \vec e_2) = 1.
$$
Thus, the vectors $e_j$ are {\it isotropic}, meaning that they have self inner product 
equal to zero. It is well known [Lam73] and not hard to show that a two-dimensional inner product 
which has isotropic vectors is isomorphic to a hyperbolic plane. 

A slightly more refined, but also fairly easy, fact [Lam73] is that if an inner product $B$ on a 
vector space $V$ allows an isotropic vector, then in fact there is a hyperbolic plane 
$P \subset V$. Then by orthogonalization we can write $V \simeq P \oplus V_1$, where
$V_1$ is the orthogonal complement of $P$ in $V$, so in particular, $\dim V_1 = \dim V - 2$.
If $V_1$ also contains an isotropic vector, we can split off another hyperbolic plane from it.
We can continue this process until we arrive at a decomposition
$$
V \simeq \sum_{j = 1}^m P_j \oplus V_m,
$$
where $\dim V_m = \dim V  - 2m$, and more importantly, $V_m$ is {\it anisotropic}, meaning 
that $V_m$ contains no non-zero isotropic vectors. 

The collection of inner products of the form $B_o \oplus (B_{1,1})^k$, where $B_o$
 is a fixed an isotropic form will be called a {\it Witt tower} of forms. 
If the anisotropic component of $B$ reduces to $\{0\}$, so that $B \simeq (B_{1,1})^m$ 
for some $m$, we call $B$ {\it split}. In other words,
the split forms constitute the Witt tower of $\{0\}$. The isometry group of the split inner product
 of rank $2m$ is denoted $O_{m,m}$. 

Finally, we note that ${\bf F}^{\times}$ acts on inner products, simply by multiplying values of an inner product $B$ by a scalar $s$ in ${\bf F}^{\times}$:
$$
(sB)(\vec v_1 , \vec v_2) = s(B(\vec v_1, \vec v_2)). 
$$
This action is clearly compatible with isomorphism, and with forming direct sums, so it 
defines an action of ${\bf F}^{\times}$ on the Witt semigroup. It is clear that, since 
$$
(s^2B)(\vec v_1 , \vec v_2) = s^2(B((\vec v_1 , \vec v_2)) = B(s\vec v_1 , s\vec v_2),
$$
the subgroup ${\bf F}^{\times 2}$ acts trivially on the Witt semigroup. It is also clear from
 the definition that the isometry groups of $B$ and of $sB$ are isomorphic:
$$
O_B \simeq O_{sB}. 
$$
Clearly, the action of ${\bf F}^{\times}/{\bf F}^{\times 2}$ on the Witt semigroup is transitive
on the generating set of rank one inner products. 

If an inner product is given by a diagonal matrix, then up to multiples, we can arrange that the 
$(1, 1)$ entry of the matrix is 1. Then if $b$ is the $(2,2)$ entry, the inner product on the plane
 spanned by the first two basis vectors is $x_1x'_1 + bx_2x'_2$. This will allow an isotropic
  vector if and only if $-b$ is a perfect square in ${\bf F}$. Otherwise, the equation
 $x^2 = - b$ defines a quadratic field  extension  ${\bf F}_b$ of ${\bf F}$, generated by the 
 element $\sqrt {-b}$. In this case, it is easy to check that if ${\bf F}_b$ is given the basis
 $(1, \sqrt {-b})$, then this inner product is exactly the bilinearization of the norm map [Lam73] from
 ${\bf F}_b$ to ${\bf F}$. Thus, every non-split 2-dimensional inner product is a scalar multiple
 of the norm form of a quadratic extension of ${\bf F}$.
 \vskip .1 in
 Consider now the case when ${\bf F} = {\bf F}_q$ is a finite field, which is the main case of 
 interest in these notes. Then ${\bf F}_q$ has a unique (up to isomorphism) quadratic extension ${\bf F}_{q^2}$, and the norm
 map from ${\bf F}_{q^2}$ to ${\bf F}_q$ is surjective. It follows that ${\bf F}_{q^2}$ with its norm
  map represents the unique non-split rank 2 inner product. Further, any rank 3 inner product 
 must contain an isotropic vector: since any 2 dimensional subspace that is not split will be 
 equivalent to ${\bf F}_{q^2}$ and so will contain a vector whose self inner product is
  the negative of any vector orthogonal to the 2D subspace. 
 
 Thus, for ${\bf F}_q$, the possible anisotropic inner products of positive rank are the two rank one inner products
 parametrized by ${\bf F}^{\times}_q/({\bf F}^{\times}_q)^2$ which is a cyclic group of order two, and the rank two inner product given by the norm map on ${\bf F}_{q^2}$. Thus, there are four Witt towers
 of inner products. Two of them consist of even rank inner products, one the split Witt tower,
 and the other one generated by the norm form on ${\bf F}_{q^2}$. In odd dimensions,
  there are two Witt towers, generated by the two rank one inner products. These have
 the form $B_1 \oplus (B_{1,1})^m$, where $B_1$ is one of the two rank one inner products. 
 Since the two rank one inner products are multiples of each other, and any multiple of a split 
 form is split, the inner products of a given odd rank are multiples of each other, and so will
 have isomorphic isometry groups. The isometry group of the non-split form of rank 2m is 
 not isomorphic to $O_{m,m}$.

 \vskip .05 in
 
 \noindent {\bf Unitary groups}
 
 \vskip .05 in
 
 There is a class of groups that become isomorphic to $GL_n$ after appropriate field extension,
 but that also share some of the attributes of orthogonal groups. 
  
 As above, suppose that $b$ is an element of ${\bf F}$ such that $-b$ is a non-square in 
 ${\bf F}^{\times}$, and let ${\bf F}_b$ denote the quadratic extension of ${\bf F}$ generated by
 $\sqrt {-b}$. For example, the complex numbers ${\bf C}$ can be described as ${\bf R}_1$.
  Denote the Galois automorphism 
 of ${\bf F}_b$ over ${\bf F}$ by $z \mapsto \overline z$; thus 
 $\overline {\sqrt {-b}} = - \sqrt {-b}$. For the finite fields ${\bf F}_q$, it is well known [Lam73] 
 that we may take $b = 1$ when $q$ is congruent to 3 modulo 4, but not when $q$ is
 congruent to 1 modulo 4.
 \vskip .05 in 
 Let $V$ now denote a vector space over ${\bf F}_b$. A function $H(\vec v_1, \vec v_2)$ of 
 two vectors in $V$ is a {\it Hermitian form}  (or more specifically, when necessary,
 an ${\bf F}_b$ Hermitian form), if it is 
 
 \vskip .09 in
 
 \hskip .1 in a) linear in the first variable, and in addition
 
 \vskip .08 in
 
 \hskip .1 in b) $H(\vec v_2, \vec v_1) = \overline {H(\vec v_1, \vec v_2)}.$
 
 \vskip .08 in
 
 \hskip .1 in c) $H$ is non-degenerate, i.e.,  given $0 \neq \vec v_1$ in $V$,
 there is $\vec v_2$ in $V$ such that $H(\vec v_1, \vec v_2) \neq 0$.
 
 \vskip .09 in
 
 Note that, if $V$ is considered as a vector space over ${\bf F}$, the ``real" part of  $H$, 
 namely   ${H(\vec v_1 \vec v_2) + H(\vec v_2, \vec v_1) \over 2}$, is an inner product. 
 Also, the ``imaginary" part  ${(H(\vec v_1 \vec v_2) - H(\vec v_2, \vec v_1))\over 2\sqrt {-b}}$ 
 is a symplectic form.
 
\vskip .05 in
  
The group of linear transformations of $V$ that preserve the Hermitian form $H$ is called
the {\it unitary group} of $H$, and is denoted $U_H$. 
 
 There is a Witt semigroup of ${\bf F}_b$ Hermitian forms, just as there is of inner products. 
 There  is an analog of Witt's Theorem, so this semigroup has the cancellation property.
  Also, every ${\bf F}_b$ Hermitian form is equivalent to one given by a diagonal matrix, 
  or equivalently, every form is a direct sum of rank one forms. Also, if 
  $N({\bf F}_b^{\times}) \subset {\bf F}^{\times}$ is the image of ${\bf F}_b^{\times}$ 
  under the norm map from ${\bf F}_b$ to ${\bf F}$,  the Witt semigroup 
  supports a natural action  of ${\bf F}^{\times}/N({\bf F}_b^\times)$, and the rank one forms are 
  parametrized by the elements of this group (which is a quotient group of 
  ${\bf F}^{\times}/{\bf F}^{\times 2}$. There is also a notion of hyperbolic plane for
  ${\bf F}_b$-Hermitian forms, and any isotropic vector (meaning, a vector with self pairing 
  equal to zero) is contained in a hyperbolic plane. 
  
  For a finite field, ${\bf F_q}$, we already know that there is only one quadratic extension,
 ${\bf F}_{q^2}$, which is also $({\bf F}_q)_b$ for any $b$ such that $-b$ is not a square
 in ${\bf F}_q$. Since the norm map from ${\bf F}_{q^2}$ to ${\bf F}_q$ is surjective, there
 is only one Hermitian form of rank one, given by the formula
  $$
  H_1(z_1, z_2) = z_1\overline{z_2}, 
  $$
  for two elements $z_j$ of ${\bf F}_{q^2}$. It follows that there is a unique isomorphism class
  of Hermitian forms in each dimension. These are organized into two Witt towers, consisting
  of the even dimensional, or split, forms, and the odd dimensional forms. In the first case we 
  denote the associated isometry groups by $U_{m,m}$, while in the second case we denote them by
  $U_{m+1, m}$.    
  
  \vskip .05 in
  
  \noindent {\bf Summary of classical groups over finite fields}
  
  The isomorphism classes of classical groups over a finite field ${\bf F}_q$ are as follows.
  \vskip .1 in
  {\bf a)} General linear groups $GL_n$, $n \geq 0$.
  \vskip .1 in
  
  {\bf b)} Unitary groups, the isometry groups of non-degenerate Hermitian forms. 
  These are organized into two Witt towers. 
  \vskip .07 in
 \hskip .5 in {\bf i)} One tower consists of the even rank forms, 
  which are split, that is, a direct sum of hyperbolic planes. 
  The isometry group of these $2m$-dimensional forms (over ${\bf F}_{q^2}$) 
  is $U_{m,m}$.
  \vskip .08 in
 \hskip .5 in  {\bf ii)} The other Witt tower consists of the odd rank forms. These are the sum of one
  line and an even dimensional form. The isometry group of this forms in $2m+1$ dimensions is $U_{m+1, m}$.
  \vskip .1 in
  
  {\bf c)} Orthogonal groups, the isometry groups of non-degenerate symmetric bilinear forms,
  aka, inner products. These are organized into four Witt towers. 
 \vskip .08 in
 \hskip .5 in {\bf i)} The isometry group of the split form of rank $2m$ 
 is $O_{m,m}$. 
  \vskip .08 in 
\hskip .5 in   {\bf ii)} The isometry group
  of the non-split form of rank $2m+2$ is $O_{m+2, m}$, with the special 
  case of $m = 0$, the anisotropic form, being $O_2$. The group $O_2$ is consists of 
  the norm units in ${\bf F}_{q^2}$, extended by the Galois automorphism of ${\bf F}_{q^2}$
  over ${\bf F}_q$. The group $O_{1,1}$ is isomorphic to ${\bf F}_q^{\times}$, extended
 by the automorphism $x \rightarrow x^{-1}$. 
 \vskip .08 in 
 \hskip .5 in  {\bf iii)} The isometry groups of forms of odd rank. Since, as we have noted, 
  the odd rank inner products are multiples of each other, they 
 have isomorphic isometry groups that we denoted by $O_{m+1, m}$.
  
  \vskip .1 in 
  
  {\bf d)} Symplectic groups, the isometry group of symplectic forms. There is one of these in 
  each even dimension, forming a single Witt tower. They are all split groups. The isometries
  of the $2n$ dimensional symplectic form is denoted $Sp_{2n}$.

\vskip .2 in

{\bf 2. Vector spaces of matrices}

\vskip .15 in

We want to introduce several special vector spaces of matrices and their Pontrjagin duals.

\vskip .1 in 

{\bf 2.1. Rectangular matrices.}
 
 \vskip .05 in 
 
 Consider the space $M_{a, b}({\bf F}) = M_{a,b}$ of $a \times b$ matrices over a field ${\bf F}$
  of coefficients. This can be thought of as the vector 
 space of linear transformations $Hom({\bf F}^b, {\bf F}^a)$ from the vector space ${\bf F}^b$
  of column vectors of length
  $b$ with entries in ${\bf F}$, to the space ${\bf F}^a$.
 
 Matrix multiplication defines mappings $M_{a,b} \times M_{b,c} \rightarrow M_{a,c}$ 
 by the standard formulas [Mostow-Sampson69], [Strang76].
In particular, the spaces $M_{a,a}$ of square matrices are associative algebras. 
The invertible elements of $M_{a,a}$
 form the group
$GL_a$ of invertible $a \times a$ matrices. 

The product group $GL_a \times GL_b$ acts on $M_{a,b}$ by multiplication on the left and right
$$
\alpha(g, h)(T) = gTh^{-1}, \eqno (2.1.1)
$$
for $g$ in $GL_a$, $h$ in $GL_b$, and $T$ in $M_{a,b}$.

On $M_{a,a}$, there is an interesting linear functional, the {\it trace map} [Lang02], [Strang76]. 
If $T = \{t_{jk}\}: 1 \leq j, k \leq a$
 indicates the entries of the $a \times a$ matrix $T$, then
$$
{\rm trace}(T) = \sum_{j = 1}^a t_{jj}, 
$$
is the sum of the diagonal entries. The trace functional is well-known, and can be easily 
checked, to have the following properties.
$$
{\rm trace}(TS) = {\rm trace}(ST), \eqno (2.1.2a)
$$
for any two $a \times a$ matrices $T$ and $S$, and
$$
{\rm trace}(g T g^{-1}) = {\rm trace}(T). \eqno (2.1.2b)
$$
for any $a \times a$ matrix $T$ and any element $g$ of $GL_a$. The equation $(2.1.2b)$
 follows from $(2.1.2a)$. 
Either of these equations characterize the trace functional up to linear multiples.

The vector space dual of $M_{a,b}$ can be realized as $M_{b,a}$, by means of matrix 
multiplication and the trace map.
 For $T$ in $M_{a,b}$ and $S$ in $M_{b,a}$, we define
 $$
 \lambda_S(T) = \lambda_T (S) = {\rm trace}(ST) = {\rm trace} (TS).
 $$
 Here the trace of $ST$ is taken in $M_{b,b}$, while the trace of $TS$ is taken in $M_{a,a}$.
  The resulting value is independent of order of multiplication.
 
 This pairing is compatible with the actions of $GL_a \times GL_b$ on these spaces: 
 for $T$ and $S$ as above and
 $g$ in $GL_a$ and $h$ in $GL_b$, we check that
 $$
 {\rm trace} ((\alpha(h, g)(S)) (\alpha (g, h)(T))) = {\rm trace}((hSg^{-1})(gTh^{-1})) = {\rm trace}((hS)(g^{-1}g) (Th^{-1}))
 $$
 $$
 = {\rm trace}(h(ST)h^{-1}) = {\rm trace}(ST).
 $$
 A fundamental result of linear algebra [Lang02], [Mostow-Sampson69], [Strang76] is that the orbits of $GL_a \times GL_b$
  acting on $M_{a,b}$
  are classified by the  {\it rank} of elements in an orbit. The rank of a matrix $T$ in $M_{a,b}$
   is defined as the dimension 
  of the linear span of the columns of $T$,
  or the dimension of the linear span of the rows of $T$. Equivalently, it is the maximal number
  of linearly independent rows,
  or maximal number of linearly independent columns, of $T$. More geometrically,
  it may be described as the dimension
 of the image space $T({\bf F}^b)$ as a subspace of ${\bf F}^a$. The rank of $T$ in $M_{a,b}$
  is a whole number between
  0 and $\min(a, b)$.
 
 Since $M_{a,b}^* \simeq M_{b,a}$, and since this isomorphism is compatible with the actions
  of $GL_a \times GL_b$ on
 these spaces, it makes sense to talk about the rank of a linear functional on $M_{a,b}$, 
 and the rank also classifies
  the $GL_a \times GL_b$ orbits on $M_{a,b}^*$.
 \vskip .1 in 
 We recall that rank satisfies some simple properties under multiplication and addition of 
 matrices, and also some monotonicity properties
 with respect to decomposition [Lang02], [Strang76]. If $T$ is in $M_{a,b}$ and $S$ is in $M_{b,c}$, then
 $$
 {\rm rank}(TS) \leq \min ({\rm rank} \ (T), {\rm rank} \ (S)).
 $$
 If $T$ and $T'$ are two $a \times b$ matrices, then 
 $$
 {\rm rank}\ (T + T') \leq {\rm rank} \ (T) + {\rm rank} (T').
 $$
 If $b = b_1 + b_2$, there is an obvious isomorphism
 $$
 M_{a,b} \simeq M_{a, b_1} \oplus M_{a, b_2},
 $$
 in which $M_{a, b_1}$ is identified to $a \times b$ matrices for which the final $b_2$ columns
  vanish identically,
 and similarly, $M_{a, b_2}$ is identified to $a \times b$ matrices for which the initial $b_1$
  columns vanish identically.
 If $T$ is in $M_{a, b}$, we let $T_1$ and $T_2$ stand for its components in $M_{a, b_1}$
  and in $M_{a, b_2}$ 
 respectively, giving us a block decomposition of $T$:
 $$
 T = \left [ \matrix {T_1 \ \ \ T_2} \right ], \eqno (2.1.3)
 $$
In this situation, we know that the rank of $T$ dominates the rank of either $T_{\ell}$ [Lang02], [Strang76]:
 $$
 {\rm rank}(T) \geq {\rm rank} (T_{\ell}) \quad {\rm for} \ \ell = 1, 2. \eqno (2.1.4)
 $$
 Also, for $T_{\ell}$ in $M_{a, b_{\ell}}$, the rank of $T_{\ell}$ is the same, whether $T_{\ell}$
  is considered
 as an element of $M_{a, b_{\ell}}$, or as the element of $M_{a, b}$ produced by extending 
 $T_{\ell}$ by adding to it the zero matrix
 of size $M_{a, b - b_\ell}$ in the appropriate block:
 $$
 {\rm rank} \ (T_1) = {\rm rank}\ (\left [ \matrix {T_1 \ \ 0 } \right ], \quad {\rm and} \quad
 {\rm rank }\ (T_2) = {\rm rank} \ \left [ \matrix { 0 \ \ T_2} \right ].
 $$
 We can also break the matrices $M_{a, b}$ into blocks by rows rather than columns, 
 and indeed, we
 can break them up by both at the same time. Thus, if we break up $a= a_1 + a_2$ and 
 $b = b_1 + b_2$ into
 sums, we have the straightforward direct sum decomposition
 $$
 M_{a, b} \simeq M_{a_1, b_1} \oplus M_{a_1, b_2} \oplus M_{a_2, b_1} \oplus M_{a_2, b_2}, 
 $$
 corresponding to the block decomposition of a matrix $T$ into pieces $T_{\ell m}$:
 $$
 T = \left [ \matrix {T_{11} \ \ \ T_{12} \cr T_{21} \ \ \ T_{22} } \right ],
 $$
 where each $T_{\ell m}$ is a matrix of the appropriate size.
 
 In this context, rank satisfies the monotonicity property
 $$
 {\rm rank} \ (T) = {\rm rank}\ \left ( 
  \left [ \matrix {T_{11} \ \ \ T_{12} \cr T_{21} \ \ \ T_{22} } \right ] \right ) \geq
 {\rm rank} \ \left (  \left [ \matrix {0 \ \ \ \ \  0\  \cr 0  \ \ \ T_{22} } \right ] \right ) 
 = {\rm rank} \ (T_{22}).
 $$
 A similar relation holds for each sub-block $T_{\ell m}$ of $T$. These relations can also be 
 extended in a straightforward way to allow for decompositions of $M_{ab}$ into more blocks,
  horizontally or vertically or both.
 \vskip .1 in
 
 {\bf 2.2. Symmetric and skew-symmetric matrices.}
 \vskip .1 in
 
 In this section, we assume for convenience that the characteristic of the base field ${\bf F}$ is not 2.
 
 Consider next the space $S^2_n({\bf F})= S^2_n$ of symmetric $n \times n$ matrices with
  coefficients in a field ${\bf F}$. We can regard a matrix $S$ in $S^2_n$ as defining
 symmetric bilinear form $B_S( \ , \ )$ on ${\bf F}^n$ by the well-known recipe
 $$
 B_S( \vec x, \vec y) = \vec y^t S \vec x = \sum_{j, k = 1}^n s_{jk} x_j y_k. \eqno (2.2.1)
 $$
 Here $\vec x$ and $\vec y$ are the column vectors
 $$
 \vec x = \left [ \matrix {x_1 \cr x_2 \cr . \cr . \cr . \cr x_n} \right ] \quad {\rm and} \quad 
 \vec y = \left [ \matrix {y_1 \cr y_2 \cr . \cr . \cr . \cr y_n} \right ],
 $$
 and the $s_{jk} = s_{kj}$ are the entries of the matrix $S$.
 
 There is a natural action $\beta$ of $GL_n$ on $S^2_n$, derived from its action on ${\bf F}^n$. This action comes 
 in two mutually dual forms; since they have the same orbits, which is the main point of interest
  here, we will be somewhat cavalier about distinguishing between the two. The actions are converted
  into each other 
 by means of the automorphism 
 $g \leftrightarrow (g^{-1})^t$ of $GL_n$, where $g^{-1}$ indicates the inverse of the invertible
  $n \times n$ matrix $g$, 
 and $g^t$ denotes the transpose of $g$. The action $\beta$ is defined by the condition
 $$
 B_{\beta (g)(S)} (\vec x, \vec y) = B_S(g^{-1}\vec x, g^{-1}\vec y) = 
 B_{(g^{-1})^tSg^{-1} }(\vec x, \vec y),
 $$
 or more succinctly,
 $$
 \beta(g)(S) = (g^{-1})^t S g^{-1}. \eqno (2.2.3)
 $$
 In analogy with arbitrary rectangular matrices, the dual vector space $S^2_n$ can be regarded
 as its own dual, by means of the pairing
 $$
 \lambda_T(S) = {\rm trace} \ (TS).
 $$
 The dual action of $GL_n$ on $S^2_n$ when it is considered as its own dual will be denoted
  by $\beta^*$. 
 We can compute that 
 $$
 (\beta^*(g)(\lambda_T))(S) = \lambda_T(\beta(g^{-1}(S)) = 
 \ {\rm trace} \ (T(g^t Sg)) \ = \ {\rm trace} \ (gTg^tS)
 \ = \ \lambda_{\beta(g^tTg)}(S).
 $$
 In other words
 $$
 \beta^*(g)(T) = g^t Tg. \eqno (2.2.4)
 $$
 
 \vskip .1 in
 
 The orbits of the action of $GL_n$ on $S^2_n$ are the isomorphism, or isometry, classes
 of symmetric bilinear forms 
 $B_S$ on  ${\bf F}^n$. These have been discussed in \S 1 above. The first invariant of the 
 form $B_S$ is its {\it rank}, which in fact coincides with the rank of $S$ as a matrix.
 There are several isomorphism classes of a given rank, depending
 on the base field ${\bf F}$.  The isomorphism classes are organized into the Witt semigroup,
 which is graded by rank, and the elements of the Witt semigroup of a given rank parametrize 
 the orbits of symmetric matrices under the actions (2.2.3) or (2.2.4).   
   \vskip .05 in 
 The notion of {\it orthogonal direct sum} of symmetric bilinear forms, as discussed in \S 1,
 corresponds to forming block diagonal sums of matrices.
 If $S_1$ is a symmetric $n_1 \times n_1$
 matrix, and $S_2$ is a symmetric $n_2 \times n_2$ matrix, then the direct sum 
 $B_{S_1} \oplus B_{S_2}$ 
 is the symmetric bilinear form on ${\bf F}^{n_1 + n_2}$ defined by the matrix
 $$
 S_1 \oplus S_2 = \left [ \matrix {S_1 \ \ \ 0 \cr 0 \ \ \ S_2} \right ]. \eqno (2.2.5)
 $$
 It is easy to see that 
 $$
 {\rm rank} \ (B_{S_1} \oplus B_{S_2}) = {\rm  rank} \ (B_{S_1}) + {\rm rank} \ (B_{S_2}).
  \eqno (2.2.5a)
  $$
 This includes the special case
 $$
 {\rm rank} \left [ \matrix {0  \ \ \ 0 \cr 0 \ \ \ S_2} \right ] = {\rm rank} \ (S_2). \eqno (2.2.5b)
 $$
 Also, a symmetric bilinear form whose defining matrix is diagonal is effectively expressed as the direct sum of one-dimensional forms.
 
 In the orthogonal direct sum $B_{S_1} \oplus B_{S_2}$, the two subspaces ${\bf F}^{n_1}$
  and ${\bf F}^{n_2}$ are
 orthogonal to each other; that is, the pairing of any vector $\vec v_1$ from ${\bf F}^{n_1}$
  with any vector
 $\vec v_2$ from ${\bf F}^{n_2}$ is zero. Somewhat more generally, we may consider
  the situation in which we have
 a direct sum decomposition ${\bf F}^{n_1 + n_2} \simeq {\bf F}^{n_1} \oplus {\bf F}^{n_2}$, 
 but there may be vectors in
 the two summands whose pairing is non-zero. In this case, we can decompose 
 the representing matrix $S$ into
 components
 $$
 S = \left [ \matrix {S_{11} \ \ \ S_{12} \cr S_{12}^t \ \ \  S_{22} } \right ] \eqno (2.2.6)
 $$
 In this situation, we have the dominance relation between ranks, that is,
 $$
 {\rm rank} \ (S) \geq {\rm rank} \ S_{ij}. \eqno (2.2.7)
 $$
 \vskip .05 in 
 
 The {\it hyperbolic plane} is the (isomorphism class represented by) the form with matrix
 $$
 \left [ \matrix {1 \ \ \ \ \ 0 \cr 0 \ \ -1} \right ] \quad {\rm or} 
 \quad \left [ \matrix {0 \ \ \ 1 \cr 1 \ \ \ 0} \right ]. 
 $$
 
  \vskip .1 in
  
  A discussion similar to the above can be carried out for skew-symmetric bilinear forms.
   In general, if we use an arbitrary 
  $n \times n$ matrix $T$ to define a bilinear form $B_T$ on ${\bf F}^n$ via the formula (2.2.1),
  then
  $$
  B_T(\vec y, \vec x) = B_{T^t} (\vec x, \vec y). \eqno (2.2.8)
  $$
  The symmetry of the bilinear forms discussed above resulted directly from the symmetry
   of the matrices $S$ under 
  discussion. We may replace the symmetry condition with a condition of skew-symmetry,
   creating a situation that is similar 
  in many ways to the symmetric case, but overall is somewhat simpler.

Consider the space $\Lambda^n({\bf F}) = \Lambda^2_n$ of $n \times n$ skew symmetric
 matrices $L$: $L^t = -L$. 
If we use a skew-symmetric matrix $L$ to define a bilinear form $B_L$ by means of the formula
 (2.2.1), then $B_L$ will
 be skew-symmetric in the sense that 
$$
B_L(\vec x, \vec y) = - B_L(\vec y, \vec x).
$$
This results directly from formula (2.2.8) and the skew-symmetry of $L$. 
We have natural actions of $GL_n$ on $\Lambda^2_n$, given by the same formulas (2.2.3)
 and (2.2.4)  as for the symmetric
 case. Again the main invariant of these actions is the rank of the matrix $L$, which again has
  more or less the same 
 interpretation in terms of the bilinear form $B_L$ as for the symmetric case. The remarks
  about forming direct sum 
 decompositions, and the inequalities of ranks described by formulas (2.2.5) through (2.2.7)
  also remain valid
  in the skew-symmetric case.

The main differences from the symmetric case are two [Lang02]. First, whereas the arithmetic of the 
coefficient field ${\bf F}$ affects 
the $GL_n$ orbit structure in $S^2_n$, in the skew-symmetric case, there is only one possible 
orbit of a given rank. That is,
all skew-symmetric bilinear forms of a given rank are isometric to each other. Second, 
skew-symmetric forms can only have
 even rank: a non-degenerate skew-symmetric bilinear form, which is commonly called a
  {\it symplectic form}, 
 can only exist
  on an even-dimensional vector space. Thus, every skew-symmetric bilinear form is isomorphic
   to the orthogonal direct sum
   of a certain number of copies of the 2-dimensional form represented by
$$
\left [ \matrix {0 \ \ -1 \cr 1 \ \ \ \ \ 0} \right ],
$$
together with copies of the one-dimensional form that is identically zero.

\vskip .1 in

{\bf 2.3. Hermitian matrices.}
\vskip .1 in
We also consider Hermitian matrices. These are square matrices  $T$, with entries $t_{j\ell}$
 in the quadratic field  extension ${\bf F}_b$, as described in \S 1. To be Hermitian, 
 the matrix $T$ should satisfy 
 $$
 T^t = \overline T,
 $$
 where $T^t$ denotes the transposed matrix to $T$, and $\overline T$ denotes the matrix
gotten by applying the Galois automorphism $t \rightarrow \overline t$ of ${\bf F}_b$ over
${\bf F}$ to the entries of $T$. Stated in terms of entries, the condition of Hermitianness is
$$
t_{\ell j}  = \overline t_{j\ell}. 
$$
We define the {\it Hermitian conjugate} matrix to $T$ by the equation
$$
T^* =  \overline {T^t} = (\overline T)^t .
$$
We let $H_n({\bf F}_{q^2})$ denote the space of $n \times n$ Hermitian matrices 
with coefficients in ${\bf F}_{q^2}$. The group $GL_n({\bf F}_{q^2})$ acts on the Hermitian matrices, in direct analogy
with the actions (2.2.3) and (2.2.4). The formulas for these analogs are
$$
\beta(g)(T) = (g^*)^{-1}Tg^{-1}, \ \ \ \ {\rm and} \ \ \ \ \beta^*(g)(T) = g^*Tg. \eqno (2.3.1)
$$
Since there is only one isomorphism class of Hermitian form in each dimension, 
there is just one $GL({\bf F}_{q^2})$ orbit in $H_n({\bf F}_{q^2})$ for matrices of a given rank.

\vskip .1 in

{\bf 2.4. Pontrjagin duals of matrix groups.}
 \vskip .1 in
 
 In the context of the preceding two subsections, if the field ${\bf F}$ of coefficients is locally
  compact, so either ${\bf R}$ or
 ${\bf C}$, or a $p$-adic field, or finite [Weil74], then the additive groups of the various spaces 
 of matrices 
 discussed above, 
 or indeed, of any finite-dimensional vector space over ${\bf F}$, is also locally compact.
  Hence, we can consider its unitary
  representations, and in particular, its Pontrjagin dual group of unitary characters [Curtis-Reiner62], [Isaacs76]. 
 
 It is well-known [Weil74], that for a finite-dimensional vector space $V$ over ${\bf F}$,
  the Pontrjagin dual group $\widehat V$ 
 can be identified with the dual vector space $V^*$. Precisely, one should fix a
  ``basic character" 
 $\chi_o$
  of the additive group of ${\bf F}$. Given this choice, the mapping
  $$
  \gamma: V^* \longrightarrow \widehat V, \eqno (2.4.1)
  $$
  defined by
  $$
  \gamma (\lambda) (\vec v) = \chi_o (\lambda (\vec v)),
  $$
  is an isomorphism of groups.  To be a little more explicit in one of the  main cases 
  of interest, namely additive groups of matrices, for each $T$ in $M_{b, a}$, we can define
   a character
  $\gamma_T$ on $M_{a, b}$ by the formula
  $$
  \gamma_T(S) = \chi_o({\rm trace} (TS)) \ \ \ \ {\rm for} \ \ \ S \in M_{a,b}. \eqno (2.4.2)
  $$
  The same formula works also for symmetric matrices and for skew-symmetric matrices, 
  and also for Hermitian matrices.
 
 The map (2.4.1) is functorial, in the sense that it is compatible with the maps induced by
  linear transformations $T: V \rightarrow U$ of vector spaces over ${\bf F}$.
  In particular, if $V = M_{a,b}$ is the space of $a\times b$ matrices, then $\widehat {M}_{a,b}$
   allows an action of
  $GL_a \times GL_b$. This is the transfer via $\gamma$ of the dual action $\alpha^*$
   on $M_{b. a}$ of formula (2.1.1). We will denote this action by $\widehat \alpha$. Similarly, we have actions $\widehat \beta$
   and $\widehat \beta^*$ of 
  $GL_n$ on $\widehat S^2_n$ and on $\widehat \Lambda^2_n$. The orbits for these actions will be
   the images under the maps
  $\gamma$ of the orbits in the various vector spaces, as described in \S 2.1 and \S 2.2.
  In particular, it makes sense to talk 
about the rank of a unitary character of one of these spaces of matrices, and to write that 
rank($\gamma_T$) = rank ($T$), for $T$ as in formula (2.4.2). The rank of a character
 of $M_{a.b}$ will 
completely determine its $GL_a \times GL_b$ orbit in $\widehat {M}_{a,b}$. Similarly, the rank
 of a character of $\Lambda^2_n$
will determine its $GL_n$ orbit in $\widehat \Lambda^2_n$. In $\widehat S^2_n$, the rank of a character
 determines its orbit
 up to a finite number of choices, independent of $n$ and the rank. For a character of $S^2_n$
 of rank $k$, we will call the isomorphism class of the symmetric bilinear form that corresponds to it, the {\it type} 
 of the character.
 
 \vskip .2 in
  
 {\bf 3. Matrix subgroups of classical groups}
 
 \vskip .15 in 
 
  This work may be regarded as part of the tradition of understanding representations
   of a reductive group by
 examining their restrictions to unipotent subgroups. A notable example of this approach
  is the theory of Whittaker 
 functions and Whittaker vectors [Knapp86], [Wallach92], [Casselman-Hecht-Milicic98] which has been used to understand
  ``typical" or ``generic"  representations.
 This was shown to be particularly effective for $GL_n$ by S. I. Gelfand [Gelfand70], who showed
 that cuspidal representations 
 of $GL_n$ over finite fields were generic, and in fact, were of the minimal dimension possible
  for generic representations. 
 
 Whittaker functions look at the occurrence of characters of the maximal unipotent subgroup. 
 Since the maximal 
 unipotent subgroup is, typically, non-abelian, characters form only a very limited subset
  of their representations.
 The approach of the theory of rank is to look at non-maximal unipotent subgroups that are
  abelian. The restriction
 of a representation of the full reductive group to such a subgroup will then decompose entirely
  into characters,
 appearing with various multiplicities, and one can look at the possibilities for such
  a decomposition. 
 There are many subgroups of classical groups that are in a straightforward way isomorphic
  to additive groups of matrices. We will describe some of them. These notes are mainly concerned with characters of unipotent radicals of
maximal parabolic subgroups. Note that these subgroups are abelian, so all their irreducible
representations are one-dimensional.  In \S 10.4, we consider an analog of Whittaker vectors
for unipotent radicals of certain parabolic subgroups, not the Borel subgroup, 
and not a maximal parabolic subgroup, and show that these can be useful for describing 
some classes of representations.   
  
 \vskip .1 in
 {\bf 3.1. The $GL_n$ case.}
\vskip .1 in 
 
 Consider the parabolic subgroup $P_{a, b} \subset GL_{a+b}$ consisting of block upper
  triangular matrices
 of the form
 $$
 P_{a,b} = \left \{ \left [ \matrix {A \ \ \ T \cr 0 \ \ \ B} \right ] \right \},
 $$
 where $A$ is in $GL_a$, $B$ is in $GL_b$, and $T$ is in $M_{a,b}$. 
 
 The group $P_{a,b}$ has a semidirect product decomposition, (aka {\it Levi decomposition})
 $$
 P_{a,b} \simeq M_L(a,b) \cdot U_{a,b} = M_L \cdot U,
 $$
 where the {\it Levi component} $M_L$ consists of block diagonal matrices:
 $$
 M_L =  \left \{ \left [ \matrix {A \ \ \ 0 \cr 0 \ \ \ B} \right ] \right \},
 $$
and $U$ consists of unipotent upper triangular matrices:
$$
U =  \left \{ \left [ \matrix {I_a \ \ \ T \cr 0 \ \ \ I_b} \right ] \right \}.
$$
Here $I_a$ and $I_b$ are the identity matrices of sizes $a \times a$ and $b \times b$
 respectively.

The group $U$ is isomorphic to the space $M_{a,b}$ of $a\times b$ matrices, via the mapping 
$T \mapsto  \left [ \matrix {I_a \ \ \ T \cr 0 \ \ \ I_b} \right ]$. Indeed, it is easy to check
 via matrix multiplication 
that matrix multiplication in $U_{a.b}$ corresponds to addition in $M_{a,b}$:
$$  
\left [ \matrix {I_a \ \ \ T \cr 0 \ \ \ I_b} \right ]  \left [ \matrix {I_a \ \ \ T' \cr 0 \ \ \ I_b} \right ] =
 \left [ \matrix {I_a \ \ \ T + T' \cr 0 \ \ \ I_b} \right ].
 $$
 Also, the action by conjugation of $M_L$ on $U$ is identified in a straightforward way
 to the action $\alpha$ of formula
 (2.1.1):
 $$
 \left [ \matrix {A \ \ \ 0 \cr 0 \ \ \ B} \right ] \left [ \matrix {I_a \ \ \ T \cr 0 \ \ \ I_b} \right ] 
 \left [ \matrix {A \ \ \ 0 \cr 0 \ \ \ B} \right ] ^{-1} = 
 \left [ \matrix {I_a \ \ \ ATB^{-1} \cr 0 \ \ \ \ \ \ \ I_b\ \ \ \ } \right ]. 
 $$
 More generally, we may consider  inside $GL_n$ a general parabolic subgroup 
 of block upper triangular matrices:
 $$
P =  \left \{ \left [ \matrix {A_1 \ \ \ T_{12} \ \ \ T_{13} \ \ \ \ . \ \ \ \ . \ \ \ \ . \ \ \ T_{1k} \cr
  \ \ 0 \ \ \ \ A_2 \ \ \ T_{23} \ \ \ T_{24}  \ \ \ . \ \ \ \  .  \ \ \   T_{2k} \cr
  \ \ 0 \ \ \ \ \ 0 \ \ \ \ A_3 \ \ \ \ \  . \ \ \ \ . \ \ \ \  . \ \ \ \ T_{3k} \cr
  \ \ \ . \ \ \ . \ \ \ . \ \ \ . \ \ \ . \ \ \ . \cr
  \ \ \ . \ \ \ . \ \ \ . \ \ \ . \ \ \ . \ \ \ . \cr
   \ \ 0 \ \ \ \ \ 0 \ \ \ \ . \ \ \ \ . \ \ \ \ . \ \ \ \ . \ \ \ 0 \ \ \ A_k } \right ] \right \} \eqno (3.1.1)
   $$
   If the matrices $A_j$ in the $j$-th diagonal block are of size $a_j \times a_j$, 
   then the matrices $T_{jk}$ 
 define a group isomorphic to the $a_j \times a_k$ matrices. Moreover, it is not hard to argue
 that all such copies of the
  $a \times b$ matrices, for given $a$ and $b$ such that $a+b \leq n$, are conjugate
  to each other. We will call such 
   subgroups of $GL_n$, {\it standard matrix subgroups}. In this context, unipotent radicals
 of the parabolic subgroups 
  with only two diagonal blocks, introduced above, produce maximal examples of such standard
 matrix subgroups
   of $GL_n$. The unipotent radicals $U_{a, b}$ are the unique examples of a standard matrix
  subgroup $U_{a, b}$ for pairs 
  $(a, b)$ with $a + b = n$. We note that, in fact, each standard subgroup is the unipotent 
radical of a maximal parabolic subgroup of a subgroup isomorphic to $GL_{a_j + a_k}$
 in $GL_n$.
  
  Terminology: we will say that a standard subgroup of $GL_n$ that is isomorphic
 to the matrices $M_{a,b}$ has {\it size} $(a,b)$.
  
  \vskip .1 in
  
  {\bf 3.2. The $Sp_{2n}$ case.}
  \vskip .1 in 
  The group $Sp_{2n}$ is the isometry group of a symplectic form on a 
  ${\bf F}^{2n}$. As noted in \S 2.2,
  all such forms are isomorphic, so this description does specify $Sp_{2n}$ uniquely. 
  
 Here is a more explicit description of $Sp_{2n}$.
 On $W = {\bf F}^{2n}  \simeq {\bf F}^n \oplus {\bf F}^n$, we use coordinates
 $$
  \vec w = \left [ \matrix {x_1 \cr x_2 \cr . \cr . \cr . \cr . \cr x_n
  \cr y_1 \cr y_2 \cr . \cr . \cr . \cr . \cr y_n} \right ]. 
 \eqno (3.2.1)
 $$
 We define the skew symmetric bilinear pairing $< \ , \ >$ on $W$ by the recipe
 $$
 <\vec w, \vec w'> = \sum_{j = 1}^n x_j y'_j - y_j x'_j.\eqno (3.2.2)
 $$
 Then $Sp_{2n}({\bf F}) = Sp_{2n}$ is realized as the group of linear transformations
  of $W$ that preserve the form  $< \ , \ >$.
 
 It is easy from the above description of $Sp_{2n}$ to check that it contains the subgroup
 $$
 P = P_n \simeq M_n\cdot U_n,\eqno (3.2.3a)
 $$
 Here $M_n$ is the group consisting of block diagonal matrices
 $$
 \left [ \matrix {A \ \ \ \ \ \ 0\ \  \cr \ 0 \ \ (A^{-1})^t} \right ],\eqno (3.2.3b)
 $$
 where $A$ is an invertible $n \times n$ matrix, and $A^t$ indicates the transpose of $A$.
  The group $U_n$ indicates the group of unipotent block upper triangular matrices
 $$
 \left [ \matrix {I \ \ \ S \cr 0 \ \ \ I} \right ] \eqno (3.2.3c)
 $$
 where $S$ indicates a symmetric $n \times n$ matrix.
 Again, the action of $A_n$ on $U_n$ is the appropriate action (the action $\beta$ of formula
  (2.2)), of $GL_n$ on symmetric matrices:
 $$
 \left [ \matrix {A \ \ \ \ \ \ 0\ \  \cr \ 0 \ \ (A^{-1})^t} \right ]  
 \left [ \matrix {I \ \ \ S \cr 0 \ \ \ I} \right ]
 \left [ \matrix {A \ \ \ \ \ \ 0\ \  \cr \ 0 \ \ (A^{-1})^t} \right ]^{-1} 
   = \left [ \matrix {\ I \ \  ASA^t \cr 0 \ \ \ \ \ I\ \ }  \right ].
   $$
 Thus, the group $U_n$ is naturally identifiable with the symmetric matrices $S^2({\bf F}^n)$,
 and its conjugacy classes under $M$ 
correspond to isomorphism classes of symmetric bilinear forms. Again, the most important
 invariant of such a form
 is its rank. As we have noted in \S 2.2, rank is not the sole determinant of a symmetric bilinear
  form. For our fields of interest - local fields and finite fields - there are a finite number
of isomorphism classes of each rank. 

\vskip .1 in

{\bf 3.3. The $O_{n, n}$ case.}
\vskip .1 in

A discussion similar to that just given for $Sp_{2n}$ applies to split orthogonal groups. 
On ${\bf F}^{2n}$, 
we again define coordinates as in formula (3.2.1), and consider the symmetric bilinear form
 (aka inner product)
  $$
  (\vec w, \vec w') = \sum_{j = 1}^n x_j y'_j + y_j x'_j.\eqno (3.3.1)
  $$
  We denote the isometry group of this  inner product as $O_{n,n} ({\bf F}) = O_{n,n}$.
 Then $O_{n,n}$ has the subgroup
  $$
  P' = P'_n  \simeq M_n \cdot U'_n,
  $$
  where $M_n$ is, just as in the case of the symplectic group, equal to group of block diagonal
  matrices of formula (3.2.3b),
   and  $U'_n$ consists of block upper triangular unipotent matrices
  $$
  U'_n = \left \{ \left [ \matrix {I \ \ \ L \cr 0 \ \ \ I} \right ] \right \},\eqno (3.3.2)
  $$
  where the matrix $L$ is $n \times n$ and skew-symmetric. Again the action by conjugation
   of $M$ on 
  $U'$ is the standard $GL_n$ action $\beta$ of formula (2.2.4) on skew-symmetric matrices. 
  Again, the conjugacy classes will be separated by
  rank, and there is only one non-degenerate symplectic form of a given rank, and this rank
   must be even. 
   \vskip .1 in
   
   {\bf 3.4. The $U_{n,n}$ case.}
   \vskip .1 in
   The even dimensional unitary groups are quite parallel to $Sp_{2n}$ and $O_{n,n}$. 
   On $({\bf F}_{q^2})^{2n}$, we again use coordinates as in formula (3.2.1), except now 
   the entries can be in ${\bf F}_{q^2})$. We use the Hermitian analog of the pairings (3.2.2)
   and (3.3.1):
   $$
   (\vec w, \vec w') = \sum_{j = 1}^n x_j \overline {y'_j} + y_j \overline {x'_j},
   $$
   The isometry group of this form is ${\bf U}_{n,n}$. (We use boldface for this ${\bf U}$ in order 
   to distinguish it from the $U$s used to indicate unipotent subgroups).
   
   The subgroup of  ${\bf U}_{n,n}$ that stabilizes the subspace defined by setting all the $y_j$
   equal to 0 is
   $$
   P'' = P''_n \simeq M''_n \cdot U''_n,
   $$
   where $M'' \simeq GL_n({\bf F}_{q^2})$ consists of block diagonal matrices of the form
   (3.2.3b), except that now the matrix $A$ belongs to $GL_n({\bf F}_{q^2})$ rather
   than to $GL_n({\bf F}_q)$. The group $U''_n$ consists of block upper triangular matrices
   analogous to those of (3.2.3c) or (3.3.2):
   $$
   \left [ \matrix {I \ \ \ T \cr 0 \ \ \ I} \right ],
   $$
   except now the matrix $T$ is Hermitian rather than symmetric or skew symmetric. 
   The action by conjugation of $M''$ on $U''$ is the natural action of $GL_n({\bf F}_{q^2})$
   on the Hermitian matrices, as described in formula (2.3.1).
  
 \vskip .2 in
 
 {\bf 4. The $U$-rank of representations of classical groups}
 
 \vskip .15 in
 
 We want to use the discussion of the preceding sections to help understand
  the representations
 of classical groups. The case of finite fields is of interest in its own right, and involves fewer
  technicalities  than for local fields. It will be the focus of this exposition.
 
 \vskip .1 in
 
 {\bf 4.1. Rank of representations of matrix groups.}
 \vskip .1 in
 
 With the understanding that the coefficient field ${\bf F}$ is finite, any matrix group such
  as $M_{a,b}$  (or $S^2_n$ or $\Lambda^2_n$) will also be a finite group. 
  Then a finite dimensional representation
 $\pi$ of $M_{a, b}$ will be a sum of characters, appearing with various multiplicities.
 We will say that a character
 $\chi$ of the group {\it appears in $\pi$} if its multiplicity is at least 1, that is, if there
 is an  eigenvector for $M_{a,b}$
 in the space of $\pi$ whose associated eigenvalue under $\pi (M_{a, b})$ is $\chi$.  
 
 As discussed in \S 2.4, each character of $M_{a, b}$ has an associated rank. 
 We will say that the {\it rank} of
 the representation $\pi$ is the maximum of the ranks of the characters that appear in $\pi$.
 {\underbar Mutatis mutandem}, this definition applies also to the groups $S^2_n$
  and $\Lambda^2_n$.
 
 \vskip .1 in
 
 {\bf 4.2.The $U$-rank of representations of classical groups.}
 \vskip .1 in
 
 {\bf Case of $GL_n$}:
 \vskip .05 in
 
 Consider a representation $\pi$ of the group $GL_n({\bf F}) = GL_n$. Look at the restriction to one of the standard 
 matrix subgroups $U$,  as described in formula (3.1.1). Then
 $$
 U \simeq M_{a,b},
 $$ 
 for some pair of whole numbers $(a, b)$ with $a + b \leq n$. Moreover, the normalizer of $U$ acts on $U$ by means of the standard action $\alpha$ (formula (2.1.1)) of the full group $GL_a \times GL_b$. Since the restriction $\pi_{|U}$ of $\pi$  to $U$ can be thought of as obtained by first restricting $\pi$ to the normalizer of $U$ in $GL_n$, and then restricting this to $U$, it follows that any two characters in $\widehat U$ belonging to the same $GL_a \times GL_b$ orbit will appear with the same multiplicity in $\pi_{|U}$. 

For a given standard matrix group $U$ and a representation $\pi$ of $GL_n$, we denote by $rank_U(\pi)$, the largest rank of characters of $U$ appearing in $\pi_{|U}$, and call it the {\it rank} of $\pi$ with respect to $U$. This number is automatically no larger than $\min(a,b)$. If $rank_{U}(\pi)$ is strictly less than $\min(a,b)$, we will say that $rank_{U}(\pi)$ is {\it small}, or that $\pi$ has {\it low rank} for $U$.

\vskip .1 in

{\bf Definition 4.2.1.}  We define the {\it $U$-rank} of $\pi$ to be the largest of the numbers $rank_{U}(\pi)$, as $U$ varies over all standard matrix subgroups $U$ of $GL_n$. We also will say that $\pi$ has {\it low (or small) $U$-rank} if it has low (or small) rank for some $U$.  
  
\vskip .1 in
 
Note that the $U$-rank of $\pi$ is bounded by $\frac{n-1}{2}$ if $n$ is odd, and by $\frac{n}{2}$ if $n$ is even. Moreover, the $U$-rank $k$ is smaller than the above upper bounds (i.e., $k < \frac{n-1}{2}$) if and only if $\pi$ is of low $U$-rank.
Finally, it is satisfying to observe that although the above definition may seem to give rise to a collection of competing candidates for the rank of a representation $\pi$, in fact, all possible values for low rank will agree.
  
  \vskip .1 in
  
 {\bf Proposition 4.2.2}. Suppose that for standard subgroups $U \simeq M_{a,b}$
   and $U' \simeq M_{c,d}$, 
  the representation $\pi$ of $GL_n$ has low rank. Then $rank_U(\pi) \ = \ rank_{U'}(\pi)$.

  \vskip .1 in

  {\it Proof}. We proceed in several steps.
  
\vskip .07 in  

  {\it Step 1.} First consider the situation when $U$ is expressed as a sum of two standard
   subgroups, $U_1$ and $U_2$, obtained by decomposing one of the blocks defining $U$ into two pieces.
  Namely, the parabolic that defines $U$ is the one on the right below, instead of the one that appears on the left, 
  $$
  \left [ \matrix {A \ \ \ X \cr 0 \ \ \ B} \right ] \longmapsto
  \left [ \matrix {A \ \ \ X_1 \ \ \ X_2 \cr 0 \ \ \ B_1 \ \ \ 0 \cr 0 \ \ \ 0 \ \ \ B_2} \right ].
  $$
 Here we have broken up the block of $B$ into two sub-blocks, with corresponding 
 decomposition of $X$ into
 $X \simeq X_1 \oplus X_2$. Note that, by breaking up the $A$ block rather than the
  $B$ block, we could subdivide $U$ vertically rather than horizontally.
 
 Now consider a representation $\pi$ of $GL_n$ that has low rank for $U$, 
 and set $rank_U(\pi) = k$. Then the characters
  of $U$ that appear in $\pi_{|U}$ all have rank at most $k$, and some of them do have rank
 $k$. Also, we know that all rank $k$  characters of $U$ will appear in $\pi_{|U}$. 
 
 If $U$ has size $(a, b)$, then $U_j$ will have size $(a, b_j)$, for two numbers $b_1$
  and $b_2$ such that $b_1 + b_2 = b$. Note that $\pi_{|U_j} = (\pi_{|U})_{|U_j}$.
  According to formulas (2.1.3) and (2.1.4), together with the discussion of \S 2.4, this implies
   that any character of $U_j$ that appears in $\pi_{|U_j}$ will have rank at most $k$. 
 Since the assumption that $rank_U(\pi)$ is low
 implies that $k < a$, the considerations of formulas (2.1.3), (2.1.4) and \S 2.4 also tell us
  that if $k \leq b_j$, then characters
 of rank $k$   will indeed appear in $\pi_{|U_j}$. If $k < b_j$, then $\pi_{|U_j}$ will have rank
  $k$, and will be of low $U_j$-rank. 
 
 In other words, the property of $\pi$ being of low rank, as well as the value of this rank,
  transfers from any standard subgroup
 $U$ to either of its components $U_j$ for which this transfer is allowed by the size
  of the matrices belonging to $U_j$.
 
\vskip .07 in
 
 {\it Step 2.} Using this argument twice, once for decomposing $U$ horizontally, as just above, 
 and once for decomposing it vertically, we arrive at the following conclusion:
 if $\pi_{|U}$ has low $U$-rank equal to $k$, then we can find a standard subgroup 
 $U'' \subset U$ such that $U''$ has size $(k+1, k+1)$ and $\pi_{|U''}$ has rank $k$, 
 which is {\underbar a fortiori} low rank for $U''$.
 
 \vskip .07 in
 
{\it Step 3.} Now we use the observation that all standard subgroups $U$, of size $(a, b)$, 
 are conjugate to each other inside $GL_n$. Thus all these subgroups will give the same measure of $U$-rank,
 and will all agree whether this rank is low or not. In particular, if one  standard subgroup of size $(k+1, k+1)$ 
 says that a representation $\pi$ of $GL_n$ has rank $k$, then all such groups will agree with this. But we have 
 just seen that, if a given standard subgroup $U$ of size $(a, b)$ with $\min(a, b) > k$ says that the representation 
 $\pi$ of $GL_n$ has rank $k$, then also there is a standard subgroup $U''$ of $U$, of size $(k+1, k+1)$, that agrees 
 with this conclusion; and vice versa. Combining these two observations, we conclude that all sufficiently large 
 standard  subgroups are in consensus about whether a representation $\pi$ of $GL_n$ has low rank, and what this rank is. 

\Qed 
\vskip .1 in

{\bf Cases of $Sp_{2n}$ and $O_{n,n}$}:
\vskip .05 in

If $G = Sp_{2n}$, then as described in \S 3.2, it contains a parabolic subgroup $P_n$
(aka, the {\it Siegel parabolic}) whose unipotent radical  $U_n$ is isomorphic to the symmetric matrices $S^2_n$. 
Since for $Sp_{2n}$ there are no rival parabolic subgroups with abelian unipotent radicals, we will simply use $U=U_n$
for purposes of defining rank. If we have a representation $\pi$ of $Sp_{2n}$, then its restriction to $U_n$ will decompose
as a sum of characters. As described in \S 2.4, each character of $U_n$ is naturally assigned a rank. 
 We can then define the {\it $U_n$-rank} (sometime we will just say $U$-rank) of $\pi$ to be the largest of the ranks of characters that appear in $\pi_{|U_n}$. If this rank is strictly less than $n$, we will say that $\pi$ is of {\it low $U_n$-rank}.
 
For a finite field of coefficients, there are two isomorphism classes of symmetric bilinear forms of a given rank.
If it happens that, for a representation $\pi$ of $Sp_{2n}$ of $U_n$-rank $k$, that only one
of the two types of rank $k$ characters appear in $\pi_{|U_n}$, then we will say that $\pi$ is {\it of the type}
of the rank $k$ characters that appear.
 
 \vskip .05 in 
 
If $G = O_{n,n}$, the split orthogonal group in even dimensions, then we also have a parabolic subgroup whose Levi 
component is $GL_n$, with unipotent radical $U=U'_n$ isomorphic to $\Lambda^2_n$, 
the skew-symmetric $n \times n$ matrices. We can again look at the rank of characters, and their conjugacy classes 
under the action of the Levi component will consist of all characters of a given rank. We note that this rank 
will always be even. We define the {\it $U'_n$ rank} (we might just call it $U$-rank) of a representation $\pi$ 
of $O_{n,n}$ to be the largest rank of the $U'_n$ characters that appear in $\pi_{|U'_n}$. 
If this rank  is less than $n - 1$, then we say that $\pi$ is of {\it low $U$-rank}.
 \vskip .1 in
 
{\bf Case of ${\bf U}_{n,n}$}:

\vskip .05 in

This goes almost word for word as for $Sp_{2n}$ and $O_{n,n}$. We omit the details.

\vskip .1 in

{\bf Cases of non-split groups}:
\vskip .05 in

We know that, among orthogonal groups over finite fields, there are isometry groups
of forms that  are not sums of hyperbolic planes: the odd dimensional orthogonal groups
$O_{n+1, n}$, the non-split even dimensional groups $O_{n+2, n}$, as described at the end of
\S 1. Similarly, the unitary groups ${\bf U}_{n+1, n}$ in odd dimensions are not split in 
this sense.

However, each of these groups contains a maximal split subgroup, unique up to
conjugacy. The orthogonal groups $O_{n+1, n}$ and $O_{n+2,n}$ contain copies
of $O_{n,n}$. The unitary groups ${\bf U}_{n+1,n}$ contain copies of $U_{n,n}$. 
We will define $U$-rank and low $U$-rank for these groups by looking at 
restrictions to split subgroups. A representation $\pi$ of $O_{n+1, n}$ or of $O_{n+2,n}$
will be said to have $U$-rank $k$, or to be of low $U$-rank, if the restriction
$\pi_{|O_{n,n}}$ has $U$-rank $k$, or has low $U$-rank. We make the parallel definition
for representations of ${\bf U}_{n+1,n}$.

 \vskip .2 in
 
 {\bf 5. Constructing representations of small $U$-rank}
 
 \vskip .15 in
 
 Although the concept of representations of a given $U$-rank was defined in \S 4, no evidence was presented that such
representations actually exist. Indeed, since the various unipotent subgroups discussed in \S 2 and \S 3 are rather small 
subgroups of rather non-abelian groups $G$, it might be reasonable to expect that 
representations of the groups $G$ would tend to be of full $U$-rank, and that the concept of rank might give little or no information about the dual set $\widehat G$ of irreducible representations. Indeed, in some sense the ``typical" representation of $G$ does have full rank. However, in this section, we will see that $G$ also possesses natural families of representations of low $U$-rank, with interesting structure.
    
  \vskip .1 in 
  
  {\bf 5.1. The Heisenberg group and the oscillator representation.}
  
  \vskip .1 in 
  
Our method for constructing representations of given $U$-rank will use dual pairs inside
the symplectic groups $Sp_{2n}$. The notion of dual pair was defined some time ago [Howe73-2], and it has proved useful for studying representations of groups over local fields, and for questions in the theory of automorphic forms (see  [Gelbart77], [Kudla86], [Kudla-Milson86], [Prasad93], and their bibliographies). However, although there have been some applications to the representations of  groups  over finite fields [Adams-Moy93], [Aubert-Michel-Rouquier96], [Aubert-Kraskiewicz-Przebinda16], [G\'erardin77], [Srinivasan79] these have been relatively limited.
  
The context for the study of representations of dual pairs is the oscillator representation of the symplectic group, which arises through its action by automorphisms of the Heisenberg group. 
This has been discussed in considerable detail in the literature [G\'erardin77], [Gurevich-Hadani07], [Gurevich-Hadani09], [Howe73-1], [Weil64]. We will review the basic construction.
  
  We start with a symplectic vector space $W$ over the coefficient field ${\bf F}$, 
  which we are taking to be a finite field of odd characteristic. Most of the discussion, however, 
  will remain relevant for infinite local fields. The space $W$ will have additional structure in various contexts below, but for now, we just assume that $W$ has been 
  coordinatized as in \S 3.2.  The symplectic form on $W$ is denoted with sharp brackets. 
  We define the {\it Heisenberg group} attached to $W$ as the set
  $$
  H(W) \simeq W \oplus {\bf F}, \eqno (5.1.1a)
  $$
  with the group law
  $$
  (w, t)\cdot (w', t') = (w + w', t + t' + {1 \over 2} <w, w'>). \eqno (5.1.1b)
  $$
   The set
   $$
   Z(H) = \{(0, t): t \in {\bf F} \},
   $$
   is the center and the commutator subgroup of $H(W)= H$.
   
   It is evident from the formula (5.1.1b) that the symplectic group $Sp(W) = Sp_{2n}$ 
   acts on $H(W)$ by automorphisms,
   by direct extension of its action on $W$, and that the action on $Z(H)$ is trivial:
   $$
   g\cdot(w, t) = (g(w), t), \ \ \ {\rm for \ every} \ g \in Sp(W), \ w \in W, \ {\rm and} \ t \in {\bf F}. 
   $$
   
Consider an irreducible representation $\rho$ of $H(W)$. Because the center $Z(H)$
 will act via operators that commute
with all of $H$, Schur's Lemma [Curtis-Reiner62], [Serre77] implies that $Z(H)$ must act by scalar operators.
 These scalars will define a 
  character $\chi$ of $Z(H)$, referred to as the {\it central character} of $\rho$. 
  If the central character is trivial, then $\rho$ 
 factors to $H(W)/Z(H) \simeq W$, which means that $\rho$ will just be a character
  of $\rho$, and in particular will be
 of dimension 1.
 
 If the central character $\chi$ of $\rho$ is non-trivial, then it turns out that $\chi$
  determines $\rho$ 
 up to unitary equivalence. The situation is described by the following celebrated result [Mackey49]:
 
 \vskip .1 in
 
 {\bf Theorem 5.1.2 (Stone-von Neumann-Mackey Theorem).} 
 Given a non-trivial character $\chi$ of the center
 $Z(H)$ of the Heisenberg group $H(W) = H$, up to unitary equivalence, 
 there is a unique irreducible representation
 $\rho_{\chi}$ of $H(W)$ with central character $\chi$. 
 Moreover, the representation $\rho_{\chi}$
  may be realized as an induced representation
 $$
 \rho_{\chi} \simeq {\rm Ind}_{L'}^H \chi', \eqno (5.1.3)
 $$
 where $L'$ is any maximal abelian subgroup of $H$, and $\chi'$ is any extension of
  $\chi$ from $Z(H)$ to $L'$.
 \vskip .1 in
 
{\bf Remarks 5.1.4.} We note that

\vskip .05 in

{\bf a)} {\it Lagrangians}. A subgroup $L'$ of $H(W)$ that contains $Z(H)$ will be commutative if and only if the image
$L = L'/(Z(H)$ of $L'$ in $W$ is {\it isotropic}: any two vectors in $L$ have symplectic
pairing equal zero. Moreover, such $L'$  will be maximal commutative if and only if $L$ is {\it Lagrangian}, i.e., maximal isotropic subspace of $W$. 
Two examples of maximal isotropic subspaces of the space 
 $W \simeq {\bf F}^{2n}$ described in formula (3.2.1), are the space $X$, 
 consisting of the length $2n$ vectors $\vec x$ for which the last $n$ coordinates equal zero, and the space $Y$, 
 consisting of the vectors $\vec y$ for which the first $n$ coordinates equal zero.

 \vskip .05 in

{\bf b)} {\it Formulas - Schr\"odinger model}. If we take the subspace $X$ to be the subspace $L$ of the previous remark, 
 then the realization (5.1.3) of $\rho_{\chi}$ can be put in the following explicit form.
 First,  we identify $X$ and $Y$ with abelian subgroups of $H(W)$, {\underbar {not}} 
 containing $Z(H)$. Thus,
 in particular, the maximal abelian subgroup $X'$ of $H$ is the direct product 
 $X' = X \oplus Z(H)$. We 
choose $X'$ to be the subgroup $L'$, and we choose the character $\chi'$ 
to be the extension of $\chi$ to $X'$ that is trivial on $X$. We have $H(W) = X' \cdot Y$, so that we can identify the space of the induced representation ${\rm Ind}_{X'} ^H \chi'$ with $L^2(Y)$. With these agreements, we can express the representation $\rho_{\chi}$ as follows: 
$$
\bullet \\\ \rho_{\chi} (\vec y) f(\vec{y'}) = f(\vec{y'} - \vec y);  \hskip 1 in \eqno (5.1.5a) 
$$
$$
\hskip .4 in \bullet \\\ \rho_{\chi}(\vec x) f(\vec{y'}) = \chi(<\vec{y'}, \vec x>) f(\vec{y'}); \hskip 1 in \eqno (5.1.5b)
$$
$$
\hskip .13 in\bullet \\\ \rho_{\chi} (0, t) f(\vec{y'}) = \chi (t)f(\vec{y'}),  \hskip 1 in \eqno (5.1.5c)
$$
for every $ \vec y, \vec{y'} \in Y, \  \vec x \in X, \ t \in {\bf F},  \  {\rm and} \ f \in L^2(Y)$.
 
 \vskip .09 in
 
 This realization of $\rho_{\chi}$ is sometimes called the {\it Schr\"odinger model} 
 attached to the pair
 $(X, Y)$ of Lagrangian subspaces of $W$.
 \vskip .1 in
 
 The Stone-von Neumann-Mackey Theorem implies that, for any element  $g$ of $Sp(W)$, 
 the automorphism of $H(W)$
 defined by $g$, since it is the identity map on $Z(H)$,  will preserve the unitary equivalence
 class of $\rho_{\chi}$. This leads in turn to a representation of $Sp(W)$ as follows.
 
Given an element $g$ in $Sp(W)$, the Stone-von Neumann-Mackey Theorem implies
 that there is a unitary operator
$\omega (g)$ on the space of $\rho_{\chi}$, such that 
$$
\omega(g) \rho_{\chi}(h) \omega(g)^{-1} = \rho_{\chi}(g(h)), \eqno (5.1.6)
$$
 for all elements $h$ of $H(W)$. Since $\rho_{\chi}$ is irreducible, the equation (5.1.6)
  determines $\omega(g)$ up to scalar multiples. Since both $\omega(g) \omega (g')$ 
 and $\omega (gg')$ will satisfy the same equation, we conclude that
 $$
 \omega (g) \omega (g') = \beta (g, g') \omega (gg'),
 $$
 for some scalar factors $\beta (g, g')$ for any two elements of $Sp(W)$. This yields 
 what is sometimes called a 
 {\it projective representation} of $Sp(W)$. 
 
 In this context, it is natural to ask, can the scalar factors $\beta (g, g')$ be eliminated,
  to yield an actual representation 
 of $Sp(W)$?  As Weil showed in [Weil64], for local fields, they cannot be eliminated,
 but they can be reduced to $\pm 1$. 
 This provides a representation of a two-fold covering group $Mp(W)$ of $Sp(W)$, 
 known as the {\it metaplectic group}. 
 
 However, in the case of finite fields ${\bf F} = {\bf F}_q$, homological arguments [Matsumoto69], [Steinberg62], or explicit constructions [G\'erardin77], [Gurevich-Hadani07], [Gurevich-Hadani09] show that
  the scalars $\beta (g, g')$ can all be reduced to 1, and we obtain an representation of the group $Sp(W)$. 
 We will refer to this representation as the
 {\it oscillator representation}, because it is an analog of the spin representation for orthogonal groups [Brauer-Weyl35], and is also related to the quantum mechanical harmonic oscillator [Howe85]. We denote the oscillator representation by $\omega$.
 
 The space on which $\omega$ acts is by definition the space of $\rho_{\chi}$. 
 Although explicit and transparent formulas
  for $\omega (g)$ for all elements $g$ of $Sp(W)$ are not available, 
 some $g$ act in simple ways in a given realization of $\rho_{\chi}$. In particular, 
 the action of the parabolic subgroup 
 $P$ of formulas (3.2.3) can be described nicely in the Schr\"odinger model on the space $L^2(Y)$, as described in formulas (5.1.5).
 
 The subgroup $M$ of formula (3.2.3) restricts to $Y$ as $GL(Y)$. Suppose we have
  an element of $M_n$, 
 as described in equation (3.2.3b) .  Then we can check that the operators
 $$
 \omega \left ( \left [ \matrix {A \ \ \ \ \ 0\ \ \ \cr \  0 \ \ (A^{-1})^t} \right ] \right )  (f) (\vec y) 
 = \left(\frac{\det(A)}{q}\right) f (A^t \vec y), \eqno (5.1.7)
 $$
$f \in L^2(Y),\, \vec{y} \in Y$, satisfy the defining equation (5.1.6) for $\omega$. Here, it turns out the specific value of the Legendre symbol\footnote{The value of $\left(\frac{x}{q}\right)$ for $x \in {\bf F}_q^\times$ is $1$ or $-1$ if $x$ is a square or not, respectively, in ${\bf F}_q^\times$.}
$\left(\frac{\det(A)}{q}\right)$ is the correct normalization [G\'erardin77], [Gurevich-Hadani07], [Gurevich-Hadani09].
 
 The elements of the subgroup $U_n$, described in equation (3.2.3c), also have a simple
  description: they act as
  multiplication operators. From the definitions, one calculates [G\'erardin77], [Gurevich-Hadani07], [Gurevich-Hadani09] that  
  $$
  \omega \left ( \left [  \matrix {I \ \ \ S \cr 0 \ \ \ I} \right ] \right ) (f)(\vec y) = 
  \chi \left ( {B_S(\vec y, \vec y) \over 2} \right ) (f)(\vec y), \eqno (5.1.8)
  $$
  $f \in L^2(Y),\, \vec{y} \in Y.$
 \vskip .1 in
 
 {\bf 5.2. Dual pairs in $Sp_{2n}$.}
 \vskip .1 in
 
 The oscillator representation is a very special representation of $Sp_{2n}$. 
 It turns out, however, that it can be used 
 to construct or describe large families of representations of all classical groups [Gurevich-Howe15]. 
 We give some examples to show how this can be done. However, we will not attempt to treat all possible examples. The examples that we will not describe here can be dealt with in the same manner and will appear elsewhere.  
  
  It turns out to be instructive to restrict the oscillator representation to subgroups of $Sp_{2n}$. 
 However, such subgroups 
  must be chosen with care. The restriction of a representation to a subgroup can be very hard
 to analyze explicitly. When the base field is a local field, it has been found that restricting the oscillator representation to subgroups of the form $G \times G'$, where each of $G$ and $G'$ is the full centralizer of the other in $Sp_{2n}$, yields interesting correspondences (the {\it local theta correspondence}) between representations of $G$ and representations of $G'$. 
  
The ideas that were used to create the local theta correspondence do not apply when the base field is finite. However,
it turns out that the notion of $U$-rank discussed in \S 4 provides a satisfying substitute, which also sheds new light on
the situation for local fields.
 
 We call pair of subgroups $(G, G')$ of a group $\Gamma$ such that each of $G$ 
 and $G'$ is the full centralizer
  of the other, a {\it dual pair} of subgroups of $\Gamma$.   
  
  \vskip .1 in
  
  {\bf Examples 5.2.1.} The following pairs of groups form dual pairs in suitable symplectic groups.

  \vskip .05 in
  
  {\bf a)} {\it The $(Sp_{2n}, O_k)$ dual pair}. If  a vector space $W$ is a equipped with a symplectic form 
  $< \ , \ >_{W}$, and 
  a vector space $V$ is equipped with a symmetric bilinear form $B$, then the tensor product
   $W \otimes V$ can be endowed with the tensor product of $< \ , \ >_{W}$ and $B$, which is a symplectic form. Thus 
  $(Sp(W), O(V))$  form a dual pair in  $Sp(W \otimes V)$. 
  If $\dim W = 2n$ and $\dim V = k$, then $\dim (W \otimes V) = 2nk$, so we may simply say that 
  $$
  (Sp_{2n}, O_k) \ {\rm form \ a \ dual \ pair \ in} \ Sp_{2nk}.
  $$
  
 \vskip .05 in 
 
{\bf b)} {\it The $(GL_k, GL_n)$ dual pair}. The groups $GL_k$ and $GL_n$ acting on the matrices $M_{k,n}$ by multiplication
  on the left and on the right, as described in formulas $(2.1.1)$, then $(GL_k, GL_n)$ form a dual pair
  in $GL_{kn}$. Moreover, it is well understood [Lang02] that, given any vector space $V$, the natural pairing between
  $V$ and its dual $V^*$ can be used to define a symplectic form on $W = V \oplus V^*$, by the recipe
 $$
 <(\vec v, \lambda), (\vec v', \lambda')> = \lambda'(\vec v) - \lambda(\vec v').
 $$
 Then the action of $GL(V)$ on $W$ defined by $$
 g\cdot (\vec v, \lambda) = (g(\vec v), \lambda\circ g^{-1}),
 $$
 embeds $GL(V)$ into $Sp(V \oplus V^*)$. Taking $V = M_{k,n}$ then embeds 
 $$
 (GL_k, GL_n) \ {\rm as \ a \ dual \ pair \ in} \ Sp(M_{k,n} \oplus M_{n,k}) \simeq Sp_{2kn}.
 $$
 
 \vskip .1 in 
 
 {\bf 5.3. Constructing representations of low $U$-rank.}
 \vskip .1 in
 
 The dual pairs described in \S 5.2 can be used to construct representations of low 
 $U$-rank for the groups
 $Sp_{2n}$, $O_{n,n}$, and $GL_n$. 
 
 Recall the {\it isotypic decomposition} of a representation.  Given any representation
  $\pi$ of a finite group $G$ on a vector space $V$, 
 and an irreducible representation $\tau$ of $G$, the {\it $\tau$-isotypic component 
 $\pi_{\tau}$} of $\pi$ is
 the sum of all the subspaces of $V$ on which $G$ acts irreducibly, by a representation 
 that is unitarily equivalent to 
 $\tau$ [Serre77]. It is well known that $V_{\tau}$ has the structure of a tensor product
 $$
 V_{\tau} \simeq V' \otimes {\tau},
 $$
 where $\tau$ here is standing for a vector space on which $G$ acts irreducibly
  by the representation $\tau$, 
 and $V'$ is an auxiliary space. 
 
 If $G'$ is a group acting on the same space
  as $G$, and commuting with $G$, then $V'$ inherits the structure of $G'$-module,
  and we can write
  $$
  V' \simeq \sum_{\sigma \in \widehat G'} V'_{\sigma},
  $$
  where   $V'_{\sigma}$ is the $\sigma$-isotypic component of $V'$ for $G'$. 
  Taking the tensor product with $\tau$
  gives us a decomposition
  $$
  V_{\tau} \simeq \sum_{\sigma \in \widehat G'} V'_{\sigma} \otimes \tau. 
  $$
  The spaces $V'_{\sigma} \otimes \tau$ are the $G' \times G$ isotypic components
   for $V_{\tau}$.
 
 \vskip .1 in
 
 {\bf Theorem 5.3.1 (The eta correspondence).} We have
 
 \vskip .06 in
 
 {\bf a)} Let $k \leq n$ be whole numbers. Consider the oscillator
  representation $\omega$ for the dual pair
 $(Sp_{2n}, O_k) \subset Sp_{2nk}$, where $O_k$ is the isometry group of either
  of the two symmetric non-degenerate bilinear forms
  in dimension $k$. Then
  
  \vskip .04 in 
  \hskip .3 in {\bf 1)} The representation $\omega_{|Sp_{2n}}$ is of $U$-rank $k$, and has type 
  corresponding to $O_k$. 
  
  \vskip .04 in 
  \hskip .3 in {\bf 2)} For any irreducible representation $\tau$ of $O_k$, let 
 $\omega_{\tau} \simeq \Omega_{\tau} \otimes \tau $ be the $\tau$-isotypic subspace of 
 $\omega$. Then, as a
  representation of $Sp_{2n}$, $\Omega_{\tau}$ contains a unique irreducible representation
   $\eta(\tau)$ of $U$-rank $k$ (and type of $O_k$). Furthermore, the mapping
 $$
 \eta: \widehat {O}_k  \hookrightarrow  \widehat {Sp}_{2n}, \eqno (5.3.2a)
 $$
 is an injection  of $\widehat {O}_k$ into the set $(\widehat{Sp}_{2n})_{U,k}$ of irreducible $U$-rank $k$ representations of $Sp_{2n}$.
 
\vskip .06 in

{\bf b)} Let $k$ and $n$ be whole numbers, with $2n \leq k$. Consider the oscillator 
representation $\omega$ for the dual 
pair $O_{k,k} \times Sp_{2n} \subset Sp_{4nk}$. Then 

\vskip .04 in 
\hskip .3 in {\bf 1)} The representation $\omega_{|O_{k,k}}$ has $U$-rank $2n$.

\vskip .04 in 
\hskip .3 in {\bf 2)} For any irreducible representation $\sigma$ of $Sp_{2n}$, let 
$\omega_{\sigma} \simeq \Upsilon_{\sigma} \otimes \sigma$ be the $\sigma$-isotypic
 component of $\omega$. Then,
as a representation of $O_{k,k}$, $\Upsilon_{\sigma}$ contains a unique irreducible
 representation $\eta(\sigma)$ of rank $2n$. Furthermore, the mapping
 $$
 \eta: {\widehat {Sp}}_{2n} \hookrightarrow  {\widehat O}_{k,k}, \eqno (5.3.2b)
 $$
 is an injection of ${\widehat {Sp}}_{2n}$ into the set $({\widehat O}_{k,k})_{U,2n}$ 
 of irreducible $U$-rank $2n$ representations of $O_{k,k}$.
 
 \vskip .06 in 
 
 {\bf c)} Let $k$ and $n$ be whole numbers, with $k \leq {n \over 2}$. Consider the oscillator
  representation for the dual pair  $GL_k \times GL_n \subset Sp_{2kn}$. Then 
  
\vskip .04 in 
\hskip .3 in {\bf 1)} The representation $\omega_{|GL_n}$ has $U$-rank $k$. 

\vskip .04 in 
\hskip .3 in {\bf 2)} For any irreducible representation $\nu$ in ${\widehat {GL}}_k$, let $\omega_{\nu} \simeq \Xi_{\nu} \otimes \nu$ be the 
 $\nu$-isotypic component of $\omega$. Then, as a representation of $GL_n$, $\Xi_{\nu}$ contains a unique irreducible
 representation $\eta (\nu)$ of rank $k$. Furthermore, the mapping
 $$
 \eta: {\widehat {GL}}_{k} \hookrightarrow  {\widehat {GL}}_{n} , \eqno (5.3.2c)
 $$
 is an injection of ${\widehat {GL}}_{k} $ into the set $({\widehat{GL}}_{n})_{U,k}$ of irreducible 
 $U$-rank $k$ representations of $GL_n$.

\vskip .15 in

For what follows, by a {\it twist} of a representation $\pi$ of a group $G$ by a character $\chi$ of $G$ (i.e., a one dimensional representation) we mean the representation $\chi\otimes\pi$. Extensive numerics and theoretical considerations lead us to formulate the following natural conjecture: 

\vskip .15 in

{\bf Conjecture 5.3.3 (Exhaustion of the $\eta$-correspondence).} We have

\vskip .06 in

{\bf a)} Every irreducible representation of $Sp_{2n}$ of $U$-rank $k$, for $k<n$, appears in the image of the mapping (5.3.2a) for (exactly) one of the two possible orthogonal groups $O_k$.

\vskip .06 in

{\bf b)} Every irreducible representation of $O_{k,k}$ of $U$-rank $2n$, for $2n<k$, appears, maybe after a twist by a character, (exactly once) in the image of the mapping (5.3.2b).

\vskip .06 in

{\bf c)} Every irreducible representation of $GL_n$ of $U$-rank $k<\frac{n-1}{2}$ appears, maybe after a twist by a character, (exactly once) in the image of the mapping (5.3.2c).

\vskip .1 in
 
{\bf A proof of Theorem 5.3.1.}

 \vskip .06 in
 
A proof of part a) of Theorem 5.3.1 is given in [Gurevich-Howe15]. A proof of part b) can be given along
 very similar lines. Here we will offer a proof of part c), which involves using a somewhat non-standard Schr\"odinger model.
 \vskip .1 in
{\it Proof.} As explained in \S 5.1, we can construct a realization of the oscillator representation
  of $Sp(W)$ for a symplectic vector
  space $W$ starting from a supplementary pair $(X, Y)$ of maximal isotropic subspaces
 of $W$. In Example 5.2.1, we 
 realize the pair $(GL_n, GL_k)$ as a dual pair in $Sp_{2nk}$ by considering the actions
  by left and right multiplication
 of $GL_k \times GL_n$ on the sum of matrix spaces $W = M_{k, n} \oplus M_{n, k}$, 
 and the pair of spaces 
 $(M_{k, n}, M_{n,k})$ constitute a pair of maximal isotropic subspaces of $W$. 
 The corresponding realization of 
 the representation $\rho_{\chi}$ of $H(W)$ would be on $L^2(M_{k,n})$, 
 and $GL_k \times GL_n$ would act by means of their permutation actions given by left and right multiplication 
 of matrices. 
 
 This is certainly a pleasant realization of $\omega_{|GL_k \times GL_n}$, but it does not display the spectral 
 behavior of the standard subgroups of $GL_n$ in as transparent a way as possible. 
 We will use a different pair of maximal isotropic subspaces that will do this better.
 
 To this end, we replace the specific spaces ${\bf F}_q^n$ and ${\bf F}_q^k$ with vector spaces
  $X$ and $Y$ over ${\bf F}_q$,
  of dimensions $n$ and $k$ respectively.
We construct the symplectic vector space 
$$
W = Hom (Y, X) \oplus Hom (X, Y), \eqno (5.3.4)
$$
equipped with the symplectic form
$$
<(A, B), (A', B')> = {\rm trace}(B'A - BA'). 
$$ 
The oscillator representation of $Sp(W)$ is then realized on space of the canonical
representation of the 
Heisenberg group $H(W)$, as described in \S 5.1.
To produce a concrete space for this representation we can choose a pair of supplementary
 Lagrangian subspaces
 of $W$. 

We will not use the pair of formula (5.3.4).
For our realization of $\omega$, we first choose a decomposition $X = X_1 \oplus X_2$
 of $X$, with $\dim X_1 = k$, 
and therefore
$\dim X_2 = n - k$. Then
$$
W = (Hom (Y, X_1) \oplus Hom (X_2, Y)) \oplus (Hom (X_1, Y) \oplus Hom (Y, X_2)) 
= L \oplus L' , \eqno (5.3.5)
$$
expresses our $W$ as a sum of two Lagrangian subspaces $L$ and $L'$.  
The Schr\"odinger model
attached to this decomposition is realized on $L^2(L')$. We will look at the action
 of $GL(X) \times GL(Y)$ 
in this model.

 Since the
groups $GL(Y)$ and $GL(X_1) \times GL(X_2) \subset GL(X)$ both preserve both Lagrangians of the decomposition
(5.3.5), they will act on   $L^2(L')$ by their permutation actions
on the points of the set $L' = Hom (X_1, Y) \oplus Hom (Y, X_2)$.
Adapting to the present context from formulas (5.1.7), we find that up to a $\pm$ sign (that will not matter for the discussion below)
$$
\omega (g') (f)(A , B) = f((g')^{-1}A, Bg'),
$$
for $g'$ in $GL(Y)$, $A$ in $Hom (X_1, Y)$ and $B$ in $Hom (Y, X_2)$.
The diagonal subgroups $GL(X_j)$ also act (up to a $\pm$ sign) via their obvious permutation actions, 
adapted from formula (5.1.7):
$$
\omega (g_1 \times g_2)(f)(A, B) = f(Ag_1, (g_2)^{-1}B),
$$
for $(g_1, g_2)$ in $GL(X_1) \times GL(X_2)$.

As a subgroup of $Sp(W)$, the standard subgroup
 $U = Hom(X_2, X_1) \simeq M_{k, n - k} \subset GL(X)$
 acts as the identity on  $L$, and also acts as the identity on
 the quotient 
 
\centerline { $W/L \simeq L'$,} 
 \noindent so it will act on  $L^2(L')$
 by multiplication operators, by an appropriate adaptation of formula (5.1.8).
 
We can define a mapping
$$
\mu: L' \simeq Hom (X_1,Y) \oplus Hom (Y,X_2) \longrightarrow Hom (X_1,X_2)$$
by taking the products of the two matrix components:
$$
\mu(A, B) = BA,
$$
with $A$ in $Hom (X_1, Y)$ and $B$ in $Hom (Y, X_2)$, as before.  Using this notation, 
the action of $T$ in $U$ on 
$L^2(Hom (X_1, Y) \oplus Hom (Y, X_2))$ may be described by
$$
\omega (T)(f)(A, B) = \gamma_{\mu(A, B)} (T) f(A, B). \eqno (5.3.6)
$$
Here $\gamma_{\mu(A,B)}$ is the character of $U$ described by formula (2.4.2).

We want to use formula (5.3.6) to describe the eigenspaces of $U$, and especially
 the rank $k$ spectrum 
of $U$. To do this, we study the structure of the set of
pairs $(A, B)$ such that $\mu(A, B) = BA$ has rank $k$.
\vskip .05 in 
The product $BA$, as an element of $Hom(X_1, X_2)$, has rank at most equal 
to the maximum of the ranks of $A$ 
and $B$. Since both $A$ and $B$ are maps to or from $Y$, neither of them can have
 rank larger than $k = \dim Y$. 
\vskip .05 in 
Let's consider the collection of $(A, B)$ such that rank($BA) = k$. If $BA$ has rank $k$, 
then $A$ must have rank 
at least $k$. Since we have selected $X_1$ so that $\dim X_1 = \dim Y = k$, 
the matrix $A$ has rank at most $k$,
and if it does have rank $k$ then it is an isomorphism
from $X_1$ to $Y$, and we can conclude that rank$(BA)$ = rank $(B)$. For $B$ 
to have rank $k$ means that $B$
is injective. If we have two injective mappings $B_1$ and $B_2$ from $Y$ to $X_2$, 
then they differ by an element of 
$GL(X_2)$. Thus, if rank($B_2$) = rank ($B_1$) = $k$, then
$$
B_2 = g_2 B_1,
$$
for some $g_2$ in $GL(X_2)$. 

Similarly, if $A: X_1 \rightarrow Y$ has rank $k$, it is bijective. Given two bijective mappings
$A_1$ and $A_2$ from $X_1$ to $Y$, there is an element $g_1$ of $GL(X_1)$ such that 
$$
A_2 = A_1 g_1.
$$
Putting these two observations together, we conclude that $GL(X_1) \times GL(X_2)$ acts transitively on 
the set of $(A, B)$ such that $BA$ has rank $k$. 
\vskip .05 in

We next look at the set $\mu^{-1}(BA)$ of points $(A', B')$ such that 
$\mu(A', B') = B'A' = BA = \mu(A, B)$ is a fixed mapping of rank $k$.

If $\mu (A', B') = B' A' = B A = \mu(A, B)$, with both having rank $k$, then we know that $A$
 and
$A'$ both have rank $k$, and are therefore both isomorphisms from $X_1$ to $Y$. 
This implies that there is an element $g'$ in $GL(Y)$ such that $A' = g' A$. 
Since the products $BA$ and $B'A'$ are equal, this means that
$B' = B(g')^{-1}$. Conversely, it is clear that $(A, B)$ and $(A', B') = (g'A, B(g')^{-1})$ satisfy
$B'A' = BA$. We conclude that $GL(Y)$ acts transitively on each set $\mu^{-1}(BA)$ 
such that the product
$BA$ is a fixed rank $k$ element of $Hom (X_1, X_2)$. In fact, this action is simply transitive, 
since if $A$ has rank $k$, 
it is surjective, and $g'A = A$ only when $g' = I_Y$ is the identity map of $Y$.

Let's look also at the stabilizer of $BA$ in $GL(X_1) \times GL(X_2)$. The image space
 $BA(X_1)$ is a $k$-dimensional
subspace of $X_2$. The stabilizer in $GL(X_2)$ of $BA(X_1)$ is a maximal parabolic
 subgroup of $GL(X_2)$, which we will denote $P_{BA}$. The restriction of an element
  in $P_{BA}$ to $BA(X_1)$ defines a surjective homomorphism
$$
r_{BA}: P_{BA} \longrightarrow GL(BA(X_1)).
$$
The kernel $Q_{BA}$ of $r_{BA}$ is of course a normal subgroup of $P_{BA}$. 
We can get a complement to $Q_{BA}$ in 
$P_{BA}$ by choosing a vector space complement $X_3$ to $BA(X_1)$ in $X_2$. 
This gives us an isomorphism
$$
P_{BA} \simeq GL(BA(X_1)) \cdot Q_{BA} \simeq GL(BA(X_1)) \cdot GL(X_3) \cdot U_{BA},
$$
where $U_{BA}$ is the unipotent radical of $P_{BA}$, and is isomorphic to 
$Hom (X_3, BA(X_1))$.

It is clear that the stabilizer of $BA$ in $GL(X_1) \times GL(X_2)$ contains $Q_{BA}$. 
It also contains a subgroup
$\Delta_{BA} = \Delta(GL(X_1) \times GL(BA(X_1)) \subset GL(X_1) \times GL(BA(X_1))$,
 defined by
$$
\Delta_{BA} = (g_1, BAg_1^{-1}(BA)^{-1}), \eqno (5.3.7)
$$
where $g_1$ is in $GL(X_1)$, and we interpret $(BA)^{-1}$ as being the inverse
 of $BA$ considered as
 an isomorphism between $X_1$ and $BA(X_1)$.
 
 Now fix a pair $(A, B)$ in $L'$, and consider the set $\mu^{-1} (BA)$ of points 
 $(A',B')$ in $L'$ 
 such that $\mu(A',B') = B'A' = BA = \mu(A,B)$. We have seen above that the group 
 $GL(Y)$ acts simply transitively on
 $\mu^{-1}(BA)$. It is also easy to check that the group $Q_{BA}$ acts trivially on
  $\mu^{-1}(BA)$. Since $GL(X_1)$
 acts simply transitively on the set of $A$ in $Hom (X_1, Y)$ of rank $k$, the group 
 $\Delta_{BA}$ also will act simply transitively on $\mu^{-1}(BA)$. 
 
 \vskip .05 in

Thus, we are in the following situation. The delta functions at the points of $L'$ 
are all eigenfunctions
for the standard subgroup $U \simeq Hom (X_2, X_1)$ of $GL(X)$. If we fix a character
 $\gamma_{\mu(A,B)}$
of rank $k$ of $U$,
then $GL(Y)$ acts simply transitively on the set of points $(A',B')$ in $L'$ whose
 associated eigencharacter
$\gamma_{\mu(A',B')}$ is the same as $\gamma_{\mu(A,B)}$. Thus, the  action of $GL(Y)$
 on the 
 $\gamma_{\mu(A,B)}$ eigenspace for $U$ is equivalent to right regular representation of 
 $GL(Y)$.
 
 On the other hand, the subgroup $\Delta_{BA}$ of the stabilizer of $BA$ in 
 $GL(X_1) \times GL(X_2)$ also
 acts simply transitively on $\mu^{-1}(BA)$, and this action of course commutes with the action
  of $GL(Y)$. It follows that, if we identify $\mu^{-1}(BA) \simeq GL(Y_1)$ by selecting $(A, B)$ as a base point, then the 
 action of $\Delta_{BA}$ gets identified to the left regular action of $GL(Y)$ on itself.
 We may conclude from the Peter-Weyl Theorem for finite groups [Serre77], [Curtis-Reiner62], that the actions of $GL(Y)$
 and of $\Delta_{BA}$ on
 $L^2(\mu^{-1}(BA))$ generate mutual commutants in the algebra of operators on 
 $L^2(\mu^{-1}(BA))$.
 
 \vskip .05 in
 
 We can use this to prove Theorem 5.3.1 c). Fix an irreducible representation $\nu$ 
 of $GL(Y)$, and consider 
 the $\nu$-isotypic component $\omega_{\nu} \simeq \Xi_{\nu} \otimes \nu$ 
 of the oscillator representation 
 for the dual pair $GL(Y) \times GL(X) \subset Sp(Hom (Y, X) \oplus Hom (X, Y))$.
  Using the decomposition 
 $X = X_1 \oplus X_2$ described above (formula (5.3.5)),  consider the action of the subgroup
 $U \simeq Hom (X_2, X_1)$ on $\omega$. Fix a character $\gamma_{BA}$ of $U$ of rank $k$, and 
 consider the $\omega_{\gamma_{BA}}$ eigenspace of $\omega$. We have
 $$
 (\omega_{\nu})_{\gamma_{BA}} \ \ \ \simeq \ \ \ \ (\Xi_{\nu})_{\gamma_{BA}} \otimes \nu. 
 $$

 Let $\Xi_{\nu} \simeq \sum_{\zeta \in \widehat{GL}(X)} m_{\zeta} \zeta$ be
  a decomposition of $\Xi_{\nu}$ into
 irreducible representations for $GL(X)$, with multiplicities $m_{\zeta}$. Then 
 $$
 (\Xi_{\nu})_{\gamma_{BA}} \simeq 
 \sum_{\zeta \in \widehat{GL}(X)} m_{\zeta} \zeta_{\gamma_{BA}},
 $$
 and so
 $$
 (\omega_{\nu})_{\gamma_{BA}} \ \ \ \simeq \ \ \ \ 
 ( \sum_{\zeta \in \widehat{GL}(X)} m_{\zeta} \zeta_{\gamma_{BA}})\otimes \nu.
  $$
 
 For any representation $\zeta$ of $GL(X)$, the eigenspace $\zeta_{\gamma_{BA}}$
  will be invariant under the group
 $\Delta_{BA}$ of formula (5.3.7), and so the subspace
 $(m_{\zeta} \zeta_{\gamma_{BA}})\otimes \nu$ of $\omega_{\nu}$ will be invariant under
 $\Delta_{BA} \times GL(Y)$. If the sum of the multiplicities $m_{\zeta}$ for which
  $\zeta_{\gamma_{BA}} \neq \{ 0 \}$
  is more than one, then there will be more than one $\Delta_{BA} \times GL(Y)$
 invariant subspace of 
 $\omega_{\nu}$. However, our description above of $\omega_{\gamma_{BA}}$ 
 shows that $\Delta_{BA}$ and
 $GL(Y)$ each generate the commutant of the other group acting on this space. 
 This implies that the product group $\Delta_{BA} \times GL(Y)$ acts irreducibly on the space
  $(\omega_{\gamma_{BA}})_{\nu} = (\omega_{\nu})_{\gamma_{BA}}$. It follows that there 
  must be a unique $\xi$ with
  $\xi_{\gamma_{BA}} \neq \{ 0 \}$, and that for this $\xi$, the multiplicity $m_{\xi} = 1$. 
This is the assertion of  Theorem 5.3.1 Part c). \Qed 

 \vskip .2 in
 
{\bf 6. Extensions of the notion of $U$-rank}
  
 \vskip .15 in
 
 An important feature of $U$-rank is that it behaves well under restriction to certain subgroups,
 and under tensor product. These stability properties allow one to define notions
  of rank that make sense beyond the range where $U$-rank makes sense.

 \vskip .1 in
 
 {\bf 6.1. Stability of $U$-rank under restriction.}

 \vskip .1 in

 For $m \leq n$, consider $GL_m$ embedded in $GL_n$ as $m \times m$ block matrices:
 $$
 \left [ \matrix {A  \ \ \ 0 \cr 0 \ \ \ I} \right ],
 $$
 where $I$ here denotes the $(n -m) \times (n - m)$ identity matrix.  
 Given a representation $\pi$ of $GL_n$,
 consider the restriction $\pi_{|GL_m}$ of $\pi$ to $GL_m$. We can ask about
  the $U$-rank of $\pi_{|GL_m}$, 
 and in particular, whether it has low $U$-rank. The following stability result is implicit
  in our analysis of rank in \S 4.
 \vskip .1 in
 
 {\bf Proposition 6.1.1.} If a representation $\pi$ of $GL_n$ has low $U$-rank
  equal to $k$, i.e. $k < {n-1 \over 2}$,
 then if $k < {m-1 \over 2}$, the restricted representation $\pi_{|GL_m}$ also has low
  $U$-rank $k$ with respect to $GL_m$.
  
 \vskip .1 in
 
 {\it Proof}. Indeed, according to Step 2 in the proof of Proposition 4.2.2, a representation 
 $\pi$ of $GL_n$ will have
  $U$-rank $k$, with $k < {n-1 \over 2}$, if and only if the restriction of $\pi$
   to a standard subgroup of size 
  $(k+1) \times (k+1)$ has
 rank $k$. But given the assumed restrictions on $k$, the subgroup $GL_m \subset GL_n$
  will contain standard subgroups
$U$ of this size, and these are also standard subgroups when considered as subgroups
 of $GL_n$. Since the rank of 
$\pi_{|U}$ is the same, regardless of whether $U$ is considered as a subgroup of $GL_m$
 or of $GL_n$, the result follows from Proposition 4.2.2. \Qed

\vskip .1 in

Similar considerations imply parallel results for the groups $Sp_{2n}$ or $O_{n,n}$. To state them 
we need to think of smaller groups of a given type as subgroups of larger ones. This can easily be done, 
in a manner analogous to our embedding of $GL_m$ in $GL_n$. Specifically, with $Sp_{2n}$ realized as described 
in \S 3.2, we define the standard embedding of $Sp_{2m} \hookrightarrow Sp_{2n}$ to be the group of $2n \times 2n$ 
matrices of the form
$$
\left [ \matrix {A \ \ \ 0 \ \ \ B \ \ \ 0 \cr 0 \ \ \ I \ \ \ 0 \ \ \ 0 \cr
C \ \ \ 0 \ \ \ D \ \ \ 0 \cr 0 \ \ \ 0 \ \ \ 0 \ \ \ I} \right ], \eqno (6.1.2)
$$
where $I$ denote the $(n -m) \times (n -m)$ identity matrix, and
$$
\left [ \matrix { A \ \ \ B \cr C \ \ \ D} \right ]
$$
belongs to $Sp_{2m}$.

With $O_{n,n}$ described as in \S 3.3, we embed $O_{m,m}$ in $O_{n,n}$ 
in the analogous way.

\vskip .1 in

{\bf Proposition 6.1.3.} We have, 

\vskip .06 in

{\bf a)} Let $\pi$ be a representation of $Sp_{2n}$, of low $U$-rank
 $k < m \leq n$. Let $Sp_{2m}$
be embedded in $Sp_{2n}$ by the standard embedding (6.1.2). Then 
the restriction $\pi_{|Sp_{2m}}$ also has low $U$-rank $k$. Conversely, 
if $\pi$ is a representation of $Sp_{2n}$ 
whose restriction to $Sp_{2m}$ has low $U$-rank $k$, then $\pi$ also has 
$U$-rank $k$ as a representation of $Sp_{2n}$.

\vskip .06 in

{\bf b)} Let $\pi$ be a representation of $O_{n,n}$, of low $U$-rank $k < m \leq n$. Let $O_{m,m}$
be embedded in $O_{n,n}$ by the standard embedding analogous to (6.1.2). Then 
the restriction $\pi_{|O_{m,m}}$ also has low $U$-rank $k$. Conversely, if $\pi$ is a representation 
of $O_{n,n}$ whose restriction to $O_{m,m}$ has low $U$-rank $k$, then $\pi$ also has $U$-rank $k$
as a representation of $O_{n,n}$.

\vskip .1 in

{\it Proof}. Although there was no need to introduce such considerations for purposes
 of defining the $U$-rank 
of a representation of $Sp_{2n}$ or $O_{n,n}$, the same ideas that led to the conclusion, 
that the $U$-rank
of a representation of $GL_n$ is determined by its restriction to standard subgroups
 of size $m \times m$,
for any $m$ between $n$ and the $U$-rank of the representation, 
apply to the present cases also, and imply the stated results. 

In particular, the inequalities between the rank of a symmetric or skew symmetric
matrix and its components under a direct sum decomposition 
(see formulas (2.2.5) through (2.2.7)) allow us to make
the following conclusion. If we have a representation $\rho$ of $U_n \simeq S^2_n$ 
that is invariant under the natural action
$\beta$ (see formulas (2.2.3) and (2.2.4)) of $GL_n$ on $S^2_n$, 
then if a character of rank $k$ appears
in $\rho$, a character of rank $\min(k, m)$ will appear in $\rho_{|U_m}$, where $U_m \simeq S^2_m$ is the upper 
triangular unipotent radical (the Siegel unipotent) for $Sp_{2m}$. Thus, 
if $\rho = \pi_{|U_n}$ is the restriction of a representation
of $Sp_{2n}$, of low $U$-rank $k < m$, then also $\pi_{|Sp_{2m}}$ will have low $U$-rank $k < m$. A parallel argument applies likewise to the orthogonal groups $O_{n,n}$. \Qed

\vskip .1 in

{\bf 6.2. Asymptotic rank.}

\vskip .1 in

 The propositions 6.1.1 and 6.1.3 permit making a definition that extends the range of $k$
 for which one can define the rank of a representation for a given group. 

  Let $G_n$ denote one of the groups $GL_n$, $Sp_{2n}$ or $O_{n,n}$. 
 In the definition below, a group $G_m$ should be taken to denote a group
  of the same type as $G_n$. 
 That is, if $G_n = Sp_{2n}$, then also $G_m = Sp_{2m}$; 
 and likewise for the other two types of groups.
 
 \vskip .1 in 
 
{\bf Definition 6.2.1 (Asymptotic rank)}.  A representation $\pi$ of $G_n$ will be said to have 
{\it asymptotic rank less than or equal to $k$}, if for all sufficiently large numbers $m$, there is  
a representation $\rho$ of $G_m$, such that 

\vskip .05 in 

{\bf i)} $\rho$ has low $U$-rank less than or equal to $k$; and

\vskip .05 in  

{\bf ii)} $\pi \, < \, \rho_{|G_n}$. 

\vskip .12 in
 
The asymptotic rank of a representation will be sometime denote by $A$-rank or $rank_A$.
 
\vskip .1 in

{\bf Remarks 6.2.2.} Note that

\vskip .06 in

{\bf a)} We may conclude  from Proposition 6.1.3 that, if $\pi$ is of $U$-rank $k$ or less,  
with $k$ in the range where $U$-rank makes sense for the given group $G _n$, then the 
asymptotic rank of $\pi$ is greater than or equal to the $U$-rank of $\pi$. 

\vskip .06 in

{\bf b)} The form of the definition is somewhat awkward. This is because we must entertain 
the possibility that the representation of $G_m$ that is needed to find a given low rank 
representation of $G_n$ has strictly larger $U$-rank than  the original representation. 
It is even conceivable {\underbar a priori} that the asymptotic rank of some representations of $G$,
even representations of low $U$-rank, could be infinite. We will see shortly  (see inequality (6.5.1)) that this is not the case. We do not know of an example of a representation with small $U$-rank for which the asymptotic rank is larger than the $U$-rank. 

\vskip .06 in

{\bf c)} Because the restriction of representations from groups to subgroups rarely 
preserves irreducibility,
it is to be expected that the restriction of  even an irreducible representation 
$\rho$ from $GL_m$ 
with $m$ large to $GL_n$ may well contain many irreducible components, 
of varying ranks. This consideration
entails the ``less than or equal" terminology of the definition.

\vskip .06 in

{\bf d)} From the structure of Definition 6.2.1, we see that if a representation $\pi$ of $G_n$ has asymptotic rank less than or equal to $k$, then it has asymptotic rank
less than or equal to $k+1$, and so on.  Thus, the subsets $(\widehat{G}_n)_{A,\leq k}$ 
of irreducible representations of asymptotic 
rank $\leq k$ increase with $k$; so they define a filtration of $\widehat {G}_n$. We will see later that
there is a number $k_{max}$ such that $(\widehat {G}_n)_{A,\leq k_{max}} = \widehat{G}_n$. 
(In fact, $k_{max} = n$ for $GL_n, \, n>1$ (for $GL_1,\, k_{max} = 0$),
and $k_{max} = 2n$ for $Sp_{2n}$ or $O_{n,n}$).
The main point of asymptotic rank is that it provides an analog
 of the filtration given by $U$-rank that extends beyond the obvious maximum values 
 of $U$-rank. Thus, it essentially  provides a refinement of $U$-rank to values higher
 than those for which $U$-rank makes sense.

\vskip .1 in

{\bf 6.3. Rank and tensor product.}
\vskip .1 in

Tensor product is a natural operation on representations of a group $G$. For abelian groups, 
tensor product of characters is the addition operation in the Pontrjagin dual group. Because of 
the relation between addition of matrices and rank, tensor product of representations 
behaves well with respect to the notions of rank defined above.

Again, we let $G_n$ stands for one of the groups $GL_n$, $Sp_{2n}$ or $O_{n,n}$.

\vskip .1 in

{\bf Proposition 6.3.1 (Compatibility with tensor product).} We have
\vskip .06 in

{\bf a)} {\it Compatibility of tensor product and $U$-rank}. Let $\pi$ and $\rho$ be representations of $G_n$, 
of low $U$-ranks $k$ and $\ell$ 
respectively, Suppose that $k + \ell$ is less than the upper bound for low 
$U$-rank of $G_n$ (i.e., ${n \over 2} - 1$
for $GL_n$, or $n - 1$ for $Sp_{2n}$ or $O_{n,n}$). Then the tensor product 
$\pi \otimes \rho$ has low $U$-rank $k  + \ell$.

\vskip .06 in

{\bf b)} {\it Compatibility of tensor product and asymptotic rank}. Let $\pi$ and $\rho$ be representations of $G_n$, of asymptotic ranks less than or equal to $k$ and $\ell$ respectively. 
Then the tensor product $\pi \otimes \rho$ has asymptotic rank less than or equal to $k + \ell$.
\vskip .1 in

{\it Proof}. To simplify the argument, let's suppose that $G_n = GL_n$. The other cases are dealt with 
similarly. With $a = {n \over 2}$ or $a = {n - 1 \over 2}$ according as $n$ is even or odd, we can determine
the $U$-ranks of $\pi$ and $\rho$ by examining the ranks of their characters of their restrictions to a 
standard subgroup $U_a$ of size $a \times a$. The assumption of part a) says that the maximum of these ranks 
are $k$ for $\pi$, and $\ell$ for $\rho$.

If we take the tensor product then 
$(\pi \otimes \rho)_{|U_a} \simeq (\pi_{|U_a}) \otimes (\rho_{|U_a})$. In
an abelian group such as $U_a$, the tensor product is formed by taking the sum in the Pontrjagin dual of 
each character appearing in $\pi$ with each character appearing in $\rho$. 
As noted, the rank of the sum of two characters is at most the sum of their ranks. Thus, if every character of $\pi$ 
has rank at most $k$, and every character of $\rho$ has rank at most $\ell$, the characters appearing in
$(\pi_{|U_a}) \otimes (\rho_{|U_a})$ can have rank at most $k + \ell$. 
Taking into account the given restrictions
on $k$ and $\ell$, this implies statement a) of the proposition. 
Taking into account the definition of asymptotic rank,
this argument also establishes statement b). \Qed

\vskip .1 in

{\bf Remark 6.3.2 (Filtrations of the representation ring).} As noted, both conditions of

\vskip .06 in

 {\bf i)} $U$-rank less than or equal to $k$; and

\vskip .06 in

{\bf ii)} asymptotic rank less than or equal to k,

\vskip .06 in

\noindent give rise to filtrations of the dual set $\widehat {G}_n$ of irreducible representations of $G_n$. 
Moreover, Proposition 6.3.1 shows that these filtrations are compatible with tensor product. This means that if we let 
$R(G_n) ={\bf Z}(\widehat {G}_n)$ be the representation ring of $G_n$ (aka the {\it Grothendieck ring}), namely the 
integer linear combinations of irreducible representations of $G_n$, equipped with tensor product for multiplication, 
then the submodules $R(G_n)_{A,\leq k} = {\bf Z}((\widehat {G}_n)_{A,\leq k} \subset R(G_n)$, of representation of asymptotic rank less or equal to $k$, give us an algebra filtration of $R(G_n)$, i.e., we have
$$
\cdots R(G_n)_{A,\leq k} \subset R(G_n)_{A,\leq k+1}\cdots,
$$
with $R(G_n)_{A,\leq i} \cdot R(G_n)_{A,\leq j} \subset R(G_n)_{A,\leq i+j}$ and $\bigcup_k R(G_n)_{A,\leq k}= R(G_n)$. We could get an analogous filtration using $U$-rank instead of asymptotic rank. 

\vskip .1 in

{\bf 6.4. Tensor rank.}

\vskip .1 in

Remark 6.3.2 suggests yet another way of defining a rank for representations. We start with representations of 
the smallest possible non-zero rank:
\vskip .06 in

$\bullet$ For $Sp_{2n}$, these are the four components of the oscillator representations, 
the representations that appear in the dual pairs $(Sp_{2n}, O_1)$, for the two possible 
1-dimensional orthogonal groups. 
\vskip .06 in

$\bullet$ For $GL_n$, these are representations that appear in the dual pair $(GL_n,GL_1)$, or these twisted by a character
(which are the rank 0 representations of $GL_n$). 
\vskip .06 in

$\bullet$ For the split orthogonal groups $O_{n,n}$, one starts with the representations of rank 2, that appear in the oscillator representation for the dual pair $(O_{n,n}, Sp_2)$. 
\vskip .05 in

And then, from the starting point of low rank representations, we take tensor products, and introduce the following notion:

\vskip .1 in 

{\bf Definition 6.4.1 (Tensor rank).} We say that 
\vskip .06 in

{\bf a)} A representation $\sigma$ in $\widehat{Sp}_{2n}$ 
has {\it tensor rank} less than or equal to $k$
 if it appears as a component of an $\ell$-fold tensor product of components
  of the oscillator representation, for $\ell \leq k$.
\vskip .06 in

{\bf b)} A representation $\nu$ in $\widehat{GL}_n$ has {\it tensor rank} less than or equal to $k$ if it appears
 as a component of an $\ell$-fold tensor 
product of components of the representations of $GL_n$ appearing
 in the oscillator representation for the dual pair 
$(GL_n, GL_1)$ for some $\ell \leq k$, or the tensor product of one of these times a character.
\vskip .06 in

{\bf c)} A representation $\tau$ in $\widehat{O}_{n,n}$ has {\it tensor rank} less than or equal to $2k$
if it appears as a component of an $\ell$-fold tensor product of components of the representations of $O_{n,n}$ 
appearing in the dual pair $(O_{n,n}, Sp_2)$, for some $\ell \leq k$.
\vskip .1 in

The tensor rank of a representation will be denoted sometime by $\otimes$-rank or $rank_\otimes$.

\vskip .1 in

{\bf Remark 6.4.2.} Note that our definition of tensor rank is given in a uniform manner in the cases of $GL_n$ and $Sp_{2n}$ (observe that $Sp_{2n}$ has no non-trivial characters). The case of $O_{n,n}$ is treated differently, as the relation between tensor rank and the operation of tensoring with the sign (= determinant) character, which is extremely interesting, requires further investigation.

\vskip .1 in

{\bf 6.5. Comparisons.}

\vskip .1 in

Since the representations used to define tensor rank have $U$-rank equal to their tensor rank, 
and since $U$-rank is compatible with tensor product by Proposition 6.3.1, 
it follows that the $U$-rank
of a representation will be less than or equal to the tensor rank. However, tensor rank, 
like asymptotic rank, potentially can, and in fact will, make sense for values greater than the maximum value of $U$-rank. 

Theorem 5.3.1 tells us exactly what are the representations of tensor rank $k$, for $k \leq n$ if $G_n = Sp_{2n}$
or $G_n = O_{n,n}$, and for $k \leq {n \over 2}$ for $G_n = GL_n$. 
These representations are parametrized by the unitary duals of the other member of the dual pair. It is clear that these representations will appear in the restrictions from larger groups of the same type. Therefore, their asymptotic rank is also at most $k$. We conclude that asymptotic rank is less than or equal to tensor rank.

Also, we can give explicit upper bounds on the tensor rank of any representation. For example, all representations
of $GL_n$ clearly appear in the dual pair $(GL_n, GL_n)$. This implies that the tensor rank of any representation
of $GL_n$ is at most $n$. We can get analogous bounds for $Sp_{2n}$ and $O_{n,n}$. Combining this with the previous 
paragraph, we can conclude that asymptotic rank is bounded on $\widehat {G}_n$ for any of the $G_n$ we have been considering.
 
 Summarizing, we have two proposals for extending the idea of rank to representations that are not of small $U$-rank.
 One is asymptotic rank, and the other is tensor rank. We have the inequalities
 $$
 U{\rm -rank} \ \ \leq \ \ {\rm asymptotic \  rank}  \ \ \leq \ \ {\rm tensor \  rank}. \eqno (6.5.1)
 $$
 Since tensor rank is easily seen to be a bounded function on $\widehat {G}_n$ for any of the $G_n$ under discussion,
 it follows also that asymptotic rank is bounded.
 
The $U$-rank is defined over a much more limited range than asymptotic rank or tensor rank, so there are clearly representations for which the inequalities in (6.5.1) are strict.  An important question is, whether this happens for reasons other than the obvious ones. In particular, are there any representations of low $U$-rank that have strictly larger tensor rank? We do not know any examples of this, and hope that it does not happen. In fact, the Exhaustion conjecture 5.3.3, if true, tells us that 

\vskip .1 in 
 
{\bf Conjecture 6.5.2 (Agreement of ranks).} For low $U$-rank representations the inequalities in (6.5.1) should be replaced by equalities.
  
 \vskip .1 in
  
   From the definitions, we can see that all three notions of rank are compatible with tensor product:
  if $\pi$ and $\rho$ are representations of $G_n$, of ranks less than or equal to
  $k$ and $\ell$ respectively, then $\pi \otimes \rho$ will have rank less than or equal to $k + \ell$, 
  for any of the three notions of rank. Thus, if we use any one of the notions of rank to define a filtration on 
  $\widehat {G}_n$, and extend this in the obvious way to the Grothendieck ring $R(G_n)$ of representations of $G_n$, 
  this gives us an algebra filtration on $R(G_n)$. One can then also pass to the associated graded ring.

 \vskip .2 in
  
  {\bf 7. The Oscillator semigroup}
  
  \vskip .15 in
  
  We have seen that the uniqueness of an irreducible representation of the Heisenberg group
  gives rise to a representation of the symplectic group. It turns out that this representation 
  can be extended further, from the symplectic group to a semigroup containing the symplectic
  group as its set of invertible elements. This semigroup has an elegant relationship with 
  symplectic geometry, and it sheds light on the oscillator representation itself, especially 
  on its restriction to dual pairs. This section is devoted to a discussion of the oscillator
  semigroup. Many of the results discussed are analogs for finite fields of the results in the 
  review paper [Howe87], and can be found in the unpublished manuscript [Howe73-2].
  
  \vskip .1 in 
  
  {\bf 7.1. Quantization of Lagrangians.}
  \vskip .1 in
  
 Let $W$ be a symplectic vector space. Recall that a subspace $L$ of $W$ is called
 {\it Lagrangian} if the symplectic form restricted to $L$ is identically zero, and $L$ is maximal
 with respect to this property. We also say that $L$ is a maximal isotropic subspace. 
 In the concrete realization of a symplectic vector space $W={\bf F}_q^{2n}$ in \S 3.2, the subspace $X$ consisting
 of vectors whose last $n$ coordinates vanish, and the subspace $Y$ consisting of vectors
 whose first $n$ coordinates vanish, are both Lagrangian subspaces.
 
 The symplectic group $Sp(W)$ will permute the Lagrangian subspaces. It is a corollary of 
 Witt's Theorem [Artin57] that $Sp(W)$ acts transitively on the set of all Lagrangian subspaces. 
 It is also easy to check that the subgroup of $Sp(W)$ that leaves the subspace $X$ invariant
 is the Siegel parabolic subgroup $P$ of formulas (3.2.3). Thus the set $\mathcal{L}(W)$ of Lagrangian
 subspaces is in natural bijection with the coset space $Sp(W)/P$.
 \vskip .08 in

{\bf The quantization vector attached to a Lagrangian.} 
If $L \subset W$ is Lagrangian, then $A_L =L \oplus {\bf F}_q$ is a maximal abelian subgroup of the Heisenberg group $H(W)$ (see (5.1.1a)). We can extend the character $\chi$ of ${\bf F}_q$ (see Theorem 5.1.2) to
a character $\tilde \chi_{_L}$ of  $A_L$ by letting $\tilde \chi_{_L}$ be trivial on $L$. As was stated in Theorem 5.1.2, the induced representation 
 $$
 Ind_{A_L}^{H(W)} \tilde \chi_{_L},
 $$
 will be a realization of the unique irreducible representation $\rho_{\chi}$ of $H(W)$
 with central character $\chi$.
 
 It follows by Frobenius reciprocity that there is a unique (up to multiple) non-zero vector in $\rho_{\chi}$
 that is invariant under $L$. We denote this vector (more precisely a line minus zero) by $q_{_L}$, and call it the
 {\it quantization} of $L$. This gives us the {\it quantization mapping}
 $$
 q: L \longmapsto q_{_L}, \eqno (7.1.1)
 $$
from $\mathcal{L}(W)$ to the projective space of $\rho_{\chi}$. 
From the relation (5.1.6) defining the
oscillator representation, we can conclude that the quantization map is equivariant with
 respect to the natural action of $Sp(W)$ on $\mathcal{L}(W)$ and the action of $Sp(W)$ on
 the projective space of $\rho_{\chi}$ induced by the oscillator representation.
\vskip .1 in

{\bf Remark 7.1.2.} Below, when we mention $q_L$, we are intending to specify it only up to a non-zero scalar multiple.
In particular, different formulas for $q_L$ may differ from each other by a scalar multiples,
which we may not specify explicitly.
\vskip .1 in

{\bf Formula for the quantization vector.} We can describe the vector $q_{_L}$ explicitly in 
the realization of $\rho_{\chi}$ given in formulas (5.1.5). This realization is constructed using a 
complementary pair $(X, Y)$ of Lagrangian subspaces of $W$. Using this pair of Lagrangians, we can 
identify the set $\mathcal{L}(W)$ with the collection of pairs $\{ (Y',B_{Y'})\}$, where $Y' \subset Y$ is any 
 subspace, and $B_{Y'}$ is any symmetric bilinear form on $Y'$. Indeed, let us denote by $p_Y$ the natural projection 
 $X \oplus Y\rightarrow Y$. Then we associate with a Lagrangian $L$ in $W$ the pair  
 $$
 Y_L = p_Y(L), \ \ B_{Y_L}(\vec y, \vec{y'}) = <\vec l, \vec{y'}>, \hskip 1 in
 $$
where $\vec y,\vec{y'} \in Y_L$, and $\vec l$ is any vector in $L$ such that $p_Y(\vec l)=\vec y$, i.e., 
such that $\vec y - \vec l$ belongs to $X$. Note that $Y_L=(X\cap L)^\perp \cap Y$, where $(X \cap L)^{\perp}$ 
indicates the usual notion of orthogonal subspace relative to a given bilinear form. That is, $(X \cap L)^{\perp}$ consists of all vectors in $W$ whose symplectic pairing with any vector in $X \cap L$ is zero. In particular, since $\vec l$ is well-defined 
up to a vector in $X \cap L$, its choice does not affect the pairing $<\vec l, \vec{y'}>$, since $\vec{y'}$ is in 
$(X \cap L)^{\perp}$. 

\vskip .05 in
 
Finally, in the realization of $\rho_{\chi}$ via formulas (5.1.5), we can describe the vectors $q_{_L}$ as follows:
$$
q_{_L} = \chi(- {1 \over 2} B_{Y_L}) \delta_{Y_L}, \eqno (7.1.3)
$$
where $\delta_{Y_L}$ is counting measure on the subspace $Y_L \subset Y$,
and $B_{Y_L}$ is the function on $Y_L$ taking the value $B_{Y_L}(\vec y, \vec y)$ at $\vec y$.
We recall that this definition is only up to a scalar multiple, and a different normalization 
of $q_{_L}$ might be convenient, according to context. 

\vskip .06 in

{\bf Example 7.1.4.} We note, in particular, the two examples:
\vskip .06 in
$\bullet$ For the Lagrangian $X$ we have
$$
q_{_X} = \textrm{the delta function at the origin of} \,\,Y.
$$
\vskip .06 in
$\bullet$ For the Lagrangian $Y$ we have
$$
q_{_Y} = \textrm{the constant function 1 (= the counting measure) on} \,\,Y.
$$
\vskip .1 in

{\bf 7.2. Conjugation action of $\omega(Sp_{2n})$.}
\vskip .1 in

The basic Stone-von Neumann-Mackey Theorem (Theorem 5.1.2) can be refined slightly,
by looking more closely at the mapping from the group algebra of $H(W)$ to the space
of operators on the space of the representation $\rho_{\chi}$. This space 
can be regarded as a Hilbert space, denoted $HS(\rho_{\chi})$, with respect to the Hilbert-Schmidt inner product
$$
B_{HS}(T', T) = {1 \over \dim \rho_{\chi}}trace (T'T^*), \eqno (7.2.1)
$$
where $T^*$ is the Hermitian transpose of the operator $T$ (with respect to the quadratic
extension ${\bf C}$ of ${\bf R}$). 

\vskip .1 in

{\bf Remark 7.2.2.} The normalization in (7.2.1) makes the identity, and all unitary operators,
into unit vectors. Note that his normalization does not make sense in the case of local fields, and a different normalization should be used there. 

\vskip .1 in

It is not difficult to verify [Gurevich-Hadani07],[Howe73-2],[Howe87] that
\vskip .1 in
{\bf Theorem 7.2.3 (Weyl transform).} The mapping 
$$
\rho_{\chi}: L^2(W) \longrightarrow HS(\rho_{\chi}), \eqno (7.2.4)
$$
given by $f \mapsto \sum_{\vec{w} \in W} f(w)\rho_\chi (w,0)$ is a canonical isometry of Hilbert spaces. In particular, the 
operators $\rho_{\chi}(\vec w, 0)$ for $\vec w$ in $W$, form an orthogonal basis for $HS(\rho_{\chi})$.
\vskip .1 in
The mapping (7.2.4) is called the {\it Weyl transform}.
\vskip .1 in 

Combining Theorem 7.2.3 with the defining relation (5.1.6) for the oscillator representation
$\omega$ yields the following result.
\vskip .1 in
{\bf Corollary 7.2.5.} The action on $HS(\rho_{\chi})$ by conjugation by $\omega (g)$, for $g$ in $Sp(W)$, is 
equivalent to the action of $g$ on $L^2(W)$ by means of its permutation action of $g$ on $W$.
\vskip .1 in
This in turn has the following appealing consequence. Recall that the {\it commutant}
${\mathcal A}_{\Gamma}$
of a given set $\Gamma$ of operators is the set of all operators that commute with each 
element of $\Gamma$. For any collection $\Gamma$ of operators, the commutant  
${\mathcal A}_{\Gamma}$ will be an algebra: closed under the operations of addition, 
multiplication by scalars, and operator multiplication. 
If  the set $\Gamma$ is closed under taking Hermitian conjugates, for example if it is a group
of unitary operators, then ${\mathcal A}_{\Gamma}$ will likewise be invariant under taking Hermitian
 conjugates. Such an algebra is called {\it self-adjoint}. We sometimes also call the commutant
of $\Gamma$ the {\it commuting algebra} of $\Gamma$.

\vskip .1 in 

{\bf Remark 7.2.6.} In fact, mapping (7.2.4) is an isomorphism of algebras, where $L^2(W)$ is equipped with the 
multiplication given by {\it twisted convolution}:
$$
(f_1 * f_2)(w) = \sum_{w_1 + w_2 = w} f_1(w_1)f_2(w_2)\chi(\frac{1}{2}<w_1,w_2>).
$$  

{\bf Corollary 7.2.7.} For any subgroup $G$ of $Sp(W)$, the commutant ${\mathcal A}_{\omega(G)}$ 
has an orthogonal basis corresponding via the map $\rho_{\chi}$ to the $G$-orbits in $W$.

\vskip .1 in

{\bf Example 7.2.8.} The full group $Sp(W)$ has two orbits when acting on $W$: the zero vector,
and the collection of all non-zero vectors. It follows that the commuting algebra 
${\mathcal A}_{\omega (Sp(W))}$ is two-dimensional. This in turn implies that $\omega$ is composed of two 
irreducible summands (These summands are the eigenspaces for
the center of $Sp(W)$, which consists of $\pm I_W$, where $I_W$ denotes the
identity operator on $W$). The algebra ${\mathcal A}_{\omega (Sp(W))}$ is the direct sum of the scalar 
operators on each summand.

\vskip .1 in 

{\bf 7.3. The signed double and the oscillator semigroup.}
\vskip .1 in 

The space $HS(\rho_{\chi})$ can also be thought of as a module for $H(W)$. Indeed,
it can be thought of an $H(W)$-module in two ways: multiplication on the left, and multiplication
on the right. This double action is described by the formula: 
$$
\tilde \rho_{\chi} (h, h') (A) = \rho_{\chi}(h) A \rho_{\chi}((h')^{-1}),
$$
for elements $h$ and $h'$ in $H(W)$. This in fact defines a representation of $H(W) \times H(W)$, and this representation
is irreducible.

It is isomorphic to the {\it outer tensor product} of $\rho_{\chi}$ with 
the contragredient representation $\rho^*_{\chi}$:
$$
\tilde \rho_{\chi} \simeq \rho_{\chi} \boxtimes \rho^*_{\chi}. 
$$
Indeed, a similar construction could be made for any irreducible representation of any finite
group. What is special about this situation is that we can again regard $\tilde \rho_{\chi}$ 
as a representation of a single Heisenberg group, the Heisenberg group of the 
{\it signed double} $2W$ of $W$. We define 
$$
2W = W\oplus W,
$$
as vector spaces, but instead of taking the direct sum of the symplectic form $< \ , \ >$ 
with itself, we put the \underbar {negative} of $< \ , \ >$ on the second factor. Thus
$2W$ is understood to have the symplectic form
$$
<(\vec w_1, \vec w_2), (\vec w'_1, \vec w'_2)>_{2W} \ 
 =  \ <\vec w_1, \vec w_1'> - <\vec w_2, \vec w'_2>. 
 $$
 for vectors $\vec w_j$, $\vec w'_j$ in $W$. 
 It is not difficult to check that the $((\vec w, t), (\vec w', t')) \mapsto (\vec w \oplus \vec w', t - t')$ 
 defines a surjective group homomorphism from $\nu: H(W) \times H(W) \rightarrow H(2W)$.
 Moreover, we have the following factorization 
\[
\begin{tikzcd}
H(W) \times H(W) \arrow[rr,"\nu"] \arrow[dr,swap,"\tilde \rho_{\chi}"] && H(2W) \arrow[dl,"\tilde \rho_{\chi}"] \\
& GL(HS(\rho_{\chi}))
\end{tikzcd}
\]
of the action $\tilde \rho_{\chi}$ to the irreducible representation (that by abuse of notation we denote also by) $\tilde \rho_{\chi}$ of $H(2W)$, with central 
character $\chi$, given by
$$
{\tilde \rho_{\chi}}(\vec w \oplus \vec w',t)(A)=\chi(t)\rho_\chi(\vec w')A\rho_\chi(-\vec w). \eqno (7.3.1)
$$
\vskip .1 in 
 When we need to distinguish between the two summands of $2W$, we will denote
 the first one by $W^+$, and the second one by $W^-$.
 \vskip .1 in 
 We clearly have an embedding 
 $$
 Sp(W) \times Sp(W) \longrightarrow Sp(2W),
 $$
 by letting each copy of $Sp(W)$ act on the respective summand of $2W$.
 
 \vskip .1 in

{\bf The operator associated with a Lagrangian in $2W$.} 
Since $HS(\rho_{\chi}) \simeq L^2(W)$ (canonically, see equation (7.2.4)) is the space for the irreducible representation 
$\tilde \rho_{\chi}$ of $H(2W)$, we have the quantization map (see equation (7.1.1))
$$
q : {\tilde L} \longmapsto q_{\tilde L}, \eqno (7.3.2)
$$  
from the variety $\mathcal{L}(2W)$ of Lagrangian subspaces of $2W$ to (the projective space of) 
 $HS(\rho_{\chi}) \simeq L^2(W)$. Thus, if $\tilde L$ is a Lagrangian subspace of $2W$,
 its quantization  $q_{\tilde L}$ will have an interpretation as an operator on $\rho_{\chi}$,
 and as a function on $W$. Sometime, we might denote this operator by $q({\tilde L})$. Note that formula (7.3.1) implies that $q({\tilde L})$ is the unique operator (up to a non-zero scalar multiple) satisfying
$$
\rho_\chi(\vec w')q({\tilde L})=q({\tilde L})\rho_\chi(\vec w), \eqno (7.3.3)
$$
for every $(\vec w,\vec w') \in  {\tilde L}$.
\vskip .1 in

{\bf A decomposition of the Lagrangian Grassmannian $\mathcal{L}(2W)$.}
 The decomposition $2W = W^+ \oplus W^-$ leads to a decomposition of the set $\mathcal{L}(2W)$, 
 according as to how a Lagrangian subspace of $2W$ sits with respect to the two summands. 
 Let $p^+$ and $p^-$ respectively be the projections of $2W$ on $W^+$ and on $W^-$. Let $L$ be a Lagrangian 
 subspaces of $2W$. Set 
 $$
 L^+ = L \cap W^+; \hskip .5 in L^- = L \cap W^-.
 $$
 It is evident that $L^{\pm}$ must be isotropic subspaces of $W^{\pm}$.
 Then, elementary symplectic linear algebra considerations [Howe87] shows that
 $$
 \dim L^+ = \dim L^-, \hskip .2 in p^+(L) = (L^+)^{\perp}, \hskip .2 in p^-(L) = (L^-)^{\perp}.
 \eqno (7.3.4)
 $$
 Here $(L^+)^{\perp}$ denotes the orthogonal subspace in $W^+$ with respect to the 
 original symplectic form, and $(L^-)^{\perp}$ is the analogous subspace of $W^-$.
 The relations (7.3.4) also imply that $L/(L^+ \oplus L^-)$ is the graph of an isomorphism $\alpha_L$ 
 between $p^+(L)/L^+$ and $p^-(L)/L^-$. These spaces both inherit symplectic forms from $W^{\pm}$,
 by restriction and reduction modulo $L^{\pm}$. The condition that $L$ be isotropic translates to
 the condition that $\alpha$ should be an isometry between $p^+(L)/L^+$ and $p^-(L)/L^-$ (but
 with  $p^-(L)/L^-$ endowed with the original form, not its negative).
 
 Thus, we get a parametrization of $\mathcal{L}(2W)$ by triples 
 $$
 \mathcal{L}(2W) \simeq \{ (L^+, L^-, \alpha) \},
 $$
  where $L^+$ and $L^-$ are isotropic subspaces of $W$, of the same dimension, 
  and $\alpha$ is an isometry between $(L^+)^{\perp}/L^+$ and $(L^-)^{\perp}/L^-$. 
  We can describe the action of  $Sp(W) \times Sp(W)$ on $\mathcal{L}(2W)$ using these parameters. 
  We have
 $$
 (g, g')(L^+, L^-, \alpha) = (g(L^+), g'(L^-), g^{-1} \alpha g'),
 $$
 for an element $(g, g')$ in $Sp(W) \times Sp(W)$. We can conclude from this that the number
 of $Sp(W)\times Sp(W)$ orbits in $\mathcal{L}(2W)$ is finite, and equal to $n+1$, if $\dim W = 2n$.
 Also, the set of elements for which $L^+ = \{0\} = L^-$ is naturally isomorphic to
 $Sp(W)$. Indeed, the graph of any element of $Sp(W)$ defines Lagrangian subspace of
 $2W$ with this property. For any infinite field, the subset of $\mathcal{L}(2W)$ filled out by $Sp(W)$ is Zariski open
  and dense. In this sense, $\mathcal{L}(2W)$ may be regarded as a ``compactification" or ``completion"
  of  $Sp(W)$.
 \vskip .1 in

{\bf The oscillator semigroup.} 
Two notable Lagrangian subspaces of $2W$ are:
 $$
 \Delta(W)={(\vec w, \vec w): \vec w \in W\},\,\, {\rm and}\,\, 
 \Delta^{-}(W)=\{(\vec w, - \vec w): \vec w \in W\}},
 $$
 called, respectively, the {\it diagonal} and the {\it anti-diagonal}. These may be thought of as the 
 graphs of the identity map of $W$ to itself, and of the negative of the identity map of $W$. 
 The vectors in $L^2(W)$ corresponding to these Lagrangians (see (7.3.2)) are: 
\vskip .04 in
$\bullet \,\,\, q_{_{\Delta (W)}} = \delta_0, \,\,\, \textrm{the delta function at the origin of} \,\,W;$ and
\vskip .04 in
$\bullet \,\,\,q_{_{\Delta^- (W)}} = \, \textrm{the constant function 1 (= the counting measure) on} \,\,W.$
\vskip .06 in
On the other hand, the operators in $HS(\rho_{\chi})$ corresponding to these two Lagrangians in $2W$ are:
$$
q({\Delta (W)}) = \textrm{the identity operator in} \,\, HS(\rho_{\chi})  = \omega (e),
$$
where $e$ here denotes the identity element of $Sp(W)$, and
$$
q({\Delta^- (W)}) = \omega (-e),
$$
where $-e$ here denotes the negative identity element in $Sp(W)$.

The Lagrangians $\Delta(W)$ and $\Delta^-(W)$ are, respectively, the parallels for $2W$ to the two Lagrangian subspaces $X$ and $Y$ of $W$, 
that were used to construct the realization (5.1.5) of  $\rho_{\chi}$. 

In view of the action of $Sp(W)$ on $\mathcal{L}(2W)$ through the embedding of $Sp(W)$ in
$Sp(2W)$ (via its action on $W^-$) and the equivariance of the quantization mapping  (7.3.2), we conclude that
$$
q(g(\Delta(W)) = \omega (g),
$$
for any $g$ in $Sp(W)$. Thus, the group $\omega (Sp(W))=\{\omega(g);\,\, g \in Sp(W)\}$ of operators is contained in
$$
{\mathfrak S}(2W) = q(\mathcal{L}(2W))=\{q({\tilde L});\,\,\, {\tilde L}\in {\mathcal L}(2W)\}, \eqno (7.3.5)
$$ 
and left and right multiplication by $\omega(Sp(W))$ preserves $q(\mathcal{L}(2W))$.
In fact, the whole set ${\mathfrak S}(2W)$ is closed under multiplication.
\vskip .1 in

{\bf Theorem 7.3.6 (The oscillator semigroup).} The set ${\mathfrak S}(2W)$ of operators is a semigroup.
\vskip .1 in

We will call ${\mathfrak S}(2W)$ the {\it oscillator semigroup} for $\rho_{\chi}$.
\vskip .1 in
A proof of Theorem 7.3.6 first appeared in [Howe73-2]. We will give here an argument which is closely related to the works [Guillemin-Sternberg79],[Weinstein80] on the symplectic ``category", and especially the work [Nazarov95] on the oscillator semigroup over a $p$-adic field. We will explicitly describe the semigroup law. Define the {\it composition} $\tilde M \circ \tilde L   \subset 2W$ of two Lagrangians $\tilde M,\tilde L \in \mathcal{L}(2W)$ by
$$
\tilde M \circ \tilde L=\{(\vec w,\vec w')\, | \,\, \exists \, (\vec w,\vec w'') \in {\tilde L}, \,\, (\vec w'',\vec w') \in {\tilde M}\}.
$$
\vskip .05 in
{\bf Proposition 7.3.7.} For every $\tilde M,\tilde L \in \mathcal{L}(2W)$ we have
\vskip .05 in
{\bf i)} $\tilde M \circ \tilde L \in \mathcal{L}(2W)$. 
\vskip .05 in
{\bf ii)} $q({\tilde M})q({\tilde L}) = \alpha({\tilde M},{\tilde L})\cdot q({\tilde M} \circ {\tilde L}),\,\, {\rm for \,\, \alpha({\tilde M},{\tilde L}) \in {\bf C}^*}$.    
\vskip .1 in
The meaning of 7.3.7 i) is that $(\mathcal{L}(2W),\circ,\Delta(W))$ is a semigroup (with $\Delta(W)$ being the identity element), while 7.3.7 ii) is telling us that indeed ${\mathfrak S}(2W)$ is a semigroup and the quantization mapping 
$$
q:\mathcal{L}(2W) \longrightarrow {\mathbb P}(HS(\rho_\chi)),
$$ 
is a (projective) semigroup representation. It follows that ${\mathfrak S}(2W)$ is in a natural way a central extension
$$
1\longrightarrow {\bf C}^* \longrightarrow {\mathfrak S}(2W) \longrightarrow \mathcal{L}(2W)\longrightarrow 1.
$$
\vskip .1 in
{\it Proof of Proposition 7.3.7.} The proof of i) is straightforward [Guillemin-Sternberg79]. For ii) one can use the characterization (7.3.3) as was done in [Nazarov95]. Take $(\vec w,\vec w') \in {\tilde M} \circ {\tilde L}$ and an element
$\vec w'' \in W$ such that $(\vec w,\vec w'') \in {\tilde L}$ and $(\vec w'', \vec w') \in {\tilde M}$. Now compute,  
$$
q({\tilde M})q({\tilde L})\rho_\chi(\vec w)=q({\tilde M})\rho_\chi(\vec w'')q({\tilde L})=\rho_\chi(\vec w')q({\tilde M})q({\tilde L}),
$$
where the first and second equalities are by (7.3.3). We conclude that ii) holds true. 
\Qed       
\vskip .12 in

{\bf Remark 7.3.8.} Every element of $\mathcal{L}(2W)$ that is not in $Sp(W)$ is in the $Sp(W)$-orbit of some triple of the form $(Z, Z, I_{Z^{\perp}/Z})$, where $Z \subset W$ is non-trivial isotropic subspace (see discussion on decomposition of Lagrangian Grassmanian above). Up to multiples, $q((Z, Z, I_{Z^{\perp}/Z}))$ is counting measure on $Z$. As an operator, it is (a multiple of) projection to the space of $\rho_{\chi}^Z$ of $Z$-invariant vectors in $\rho_{\chi}$. Note that $Z^{\perp}/Z$ inherits a symplectic structure from $W$, and the space $\rho_{\chi}^Z$ forms the standard irreducible representation for $H(Z^{\perp}/Z)$. The $Sp(W)$ equivariance of the quantization map then tells us that $q(Z, Z', \alpha)$ is a partial isometry from $\rho_{\chi}^Z$ to $\rho_{\chi}^{Z'}$, implemented by means of some element $\beta$ in $Sp(W)$ such that $\beta(Z) = Z'$, and such that the map from $Z^{\perp}/Z$ to $Z'^{\perp}/Z'$ induced by $\beta$ is $\alpha$.

\vskip .2 in

{\bf 8. The oscillator semigroup, dual pairs, and old vs. new spectrum}

\vskip .15 in 

One of the most interesting phenomena associated to the oscillator representation $\omega$
of a symplectic group $Sp(W)$ is the insight it gives into representations
 of other classical groups via the study of dual pairs of subgroups. The main observation
 is that, if $(G, G')$ is a dual pair in $Sp(W)$, then $\omega (G)$ and $\omega (G')$ 
 approximately generate mutual commutants inside the algebra $HS(\rho_{\chi})$
 of operators on the space of the underlying representation $\rho_{\chi}$ of the
 Heisenberg group $H(W)$. Unfortunately, except in a few cases, it is not strictly true
 that $\omega(G)$ and $\omega (G')$ generate the full commutant of the other. The 
 challenge then is, to make precise sense in which they ``approximately" do so. 
 The notion of $U$-rank provides an approach to this issue. As indicated in Theorem 5.3.1,
 in a large number of cases, if one restricts attention to the representations of $G$ of maximal
 possible $U$-rank, then $\omega (G')$ generates the commutant of $\omega (G)$ on this 
 subspace.

\vskip .1 in

{\bf 8.1. The oscillator semigroup and dual pairs.}

\vskip .1 in

Another approach to the issue is to enlarge the collection of operators under consideration, by looking at the oscillator semigroup. It turns out [Howe73-2] that this always gives us a description of the commutant of the image under the oscillator representation of either member of a dual pair. 
 
 Let ${\mathfrak S}(2W)=q(\mathcal{L}(2W))$ be the oscillator semigroup in $HS(\rho_{\chi})$, for the standard
  representation $\rho_{\chi}$ of the Heisenberg group $H(W)$, for the symplectic
  vector space $W$ (see equation (7.3.5)). The symplectic group $Sp(W)$, or more correctly, its image
  $\omega (Sp(W))$ in the operators on $\rho_{\chi}$, is a subgroup of 
  ${\mathfrak S}(2W)$, and so can act on it by multiplication on the left, by multiplication on the right, and by conjugation. 
  
Let $(G, G')$ be a dual pair of subgroups of $Sp(W)$. Consider
the set ${\mathfrak S}(2W)^G$ of elements of ${\mathfrak S}(2W)$ that are invariant under conjugation by
$G$, or equivalently, that commute with $G$, as operators on the space $\rho_{\chi},$
$$
{\mathfrak S}(2W)^G = \{q({\tilde L}) \in {\mathfrak S}(2W);\,\,\,\omega(g) q({\tilde L})=q({\tilde L})\omega(g),\,\, 
{\rm for \, every}\,\, g \in G\}.
$$
Evidently, since $G'$ is by definition the set of elements in $Sp(W)$ that commute with $G$,
we have $\omega(G') \subset {\mathfrak S}(2W)^G$. Also, it is easy to see that ${\mathfrak S}(2W)^G$ is a subsemigroup
of ${\mathfrak S}(2W)$. In particular, it is a semigroup. It bears a similar relationship to $G'$,
as a ``compactification" or ``completion", that $\mathcal{L}(2W)$ bears to $Sp(W)$.
For our current purposes, our main interest in ${\mathfrak S}(2W)^G$ is the following result.

\vskip .1 in

{\bf Theorem 8.1.1.}  Given a dual pair $(G, G')$ in $Sp(W)$, the semigroup ${\mathfrak S}(2W)^G$ spans the commuting algebra of $\omega (G)$ in $HS(\rho_{\chi})$.

\vskip .1 in

{\it Proof}. The proof of this result, like that of many others in the theory of the classical groups,
devolves into considering a collection of cases. We will discuss here the most interesting case,
that of the pair $(G, G') = (Sp(V), O(V'))$, corresponding to a symplectic vector space $V$,
 a space $V'$ endowed with an inner product, both contained in $Sp(W)$, 
 with $W = V \otimes V'$.
 
By duality, we can also describe $V \otimes V'$ as $Hom (V', V)$, or as $Hom (V, V')$
(with the ``Hom"s here signifying ${\bf F}_q$ linear maps). We let $< \ , \ >$ denote the symplectic
 form on $V$, and let $(\ , \ )$ denote the inner product on $V'$. We can realize the isomorphism
$Hom (V', V) \simeq Hom (V, V')$ by the mapping 

$$
T' \longleftrightarrow T,
$$
where $T$ and $T'$ are related by the condition
$$
<T'(\vec v'), \vec v> = (\vec v', T(\vec v)), \eqno (8.1.2)
$$
for $T'$ in $Hom(V', V)$, $\vec v$ in $V$ and $\vec v'$ in $V'$.

The relation (8.1.2) implies that $T'T$ is an operator on $V$ such that the bilinear form
$$
B_{T'T}(\vec u, \vec v) = <T'T(\vec u),\vec v>,
$$
is symmetric (equal to $(T(\vec u),T(\vec v)$). Similarly, the operator $TT'$ induces on $V'$ the bilinear form
$$
\Lambda_{TT'}(\vec u', \vec v') = (TT'(\vec u'), \vec v'),
$$
which is skew-symmetric (equal to $<T'(\vec u'), T'(\vec v')>$). 

\vskip .06 in 

{\bf Remark 8.1.3.} The collection of all operators $\beta$ on $V$ such that $<\beta(\vec u),\vec v>$
defines a symmetric bilinear form is naturally identified to $\mathfrak{sp}(V)$,  the symplectic Lie algebra
associated to $< \ , \ >$ [Jacobson62]. Likewise, the collection of all operators $\beta'$ on $V'$ such that the bilinear 
form $(\beta (\vec u'), \vec v')$ is skew symmetric may be identified to the Lie algebra $\mathfrak{so}(V)$, of 
the orthogonal group of the inner product $( \ , \ )$ [Jacobson62]. Thus, we have mappings
\[
\begin{tikzcd}
& T' \in Hom(V',V) \arrow[dl] \arrow[dr] & \\
T'T \in \mathfrak{sp}(V) & & TT' \in \mathfrak{so}(V')
\end{tikzcd}
\]
These maps are known in the symplectic geometry literature as {\it moment mappings}.

\vskip .06 in

We need to show that the commutant algebra ${\mathcal A}_{\omega(G)}$ of $G=Sp(V)$ in $HS(\rho_\chi)$ is equal to the subalgebra spanned by ${\mathfrak S}(2W)^G$. 

Using the Weyl transform (7.2.4) we can view those algebras inside $L^2(W)$. We know from Corollary 7.2.7 that the commutant ${\mathcal A}_{\omega(Sp(V))}$ of $Sp(V)$ is spanned by the characteristic functions of the $Sp(V)$ orbits in $W \simeq Hom (V',V)$. We will show that these are in fact in the span of ${\mathfrak S}(2W)^G$. 

Indeed, Witt's Theorem implies that the $Sp(V)$ orbit of $T'$ in $Hom (V', V)$ is determined by the following data:
\vskip .07 in

i) The skew-symmetric bilinear form $\Lambda_{TT'}$ on $V'$;
\vskip .05 in
ii) The kernel of $T'$.
\vskip .07 in 

On the other hand, equation (7.1.3) tells us that the elements of ${\mathfrak S}(2W)$ are given by ``exponential quadratic 
characters" on subspaces $Z$ of $W$, that is, by functions of the form 
$\chi(B(\vec z, \vec z))$,
where $B$ is a symmetric bilinear form on $Z$. The elements
of ${\mathfrak S}(2W)^G$ arise when the subspace $Z$ is invariant under $Sp(V)$, and the
bilinear form $B$ is invariant under $Sp(V)$ acting on $Z$.
A subspace $Z$ of $W = V \otimes V'$ that is invariant under $Sp(V)$ has the form
$V \otimes Z'$, where $Z'$ is a subspace of $V'$. A symmetric bilinear form on $Z$
that is invariant under $Sp(V)$ has the form $< \ , \ > \otimes \Lambda'$, where
$\Lambda'$ is a skew-symmetric bilinear form on $Z'$. 
The values of an element ${\mathfrak S}(2W)^G$ on an $Sp(V)$ orbit $\mathcal{O}$
in  $W$ is of course equal to zero unless $\mathcal{O}$ is contained in $Z$. If $\mathcal{O}$ is contained in $Z$, then the value
of the element of ${\mathfrak S}(2W)^G$ corresponding to the skew form $\Lambda'$ on $Z'$ 
is the value of the pairing of $\Lambda'$ with the restriction to $Z'$ of the form 
$\Lambda_{TT'}$ defining the orbit. Thus the pairing of elements of ${\mathfrak S}(2W)^G$
with orbits essentially amounts to Fourier analysis on the skew-symmetric forms on
$Z'$. Since $\Lambda'$ is allowed to vary over all possible skew-symmetric forms,
it follows that sums of these elements can produce the point mass at the restriction to
$Z' $ of $\Lambda_{TT'}$. Since we can repeat this analysis for each possible invariant
subspace that could possibly contain a given orbit, we can also separate, i.e., span the characteristic function of, orbits that
correspond to the same form $\Lambda_{TT'}$, but have different kernels. It follows
that ${\mathfrak S}(2W)^G$ will separate $Sp(V)$ orbits in $W$, as claimed.

\Qed
\vskip .1 in 

{\bf Remark 8.1.4.} Notice that, although ${\mathfrak S}(2W)^G$ spans the commutant of $\omega (G)$,
and ${\mathfrak S}(2W)^{G'}$ spans the commutant of $\omega(G')$, it will typically be the case that
${\mathfrak S}(2W)^{G}$ and ${\mathfrak S}(2W)^{G'}$ do not commute with each other.
The linear spans of  $\omega (G)$ and $\omega (G')$ are (usually) slightly too small
to be mutual commutants, and the linear spans of ${\mathfrak S}(2W)^{G}$ and ${\mathfrak S}(2W)^{G'}$
are (usually) slightly too large to be mutual commutants; that is, they usually will not commute 
with each other.
\vskip .1 in

{\bf 8.2. Old and new spectrum.}
\vskip .1 in

Theorem 8.1.1 reveals that $\omega (G')$ usually will not generate the full
commutant of $\omega (G)$, and thus the restriction of $\omega$ to $G \times G'$
will not produce a well-defined matching between representations of $G$
and representations of $G'$. Nevertheless, it does allow one to establish a partial
result of that nature. Again, we take a pair $(Sp_{2n}, O_k)$ as an illustration
of what happens. As above, we let $V$ and $V'$ be the spaces on which the two
groups operate, so that $W = Hom (V', V)$ is the symplectic space such that
$Sp_{2n} \times O_k \subset Sp(W)$.

As discussed in \S 1, the inner product $( \ , \ )$ on $V'$ of which $O_k$ is the isometry
group is the orthogonal sum of an anisotropic inner product and a certain number $\ell$ of
hyperbolic planes. This in turn implies that the restriction of $\omega$ to $Sp_{2n}$ is
the tensor product of the representation $\omega_a$ coming from the anisotropic form,
 together with $\ell$ copies of the representation $\omega_{1,1}$ of $Sp_{2n}$ attached
to the dual pair $(Sp_{2n}, O_{1,1})$.
One way to interpret the discussion of \S 7.2 is, that $\omega_{1,1}$ is exactly
the permutation representation of $Sp_{2n}$ acting on $V$, 
 its defining symplectic vector space.
Since this representation contains the trivial representation (in fact, with multiplicity 2),
we can conclude that, if we tensor any representation $\tau$ of $Sp_{2n}$ with
 $\omega_{1,1}$,
the resulting representation $\tau \otimes \omega_{1,1}$ will contain a copy of $\tau$.
Thus, the subset of $\widehat {Sp}_{2n}$ that appears with non-zero multiplicity
in $\omega_a \otimes \omega_{1,1}^{\otimes^{\ell}}$ is strictly increasing with $\ell$.

Denote by $\widehat{Sp}_{2n}(a, \ell)$ the subset of $\widehat{Sp}_{2n}$ consisting of representations that appear with non-zero multiplicity in $\omega_a \otimes \omega_{1,1}^{\otimes^{\ell}}$. Then according to the previous paragraph, we have
$$
\widehat {Sp}_{2n}{(a, \ell-1)} \subset \widehat {Sp}_{2n}(a, \ell).
$$
We can call the set $\widehat {Sp}_{2n}(a, \ell-1)$ of previously appearing
representations the ``old spectrum", and the difference set
$$
\widehat {Sp}_{2n}{(a, \ell)} \smallsetminus \widehat {Sp}_{2n}{(a, \ell-1)},
$$
the ``new spectrum". In addition, let us denote by 
$$
\left (\omega_a \otimes \omega_{1,1}^{\otimes^{\ell}}\right )_G^{new}
\subset \omega_a \otimes \omega_{1,1}^{\otimes^{\ell}},
$$ 
the subspace that is the sum of the isotypic components for representations of 
$G$ in the ``new spectrum". With this terminology, we can state the following result.
\vskip .1 in

{\bf Proposition 8.2.1.} The operators 
$\left (\omega_a \otimes \omega_{1,1}^{\otimes^{\ell}} \right )(G')$
generate the commutant of $G$ acting on 
$\left (\omega_a \otimes \omega_{1,1}^{\otimes^{\ell}}\right )_G^{new}$. 
In particular, the action of $G \times G'$ on 
$\left (\omega_a \otimes \omega_{1,1}^{\otimes^{\ell}}\right )_G^{new}$
defines a bijection of some subset of $\widehat{G'}$ to 
$\widehat {Sp}_{2n}{(a, \ell)} \smallsetminus \widehat {Sp}_{2n}{(a, \ell-1)}$.
\vskip .1 in

{\it Proof.} For purposes of this argument, we will abbreviate 
$\omega_a \otimes \omega_{1,1}^{\otimes^{\ell}} = \omega$. We know from Theorem 8.1.1
that the commutant of $\omega (G)$ in the operators on the space of $\omega$ is
spanned by ${\mathfrak S}(2W)^G$. Besides $\omega (G')$, the other elements of 
${\mathfrak S}(2W)^G$ are gotten (see Remark 7.3.18) by multiplying with projections to the fixed vectors for
isotropic subspaces of $W$. If such a subspace is invariant for $Sp(V)$, it must have the form
$Z = V \otimes Z'$, where $Z'$ is an isotropic subspace of $V'$. The space
$\rho_{\chi}^{Z}$ of fixed vectors for $Z$ will be the standard irreducible representation
for the Heisenberg group $H(Z^{\perp}/Z)$, and the pair of groups
$(Sp_{2n}, O(Z'^{\perp}/Z')$ will be a dual pair inside $Sp(Z^{\perp}/Z)$.  The group
$O(Z'^{\perp}/Z')$ is a smaller group in the same Witt tower as $G'=O(V')$. The representation
of $Sp(V)$ on the space $\rho_{\chi}^{Z}$ will be the representation attached to the
dual pair $(Sp_{2n}, O(Z'^{\perp}/Z'))$, and therefore, the set of irreducible representations
of $Sp(V)$ appearing will belong to $\widehat {Sp}_{2n}{(a, \ell-1)}$, that is, they will
be old spectrum. It follows that the operators of  ${\mathfrak S}(2W)^G \smallsetminus \omega (G')$ will act 
by zero on the new spectrum. Therefore, on the new spectrum, the commutant must come
from $\omega (G')$, as desired. \Qed
\vskip .1 in

{\bf Remarks 8.2.2.} We note that
\vskip .05 in

{\bf a)} A result similar to Proposition 8.2.1 holds, with the appropriate modifications, for any of the previously described classical dual pairs $(G,G')$ , when $G'$ varies in a Witt tower. For example in the $(GL_n,GL_k)$ case.
\vskip .07 in

{\bf b)} Proposition 8.2.1 gives another sense in which $\omega (G')$ is ``approximately" the commutant of $\omega (G)$, and it does establish a
matching between certain representations of the two groups. However, it
does not provide any description of what representations of $G$ or of $G'$ are involved, except the very indirect condition of ``new spectrum". However, it is a natural extension of the $\eta$-correspondence (see Theorem 5.3.1). Indeed, when it applies, the notion of $U$-rank provides a clear criterion for a representation of the group $G$ to belong to the new spectrum. More precisely, it gives a sufficient criterion for a representation to belong to the new spectrum, and then Theorem 5.3.1 shows that this must account for the entire new spectrum, since the full group algebra of $G'$ is required to generate the commutant for these representations. It is also worth remarking that the version of Theorem 5.3.1 given in [Gurevich-Howe15] establishes a close relationship between corresponding representations, including a relationship between their dimensions.
\vskip .07 in

{\bf c)} On the other hand, Theorem 5.3.1 gives no insight into the nature of what
we have labeled here the ``old spectrum". In particular, it does not establish that it consists of representations that have appeared in correspondences for dual pairs $(G, G'')$, where $G''$ is a smaller member of the same Witt tower as $G'$.
\vskip .07 in

{\bf d)} The content of Sections  \S 7 and \S 8 is based on the paper [Howe73-2], which was circulated as a preprint in the 1970s, but never published.

\vskip .2 in 
 
{\bf 9. A description of the new spectrum for the dual pair $(GL_n,GL_k)$}

\vskip .15 in 

In this section we describe the new spectrum (that was defined in \S 8.2) for the dual pairs $(GL_n, GL_k)$. This description holds for all $k$, in particular, for $k > {n \over 2}$, so it provides an extension of Theorem 5.3.1 c), and yields a complete description of the irreps of arbitrary tensor rank (which can go up to $n$ for $GL_n$) using the eta correspondence. We will sketch the main ideas. Complete details will appear in [Gurevich-Howe17].

\vskip .1 in

{\bf 9.1. Intrinsic characterization for the new spectrum and tensor rank for $GL_n$.}

\vskip .1 in

We begin with an alternative characterization of the new spectrum and tensor rank that is available for $GL_n$. 

First we note that,

\vskip .07 in

{\bf Remark 9.1.1.} a representation $\rho$ of $GL_n$ has tensor rank $k$ if and only if, up to tensoring by a character, it appears in the $k$-fold tensor product of the natural permutation representation of $GL_n$ on ${\bf F}_q^n$, and does not appear in the $l$-fold tensor product of that representation for $l < k$. 
This $k$-fold tensor product is just the natural permutation representation on $k$ copies of ${\bf F}_q^n$. This action can be identified with the permutation representation of $GL_n$ on space of $n \times k$ matrices $M_{nk}$, by multiplication on the left.

\vskip .1 in

Next, let $\vec e_j$ be the $j$-th standard basis vector of ${\bf F}_q^n$, whose $j$-th coordinate is 1, and all other coordinates are equal to 0. 
Let $H_k$ be the subgroup of $GL_n$ defined by
$$
H_k = \{ g \in GL_n; \,\, ge_j=e_j \,\, {\rm for \, every }\,\, 1\leq j \leq k \}. \eqno (9.1.2)
$$
Then $H_k$ is conjugate to the subgroup that leaves any $k$-dimensional subspace of ${\bf F}_q^n$ point-wise fixed. Evidently $H_k \subset H_{k-1}$.
Now, we can write down the intrinsic criterion.
 
 \vskip .1 in
 
{\bf Proposition 9.1.3 (Intrinsic characterization for new spectrum tensor rank for $GL_n$).} Assume $0 \leq k \leq n$. A representation $\rho$ of $GL_n$ is in the new $k$-spectrum (respectively, has tensor rank $k$) if and only if it admits an invariant vector (respectively, eigenvector) for $H_k$ but not for $H_l$ for $l < k$.

\vskip .1 in
 
{\it Proof}. It is enough to prove the assertion on new $k$-spectrum (see Remark 9.1.1). 

By definition, a representation $\rho$ of $GL_n$ is in the new $k$-spectrum if and only if it appears in the $k$-fold tensor product (and not before) of the natural permutation representation of $GL_n$ on its standard vector space ${\bf F}_q^n$. As remarked earlier, this action can be identified with the permutation representation of $GL_n$ on space of $n \times k$ matrices $M_{nk}$, by multiplication on the left. Of course, $GL_k$ will act also on $M_{nk}$, by multiplication on the right. The two actions commute with each other. A $GL_n \times GL_k$ orbit is defined by the rank of its elements. The rank can vary from 0 to $\min (k, n)$, which we assume is $k$.
 
For $1 \leq \ell \leq k$, the group $H_{\ell}$ is the stabilizer in $GL_n$ of the matrix $J_{\ell}$, whose first $\ell$ diagonal entries are 1, and all other entries are 0. The matrix $J_{\ell}$ 
 evidently has rank $\ell$. 
 
It is a standard fact that the representation of a group on the functions on a coset space $G/H$ is identifiable with the induced representation $Ind_H^G {\bf 1}$, from the trivial representation of the stabilizing subgroup $H$. Because of the containments $H_{\ell - 1} \supset H_{\ell}$, it follows that the representation of $GL_n$ on the functions on $M_{nk}$ is contained in a certain number of copies of the induced representation $Ind_{H_k}^{GL_n}{\bf 1}$. 
 
 Frobenius reciprocity now tells us that a representation $\rho$ in $\widehat {GL}_n$ will appear
 in the functions on $M_{nk}$ if and only if it contains a fixed vector for $H_k$. This completes the proof of the proposition \Qed
 
\vskip .1 in

{\bf 9.2. The new spectrum.}

\vskip .1 in

The following description of isotypic components will be useful. Let $P_{k, n-k}$ denote the block upper triangular parabolic subgroup, whose block diagonal Levi component is $GL_k \times GL_{n-k}$. Then $H_k$ is a normal subgroup of $P_{k, n-k}$, and further we have
$P_{k,n-k} \simeq GL_k \cdot H_k$. For a representation $\nu$ of $GL_k$, let $\tilde \nu$ be the representation of $P_{k,n-k}$ pulled back from $\nu$ via the quotient mapping $P_{k,n-k} \rightarrow GL_k$. Clearly $\tilde \nu$ will be trivial on $H_k$.

\vskip .07 in

{\bf Proposition 9.2.1 (Formula for an isotypic component in the new spectrum).} Consider the joint
action of $GL_n \times GL_k$ on $GL_n/H_k$, with $GL_n$ acting by translations on the left, and
$GL_k$ (which normalizes $H_k$) acting by translations on the right. Fix an irreducible representation
$\nu$ of $GL_k$. Then the $\nu$ isotypic component of $L^2(GL_n/H_k)$ is the representation
$$
\left (Ind_{P_{k, n-k}}^{GL_n} \tilde \nu \right) \otimes \nu.
$$

\vskip .1 in 

{\it Proof.} The group $GL_k$ normalizes $H_k$, and the product $GL_k\cdot H_k$ is the parabolic subgroup $P_{k, n-k}$. 
The result follows from the standard formula for  $Ind_{P_{k,n-k}}^{GL_n} {\tilde \nu}.$  \Qed

\vskip .1 in 

Propositions 9.1.2 and 9.2.1 together imply,

\vskip .07 in

{\bf Corollary 9.2.2.} Let $\eta(\nu)$ be the irreducible representation of $GL_n$ in the new spectrum for the dual pair $(GL_n, GL_k)$ corresponding
 to an irreducible representation $\nu$ of $GL_k$. Then, 
 $$
\eta(\nu) < Ind_{P_{k,n-k}}^{GL_n} {\tilde \nu}.
 $$
 
 \vskip .1 in 
 
{\it Proof.} According to Proposition 9.1.3, the new spectrum consists of representations of $GL_n$ that, 
up to twisting by a character, appear in $L^2(GL_n/H_k)$, but do not have any fixed vectors for $H_{k-1}$. Since, by definition, $\eta (\nu)$ must lie in the $\nu$-isotypic component of the oscillator representation for the pair $(GL_k, GL_n)$ and must also lie in the
new spectrum for $GL_n$, the result follows from Proposition 9.2.1.  \Qed

 \vskip .1 in

Finally, the intrinsic criterion of Proposition 9.1.3 can assist in the characterization of the new spectrum (as described in \S 8.2) for the dual pair $(GL_n, GL_k)$ with for all $k$. 

Here is the formulation.

\vskip .07 in

{\bf Theorem 9.2.3 (The domain of the eta correspondence for $(GL_n,GL_k)$).} Assume $0 \leq k \leq n$. Then, a representation $\nu$ of $GL_k$ contributes, via the $(GL_n, GL_k)$ duality, to the new $k$-spectrum for $GL_n$ if and only if $\nu$ appears in the new $l$-spectrum of $GL_k$ for some $l \geq 2k - n$. 

\vskip .1 in 

{\it Proof of the necessity statement}. If $\nu$ appears in the new $l$-spectrum for $l < 2k-n$ then $\nu$ has a fixed vector for $H'_{2k - n - 1}$, the stabilizer of the first $2k - n - 1$ basis vectors in ${\bf F}_q^k$. This in turn implies that $\nu$ is contained in an induced representation $Ind_{P'_{2k-n-1, n-k+1}}^{GL_k} (\tilde \nu \otimes {\bf 1}_{n-k+1})$ from the maximal parabolic subgroup $P'_{2k-n-1, n-k+1} \simeq GL_{2k-n-1} \cdot H'_{2k - n - 1}$ of $GL_k$, and from a representation that is trivial on $H'_{2k - n - 1}$, tensored with an appropriate representation
$\tilde \nu$ of $GL_{2k - n - 1}$. Combining this with Corollary 9.2.2, we conclude that any candidate for $\eta(\nu)$ will be contained in an induced representation
$$
Ind_{\tilde P_1}^{GL_n} (\tilde \nu \otimes {\bf 1}_{n-k+1} \otimes {\bf 1}_{n-k}) = Ind_{\tilde P_1}^{GL_n} \tau_1,
$$
where $\tilde P_1 = P_{2k-n-1, n-k+1, n-k}$ is the upper triangular parabolic subgroup with Levi component $\tilde L_1 = GL_{2k - n - 1} \times GL_{n - k +1} \times GL_{n-k}$, and $\tilde \nu$, ${\bf 1}_{n-k+1}$, and ${\bf 1}_{n-k}$, are the factors of the inducing representation $\tau_1$ on the three factors of $\tilde L_1$. 

Consider the upper triangular parabolic $\tilde P_2 = P_{2k - n - 1, n-k, n - k+1}$, with Levi component $\tilde L_2 = GL_{2k - n - 1} \times GL_{n - k} \times GL_{n-k+1}$, and the representation 
$$
\tau_2 = \tilde \nu \otimes {\bf 1}_{n-k} \otimes {\bf 1}_{n-k+1},
$$
of $\tilde L_2$. It is clear that $\tilde P_1$ and $\tilde P_2$ are associate (namely, their Levi components are conjugated), and that the element that conjugates $\tilde L_1$ to $\tilde L_2$ by interchanging the last two factors will also conjugate $\tau_1$ to $\tau_2$.

The philosophy of cusp forms (see Section 10.1) tells us that $Ind_{P_1}^{GL_n} \tau_1$ and $Ind_{P_2}^{GL_n} \tau_2$ are equivalent. But since $P_2$ contains $H_{k-1}$, we see that $Ind_{P_2}^{GL_n} \tau_2$ will consist of representations that belong to the new $l$-spectrum for $l$ at most $k-1$. Thus, so also does $Ind_{P_1}^{GL_n} \tau_1$, which amounts to the statement of the necessity part of the theorem. \Qed

\vskip .1 in 

The main idea involved in our proof of the sufficiency part of Theorem 9.2.3 will be given in \S 10.

\vskip .1 in 

{\bf Remark 9.2.4.} In a recent preprint [Guralnick-Larsen-Tiep17] Guralnick, Larson and Tiep study representations from a point of view related to the notion of tensor rank 
of \S 6. They announced that, using different techniques, i.e., a substantial amount of character theory, they are also able to prove Theorem 9.2.4.

\vskip .2 in
 
{\bf 10. New Spectrum and Tensor Rank for $GL_n$ via the Philosophy of Cusp Forms}

\vskip .15 in 

In Theorem 9.2.3, we have seen an easy argument showing that only irreps from the set 
$(\widehat{GL}_k)_{new,\geq 2k-n}$, of irreps of $GL_k$ belonging to the new $l$-spectrum for $l$ at least $2k-n$, can give rise to irreps of $GL_n$ from the new $k$-spectrum, i.e., from the set $(\widehat{GL}_n)_{new,k}$. 

It is harder to show that each member of $(\widehat{GL}_k)_{new,\geq 2k-n}$ do give rise to uniquely defined representations of $GL_n$ of rank $k$. Our argument for this passes through a fairly precise description of the tensor rank of representations, and of the $\eta$ correspondence, 
which invokes  the ``philosophy of cusp forms" [Harish-Chandra70], [Howe-Moy86], and is  somewhat long. Here we will give a sketch of the main ideas. More details will appear in [Gurevich-Howe17]. 

In this approach, the basic case to understand is the {\it spherical principal series}, i.e., the representations of $GL_n$ that appear in $L^2(GL_n/B)$, where $B$ is a Borel subgroup; or in other words, in the action of $GL_n$ on the variety of complete flags (nested sequences of subspaces) inside ${\bf F}_q^n$. In turn, the spherical principal series can be understood in terms of the representation theory of the symmetric group $S_n$ [Green55], [Howe-Moy86], [Zelevinsky81].

Let us start with the basics of the philosophy of cusp forms.

\vskip .1 in 

{\bf 10.1. The philosophy of cusp forms.} 

\vskip .1 in

Recall that, a representation $\pi$ of $GL_n$ is called {\it cuspidal} if it does not contain a fixed vector for the unipotent radical of any parabolic subgroup.
Given this definition, it is easy to show that any irreducible representation is contained in a representation induced from a representation $\rho$ of a parabolic subgroup $P$ that 
\vskip .05 in

{\bf i)} is trivial on the unipotent radical $U_P$ of $P$, and

\vskip .05 in

{\bf ii)} is a cuspidal representation of the quotient $P/U_P$.

\vskip .05 in

\noindent Note that $P/U_P$ is reductive, and in fact, for the case of $GL_n$,
is a product of $GL_m$s for $m < n$,
so that the philosophy of cusp forms provides an inductive construction (in two senses) 
of all irreducible representations. 

The main ingredients needed to carry out this construction explicitly are

\vskip .05 in

{\bf a)} knowledge of the cuspidal representations of $GL_m$, and

\vskip .05 in

{\bf b)} decomposing the representations induced from cuspidal representations.

\vskip .05 in

\noindent With regard to a), the cuspidal representations of $GL_m$ are well known. They were already accounted for by Green  [Green55], and further insight came from the work of S. I. Gelfand [Gelfand70].  
We will not need to deal explicitly with this part of the theory here. 
With regard to b), there is a very general result that relates different induced-from-cuspidal
representations. Recall that two parabolic subgroups $P_1$ and $P_2$ are called
{\it associate} if their Levi components are conjugate. Given parabolic subgroups
$P_1$ and $P_2$, with cuspidal representations $\rho_1$ and $\rho_2$ respectively,
then the representations of $GL_n$ induced from the $\rho_a$, for $a = 1,2$ are either
equivalent, which happens when the $P_a$ are associate, and for some element 
$g$ in $GL_n$, we have $L_{P_2} = g(L_{P_1})g^{-1}$, and moreover 
$\rho_2 \simeq g \rho_1 g^{-1}$;  or they are completely disjoint - they have no irreducible
components in common (more precisely, the double cosets $P_1 g P_2$ that support
intertwining operators between the two induced representations are exactly 
the ones for which $g$ can be chosen to conjugate $L_{P_1}$ to $L_{P_2}$, 
and to be an equivalence of the representations $\rho_1$ and $\rho_2$ [Howe-Moy86]).
Thus, associativity classes of cuspidal representations of parabolic subgroups
define a partition of the unitary dual of $GL_n$ into disjoint subsets.

As noted already, for  $GL_n$, any Levi component is a product of $GL_m$s, with $m \leq n$.
The collection $D_P  = \{m_j\}$ of the sizes of the $GL$ factors of a Levi component define 
a partition of $n$. Up to conjugation, we can assume that $P$ consists of block upper triangular
matrices. Given this, we let $m_j$ be the size of the $j$-th block from the upper left corner 
of the matrices. Then the $m_j$ form a partition of $n$. 
We will refer to this as the {\it $P$ partition}. Also, we will say that a cuspidal
representation of $GL_m$ has {\it cuspidal size} $m$. Thus, the cuspidal sizes of a cuspidal
representation of a parabolic subgroup $P$ also define a partition, the same as the
$P$ partition. Up to associativity of parabolic subgroups in $GL_n$, we can arrange that
the block sizes $m_j$ of $P$,  equivalently, cuspidal sizes of $\rho$, decrease as $j$ increases. We then also associate to  the partition, and to $P$, a Young diagram $D_P$ [Fulton97], whose $j$-th row has length $m_j$.

If the block sizes are all equal to 1, then the parabolic is (conjugate to) the Borel subgroup 
$B$ of upper triangular matrices. We will refer to constituents of representations induced from
characters of $B$ as the {\it split principal series}. Constituents of the representation 
induced from the trivial character of $B$ will be referred to as {\it spherical principal series}.
Any representation induced from a (one dimensional) character of a parabolic subgroup
will have constituents all belonging to the split principal series.

On the other hand, if the block sizes are all greater than 1, then no constituent will contain 
a character of any parabolic subgroup. We will refer to the components of representations 
induced from cuspidal representations of parabolics with block sizes of 2 or larger (i.e., no
blocks of size 1) as {\it unsplit} representations.

The general representation is gotten by combining unsplit representations and split 
principal series. More precisely, given the parabolic $P_D$ as described above,  
let $P_{r,s} \supset P_D$  be the maximal parabolic subgroup, 
with Levi component $GL_r \times GL_s$, where $GL_r$ contains 
all the blocks of $P_D$ of size greater than 1, and $GL_s$  contains all the blocks of $P_D$ of size 1. 
Then a constituent of a representation of $GL_r$ induced by a cuspidal representation 
of $GL_r \cap P_D$ will be  an unsplit representation of $GL_r$. On the other hand, 
a constituent of the representation of $GL_s$ induced from a cuspidal representation of 
$GL_s \cap P_D$ (which is the upper triangular Borel subgroup of $GL_s$)
will be a split principal series of $GL_s$. We call these constituents $\sigma_r$ 
and $\tau_s$ respectively. Since both $GL_r$ and $GL_s$ are quotients of
$P_{r,s}$, we can regard both $\sigma_r$ and $\tau_s$ as representations of $P_{r,s}$.
Because $GL_r \times GL_s$ is also a quotient of $P_{r,s}$, the tensor product
$\sigma_r \otimes \tau_s$ will also define a representation of $P_{r,s}$, and this representation
will be irreducible. Now if we look at the induced representation
$$
\pi = \pi (\sigma_r, \tau_s) = Ind_{P_{r,s}}^{GL_n} (\sigma_r \otimes \tau_s),
$$
the philosophy of cusp forms tells us that $\pi$ is irreducible, and further says
that all irreducible representations of $GL_n$ arise in this way, (including the cuspidal 
representations, which correspond to the situation when $P_D = GL_n$).

\vskip .15 in

There is a second partition, that permits a more refined understanding of the split principal series,
essentially reducing it to understanding the spherical principal series. 
A character $\chi$ of the standard upper triangular Borel subgroup $B$ is given by a collection 
$\psi_j$, $1 \leq j \leq n$,  of characters of $GL_1 \simeq {\bf F}_q^{\times}$, 
where $\psi_j$ is the restriction of the character $\chi$  to the $j$-th diagonal entry of an element of $B$. 
Up to associativity, these characters can be reordered as desired. 
Thus, we may arrange that all diagonal entries $j$ for which the $\psi_j$
are equal to a given character of $GL_1$ are consecutive. Given this, we can consider 
a block upper triangular parabolic subgroup  $P_{\chi}$ such that, in each diagonal block, the characters 
$\psi_j$ are equal, and the $\psi_j$ contained in different diagonal blocks are different.
We may also assume when convenient that the sizes of these blocks decrease from top to bottom.
This associates a well-defined Young diagram $D_{\chi}$, and a parabolic subgroup $P_{\chi}$,
to a given character $\chi$ of $B$. 

 Let the $i$-th block from the top of $P_{\chi}$  be $GL_{m_i}$. 
Let $\tau$ be a constituent of the representation of $P_{\chi}$ induced from the character $\chi$ of $B$.
Then, the philosophy of cusp forms tells us that:

\vskip .05 in 

\noindent {\bf Scholium 10.1.1.} we have,

\vskip .05 in 

{\bf a)} the representation $Ind_{P_{\chi}}^{GL_n} \tau$, induced from $P_{\chi}$ to $GL_n$, will be irreducible;

\vskip .05 in

{\bf b)} $\tau \simeq \otimes_i \tau_i$, where $\tau_i$ is a constituent of the representation
of $GL_{m_i}$ induced from the restriction of $\chi$ to $B\cap GL_{m_i}$; and

\vskip .05 in 

{\bf c)} this process gives a bijection  to the set of constituents of $Ind_B^{GL_n} \chi$ 
from the set of constituents of $Ind_B^{P_{\chi}} \chi$,
and this set is the product of the sets of constituents of the $Ind_{B \cap GL_{m_i}}^{GL_{m_i}} \chi$.
\vskip .03 in 

Moreover, because of the way $P$ was defined, each representation 
$Ind_{B \cap GL_{m_i}}^{GL_{m_i}} \chi$ has the form
$$
(\chi_i \circ \det)\otimes(Ind_{B \cap GL_{m_i}}^{GL_{m_i}} {\bf 1}),
$$
where $\chi_i$ here indicates the common character of $GL_1$ assigned to the diagonal entries 
of $B \cap GL_{m_i}$, and as usual, $\det$ is the determinant homomorphism from $GL_{m_i}$
to $GL_1$. This means that each constituent $\tau_i$  of each $Ind_{B \cap GL_{m_i}}^{GL_{m_i}} \chi$
has the form $\tau_i \simeq (\chi_i\circ det) \otimes \tilde \tau_i$, where $\tilde \tau_i$ is a constituent
 of the spherical principal series for $GL_{m_i}$.
 
 \vskip .15 in 
 
{\bf 10.1.2. Standard realization.} Combining the decomposition of Scholium 10.1.1 with the split-unsplit decomposition noted above,
 we see that every representation $\pi$ of $GL_n$ is induced irreducibly from a parabolic subgroup 
 $Q_{\pi} = Q$ with Levi component $L_{Q} \simeq \prod_j GL_{m_j}$, and a representation
 $\rho_{\pi} = \rho$ of $Q$, trivial on the unipotent radical $U_Q$, and given on $L_Q$ by a tensor product $\rho \simeq \otimes_j \sigma_j$, where:
 
\vskip .05 in
 
{\bf i)} $\sigma_1$ is an unsplit representation of $GL_{m_1}$; and
 
\vskip .05 in
 
{\bf ii)} for $j \geq 2$, $\sigma_j \simeq (\chi_i \circ \det) \otimes \nu_j$, where the $\chi_j$ are distinct characters of $GL_1$, 
 
\ \ \ \ and $\nu_j$ is an irreducible spherical principal series representation of $GL_{m_j}$.
 
 \vskip .07 in 
\noindent  We will call this realization of $\pi$ the {\it standard realization}, and will refer to $Q_{\pi}$ and
 the constituents $\sigma_j$ of $\rho_{\pi}$ as the {\it standard data} for realizing $\pi$.
 
 \vskip .1  in

 {\bf 10.2. The spherical principal series and tensor rank.}
 
 \vskip .1 in 
 
Since most of the constituents $\sigma_j$ of the representation $\rho_{\pi}$ (described just above) are spherical principal series twisted by characters, the next step in understanding $\pi$ is to focus on the spherical principal series (SPS). We will see that an SPS representation is of tensor rank $k$ exactly when it appears in the new $k$-spectrum, and moreover, this common $k$ can be read easily from the parametrization of the SPS representation by a Young diagram.

\vskip .1 in

The spherical principal series of $GL_n$ have been studied extensively [Green55], [Howe-Moy86], [Zelevinsky81]. Our principal reference for this account will be [Howe-Moy86]. 

The spherical principal series can be helpfully studied through the family of induced representations $I_P = Ind_P^{GL_n} {\bf 1}$, for all parabolic subgroups, where ${\bf 1}$ is the trivial representation of the relevant group. 
Up to conjugation, it is enough to consider the parabolics that contain the Borel subgroup 
of upper triangular matrices, i.e., the block upper triangular parabolic subgroups.
Also, it is standard that if $P$ and $P'$ are associate parabolics, then the representations
$Ind_P^{GL_n}{\bf 1}$ and $Ind_{P'}^{GL_n}{\bf 1}$ are equivalent. Thus, we can select a representative
from each associativity class of parabolics. We do this in the usual way, by requiring/assuming 
that the block sizes $m_i$ of the diagonal blocks $GL_{m_i}$ of $P$, listed from top to bottom,
are decreasing with increasing $i$. This again gives us a partition (our third partition) of $n$,
with an associated Young diagram $D_P$. (We note that, if all of the above discussion is
referenced to a fixed original $n$, then the successive partitions we have been describing are
partitions of parts of the preceding partition.) 

Let $P = P_D$ be a parabolic as in the preceding paragraph, with blocks whose sizes 
decrease down the diagonal, associated to the Young diagram $D$.
Consider the induced representation 
$$
I_D = I_{P_D} = Ind_{P_D}^{GL_n} {\bf 1}.
$$
All the constituents of $I_D$ are spherical principal series representations. We can be
somewhat more precise. Recall that the  set of isomorphism classes of representations of a group $G$
form a free abelian semigroup on the irreducible representations, and as such, has a natural order structure, given by dominance of all coefficients in the expression of a given representation as a sum of irreducibles. Equivalently, one representation dominates a second one if the second one is equivalent to a subrepresentation of the first.

The set of partitions also has a well-known order structure, the {\it dominance order} [Ceccherini-Silberstein-Scarabotti-Tolli10]. We know the following facts (see [Howe-Moy86]).

\vskip .1 in

\noindent {\bf Scholium 10.2.1.} We have,

\vskip .05 in 

{\bf i)} The map $D \mapsto I_D$ is order preserving from the set of partitions
with its reverse dominance order to the semigroup of spherical principal series representations of $GL_n$.

\vskip .05 in 

{\bf ii)} $I_D$ contains an irreducible constituent $\nu_D$ with multiplicity one, and with the property that $\nu_D$ is not contained in any $I_{D'}$ with $D' > D$ in the dominance order.

\vskip .05 in

{\bf iii)} The $\nu_D$ constitute the full collection of irreducible spherical principal series representations.

\vskip .15 in

\noindent {\bf Remarks 10.2.2.} We have,

\vskip .05 in  

{\bf a)} The representation $\nu_D$ can be distinguished by its dimension: it is the only constituent of $I_D$ whose dimension, as a polynomial in $q$, has the same degree as the cardinality of $GL_n/P_D$. It also can be distinguished as being the only constituent that contains certain characters of 
$\overline U_D$, the unipotent radical opposite to the unipotent radical of $P$. For this, see Section 10.5 below.

\vskip .05 in 

{\bf b)} The facts stated in Scholium 10.2.1 are parallels of similar facts for the symmetric group $S_n$. 
For a given partition $D$ of $n$, let $S_D$ denote the stabilizer of $D$ in $S_n$, 
and let $\nu^o_D$ be the irreducible representation of $S_n$ associated to the partition $D$. 
Also, let $I^{\circ}_D$ be the representation of $S_n$ induced from the trivial representation of $S_D$.
Then the analogs of facts i), ii) and iii) are valid also for the $I^{\circ}_D$.
In particular, $\nu^o_D$ is contained in $I^o_D$ with multiplicity one.

\vskip .05 in

{\bf c)} Moreover, the Bruhat decomposition for $GL_n$ [Borel91] implies that the intertwining number of 
$I_{D_1}$ and $I_{D_2}$ for two partitions $D_1$ and $D_2$ is equal to
the intertwining number of $I^{\circ}_{D_1}$ and $I^{\circ}_{D_2}$. In this sense,
the description of the spherical principal series representations of $GL_n$ 
is essentially the same as the representation theory of the symmetric group.

\vskip .05 in 

{\bf d)} For example, if $\nu^o_D = \sum_E c_{DE} I^o_E$, then also 
$\nu_D =  \sum_E c_{DE} I_E$. That is, the expression of the irreducible spherical 
principal series representations $\nu_D$ as linear combinations of the induced representations 
$I_D$ involves exactly the same coefficients as the parallel expressions for $S_n$.

\vskip .15 in

It is easy to determine the new $k$-spectrum to which a spherical principal series representations belongs and its tensor rank. 

\vskip .1 in

{\bf Proposition 10.2.3.} For a partition $D$ of $n$, normalized as above, the spherical principal series representation $\nu_D$ belongs to the new $k$-spectrum and is of tensor rank $k$ for $k=n-d_1$, where $d_1$ is the longest row of the associated Young diagram.
\vskip .1 in 

For the complete proof of Proposition 10.2.3, we refer to [Gurevich-Howe17]. Here we mention that the argument relies on the Pieri rule for representations of $S_n$, and its parallel for the spherical principal series. Here is the Pieri Rule for the symmetric group [Ceccherini-Silberstein-Scarabotti-Tolli10].

\vskip .1 in

{\bf Scholium 10.2.4.} Let $\hat D$ be a partition of $k$. Then the induced representation
$$
Ind_{S_k \times S_{n-k}}^{S_n} (\nu^o_{\hat D}\otimes {\bf 1}),
$$
where ${\bf 1}$ here is the trivial representation of $S_{n-k}$, is a multiplicity free sum 
of representations $\nu^o_D$ of $S_n$, where the partition $D$ satisfies the conditions that: 

i) $\hat D \subset D$ (in the sense of containment of Young diagrams); and

ii) $D - \hat D$ is a skew row.
\vskip .1 in

For the notions of containment of Young diagrams and of skew row see [Ceccherini-Silberstein-Scarabotti-Tolli10]. 

\vskip .1 in

{\bf 10.3. New spectrum and tensor rank for general irreps of $GL_n$.}

\vskip .1 in

Combining Scholium 10.2.4 with the description of a representations by means of its standard data, as described at the end of \S 10.1,  allows us to  determine when a general irrep $\pi$ of $GL_n$ appears in the new spectrum and what is its tensor rank. Indeed, 

\vskip .1 in

{\bf Proposition 10.3.1.} Let $\pi$ be an irreducible representation of $GL_n$. 
Let $Q_{\pi}$ be the parabolic subgroup associated to $\pi$ in Paragraph 10.1.2, and let  $\{\sigma_j\}$ be the standard data for realizing $\pi$. As described there, for $j \geq 2$, the $\sigma_j \simeq (\chi_j \circ \det) \otimes \nu_j$, where $\nu_j$ is a spherical principal series for $GL_{m_j}$. Let us assume the $\chi_1$ is he trivial character and $\nu_1$ the corresponding SPS representation. Then, 

\vskip .05 in

{\bf (a)} $\pi$ appears in the new $k$-spectrum for $k$ equal $n$ minus the size of the longest row of the Young diagram parameterizing $\nu_1$.

\vskip .04 in

and

\vskip .04 in

{\bf (b)} the tensor rank of $\pi$ is equal to $n$ minus the size of the longest row over all the Young diagrams that parameterize the $\nu_j, \,\, j\geq 2$.

\vskip .1 in

For the details of proof, we again refer to [Gurevich-Howe17].

\vskip .1 in

{\bf 10.4. Explicit description of the eta correspondence for the pair $(GL_n, GL_k)$.}

\vskip .1 in

Proposition 10.3.1, along with the results of Sections 10.1 and 10.2, allow us to give a succinct description of the $\eta$ correspondence for the pair $(GL_k, GL_n)$. 

\vskip .05 in

Take a representation $\nu$ of $GL_k$, of the new $l$-spectrum for $l$ at least $2k-n$. Form the standard block upper triangular parabolic subgroup $P_{k, n - k} \simeq GL_k \cdot H_k$. The group $H_k$ (which pointwise fixes the first $k$ coordinate subspace of ${\bf F}_q^n$) is a normal subgroup of $P_{k, n-k}$. Extend $\nu$ to a representation of $P_{k,n-k}$, that we will also denote $\nu$, by letting $H_k$ act trivially. Form the induced representation
$$
Ind_{P_{k,n-k}}^{GL_n} \nu.
$$
Then $\eta (\nu)$, that we want to describe, will be the unique constituent of $Ind_{P_{k,n-k}}^{GL_n} \nu$ from the new $k$-spectrum.

To get an explicit formula of $\eta (\nu)$, and to be sure that it determine a well-defined representation, we invoke the philosophy of cusp forms, use it to give a more refined description of $\nu$, and then to describe $\eta (\nu)$.

First, suppose that we start with a spherical principal series representation
$\hat \nu = \nu_{\hat D}$ of $GL_k$, where $\hat D$ is a partition of $k$. The condition that $\hat \nu$ be from the new $l$-spectrum for $l$ at least $2k-n$, implies that no part of $\hat D$ has size greater than $n-k$.
As just described, $\eta (\hat \nu)$ is a constituent of 
$Ind_{P_{k, n - k}}^{GL_n}  \hat \nu$.
The decomposition of this induced representation is described by the parallel of Scholium 10.2.4
for the spherical principal series of $GL_n$. It consists of spherical principal series
representations $\nu_D$, for partitions $D$ of $n$, such that
\begin{itemize}
\item the Young diagram of $D$ contains the diagram for $\hat D$, and 
\item the difference $D - \hat D$ is a skew row. 
\end{itemize} 
It should not be hard for the reader to conclude that the only such $D$ that has no rows 
longer than $n-k$ is the one $\eta (\hat D)$ corresponding to adjoining one part of size $n - k$ 
to the parts of $\hat D$. We then have 
$$
\eta (\nu_{\hat D})=\nu_{\eta (\hat D)},
$$
where $\nu_{\eta (\hat D)}$ is the spherical principal series representations of $GL_n$ attached to $\eta (\hat D)$ (see Scholium 10.2.1.). 

\vskip .05 in 

Now let $\rho$ be any representation of $GL_k$, from a new spectrum which is at least $2k-n$. Using the philosophy of cusp forms (see Section 10.1), we can realize $\rho$ as an induced representation
$$
\rho = Ind_{P_{\ell, k - \ell}}^{GL_k} (\sigma \otimes \hat \nu),
$$
where:

i)  $P_{\ell, k - \ell}$ is the standard upper triangular parabolic subgroup of $GL_k$, 
with Levi component $GL_{\ell} \times GL_{k - \ell}$; 

ii) $\sigma$ is an irreducible representation of $GL_{\ell}$ whose standard data includes no spherical principal series; and

iii) $\hat \nu$ is an irreducible spherical principal series of $GL_{k - \ell}$. 

\noindent Let $\hat D$ be the partition of $k-\ell$ that parametrizes $\hat \nu$. Then the new spectrum condition implies that $\hat D$ contains no part larger than $n-k$. 

By transitivity of induction, we know that 
$$
Ind_{P_{k, n-k}}^{GL_n} \rho \simeq Ind_{P_{\ell, k - \ell, n-k}}^{GL_n} (\sigma \otimes \hat \nu \otimes {\bf 1}),
$$
where $\bf 1$ here denotes the trivial representation of the third factor of the Levi component $GL_{\ell} \times GL_{k - \ell} \times GL_{n-k}$ of $P_{\ell, k - \ell, n-k}$.

Let $P_{\ell, n - \ell} \supset P_{\ell, k - \ell, n-k}$ be the parabolic subgroup with Levi component 
$GL_{\ell} \times GL_{n - \ell}$, and such that 
$GL_{n - \ell} \cap P_{\ell, k - \ell, n-k} = P_{k - \ell, n - k}$ is the parabolic subgroup of $GL_{n-\ell}$ with Levi  component $GL_{k - \ell} \times GL_{n - k}$, consisting of the last two factors of the Levi component
of the Levi of $P_{\ell, k - \ell, n-k}$. Then, the representation
$Ind_{P_{\ell, k - \ell, n-k}}^{P_{\ell, n - \ell}}  (\sigma \otimes \hat \nu \otimes {\bf 1})$,
will factor to a representation of $GL_{\ell} \times GL_{n - \ell}$, and this representation will have the form 
$$
\sigma \otimes Ind_{P_{k - \ell, n - k}}^{GL_{n - \ell}} (\hat \nu \otimes {\bf 1}).
$$
Again by the $GL_{n-\ell}$ analog of Scholium 10.2.4, we know that
\begin{itemize}
\item $Ind_{P_{k - \ell, n - k}}^{GL_{n - \ell}} \hat \nu \otimes {\bf 1}$ decomposes multiplicity freely,
and
\item with a unique component $\eta (\hat \nu) = \nu_{\eta (\hat D)}$, a spherical principal series of $GL_{n - \ell}$,
such that $\eta (\hat D)$ has no parts greater than $n - k$ (so that $\eta (\hat \nu)$ is a member of the new $\ell + k - n$ spectrum).
\end{itemize}
If we now form 
$
Ind_{P_{\ell, n - \ell}}^{GL_n} (\sigma \otimes \eta (\hat \nu)),
$
then since $\sigma$ has no spherical principal series in its standard data,
we know that this representation will be irreducible, and a member of the new $k$-spectrum. It is $\eta (\rho)$:
$$
\eta (\rho) = Ind_{P_{\ell, n - \ell}}^{GL_n} (\sigma \otimes \eta (\hat \nu)). \eqno (10.4.1)
$$
This completes the explicit description of the eta correspondence for the pair $(GL_n, GL_k)$.

\vskip .2 in 

{\bf 10.5. Generalized Whittaker rank}

\vskip .1 in

In this subsection, we would like to offer a complement to the theme of $U$-rank. This 
is not necessary for understanding tensor rank, but perhaps it can play a role in 
providing suitable spectral analogs that extend $U$-rank. Also, it provides an argument,
independent of the Bruhat decomposition, for Scholium 10.2.1 ii).

We analyzed  the $\eta$ correspondence for $(GL_n, GL_k)$ in the range $k < \lfloor{n \over 2}\rfloor$ 
using the notion of $U$-rank, as discussed in \S 5 above. This invoked the spectral behavior of 
representations under restriction to the unipotent radicals of maximal parabolic subgroups.
Since these unipotent radicals are abelian, this analysis comes down to studying the characters 
of these groups. In the literature of representation theory, considerable attention has been paid 
to the appearance in representations of characters of maximal unipotent subgroups [Jacquet67], [Kostant78], [Shalika74]. In this context  also, there is a natural notion of rank, but it has not been emphasized; 
rather the interest has been in studying representations that allow an eigenvector (in some distributional sense) for the maximal unipotent subgroup, with respect to a character of maximal rank. Such vectors are known in the literature as {\it Whittaker vectors}.

For $GL_n$, and for any parabolic subgroup $P$, with unipotent radical $U_P$, there is 
a straightforward way to define the rank for characters of $U_P$. Indeed, given a Levi component
$L_P$ for $P$, the center of $L_P$ acts on (the Lie algebra of) $U_P$, and the eigenspaces for 
this action are standard matrix subgroups, in the sense of \S3.1. Moreover, the complement to the 
commutator subgroup/subalgebra of $U_P$ is spanned by a sum of these eigenspaces. 
Therefore, given a character of $U_P$, the restriction of this character to any one of the eigenspaces
has a well-defined rank.

Let $D = \{d_1 \geq d_2 \geq d_3 \ . \ . \ . \geq d_{\ell}\}$ be a partition of $n$, and let $P_D$
be the corresponding block upper triangular parabolic subgroup, with blocks of size $d_i$ along
the diagonal. Let $L_D$ be the Levi component of  $P_D$ consisting of block diagonal matrices.
Let $U_D$ be the unipotent radical of $P_D$, and let $\overline U_D$ be the opposite unipotent
subgroup. Then $\overline U_D$ is normalized by $L_D$, which therefore acts by conjugation
on the characters of $\overline U_D$.

\vskip .11 in

{\bf Proposition 10.5.1.} We have,

\vskip .05 in

{\bf a)} The characters of $\overline U_D$ of maximal rank form a single
orbit under conjugation by $L_D$.

\vskip .05 in

{\bf b)} A character of maximal rank of $\overline U_D$ appears with multiplicity one in the 
induced representation $I_D =  Ind_{P_D}^{GL_n} {\bf 1}$.

\vskip .05 in

{\bf c)} Let $D' > D$ be a partition that (strictly) dominates $D$. Then the representation $I_{D'}$, when
restricted  to $\overline U_D$, does not contain any characters of maximal rank.

\vskip .07 in

{\it Proof.} a) Consider two adjacent diagonal blocks, of sizes $d_i \geq d_{i+1}$. 
The (standard) group $U(i)$ of size $(d_i, d_{i+1})$ of unipotent matrices, 
whose only non-zero off-diagonal entries are in the  $d_i \times d_{i+1}$ rectangular block 
directly below the diagonal block  of size $d_i$ and directly to the left
 of the diagonal block of size $d_{i+1}$,
 is conjugated to itself by the action of $GL_{d_i} \times GL_{d_{i+1}}$.
As discussed in \S 3.1, matrices with entries in this block may be regarded as maps from ${\bf F}_q^{d_{i}}$
to ${\bf F}_q^{d_{i+1}}$, with the conjugation action being multiplication on the right by $GL_{d_i}$ 
and on the left by $GL_{d_{i+1}}$. The elements of $U(i)$ of maximal rank (i.e., rank $d_{i+1}$)
define  surjective maps, and any two such maps differ by multiplication on the right 
by an element of $GL_{d_i}$.  Dualizing, we see through considerations explained in \S2.1,
 that the same is true for characters of maximal rank of $U(i)$. Because of this, we see that, 
 if we start at the bottom of the matrix, and work our way to the top, that by conjugating
 by $L_D$, we can transform any character of maximal rank of $\overline U_D$  to any other. 
 
 \vskip .05 in
 
 Remark: The analogous result for a parabolic subgroup whose diagonal blocks are not
 ordered by size could be false.
 \vskip .05 in

 b) Now consider the representation $I_D = Ind_{P_D}^{GL_n} {\bf 1}$. The restriction of
 $I_D$ to $\overline U_D$ consists of a sum of permutation representations, one for each
 double coset $\overline U_D g P_D$, for $g$ in $GL_n$. Thanks to the Bruhat decomposition [Borel91], we may select $g$ to have the form $g = hw$, where $h$ is in $L_D$, and $w$ is an element of the symmetric group of permutation matrices, which is the Weyl group for $GL_n$. Since $h$ normalizes $\overline U_D$, and we are dealing with a property that is invariant under automorphisms of $\overline U_D$, we may ignore the effect of $h$, and just deal with the double coset $\overline U_D w P_D$. The action of $\overline U_D$ on the given double coset $\overline U_D wP_D$ is just the action by left translations on the coset space $\overline U_D/(\overline U_D \cap wP_Dw^{-1})$.
 
 In the action of $\overline U_D$ on any given double coset, 
 a given character of $\overline U_D$  will appear with multiplicity 1 or 0. Since the action
 of $\overline U_D$ on $\overline U_D w P_D$ is equivalent to the action by left
 translations on the coset space $\overline U_D/(\overline U_D \cap wP_Dw^{-1})$, 
 a given character  $\phi$ of $\overline U_D$ will appear with multiplicity 1
  if and only if $\overline U_D \cap wP_Dw^{-1}$
 is contained in the kernel of $\phi$. If $w = I$, the identity element, then 
 $\overline U_D \cap gP_Dg^{-1} = \overline U_D \cap P_D = \{I\}$, the identity element,
 so the action of $\overline U_D$ on this double coset is the regular representation
 of $\overline U_D$, and the multiplicity of any character of $U_D$ in this representation is 1. 
 Thus, we want to show that the multiplicity of any maximal rank character in the
 action on any other double coset is 0.
 So consider some other permutation $w$, and consider 
 $\overline U_D \cap wP_Dw^{-1} = w((w^{-1} \overline U_D w) \cap P_D)w^{-1}$. 
 
 Let $E_{ab}$ be the matrix units of the $n \times n$ matrices - the matrices 
 with exactly one non-zero entry, in the $a$-th row and $b$-th column, which is a 1. 
 For $a \neq b$, set $$\Gamma_{ab} = \{ I + \alpha E_{ab}: \alpha \in {\bf F}_q\}.$$
 Then $\Gamma_{ab}$ is a one parameter group of transvections.
 We can compute that 
 $$w^{-1} \Gamma_{ab}w = \Gamma_{w^{-1}(a)w^{-1}(b)},$$
 where, on the right side of this equation, $w$ is indicating the permutation of indices defined by $w$.
 
 Note that we are free to multiply $w$ on the left and the right by permutations
 that preserve the partition $D$. Thus, all that matters about $w$ is how many indices
 in each block it distributes to a given other block, not which particular indices
 in a given block go to which indices in another  given block. 
 
 Consider the smallest index $a$ that gets moved by $w$ from its block to a lower block.
 (Note that all the blocks above the block of $a$ must then be preserved by $w$. For, if $a$
 does not belong to the first/highest block, then all the indices in this block must be 
 mapped to other indices in this block; so the block is preserved. This argument extends
 until we reach the block containing $a$.) This implies that, up to multiplication by an element of $P_D$,
  we can take $w$ to act as the identity on these blocks. For purposes of this argument, 
  these blocks might just as well not be there.
 
 One of the one-parameter  groups $\Gamma_{c\ell}$ will belong to the group $U(i)$ defined in
 the proof of part a) if and only if 
 
 $$\sum_{k \leq i-1} d_k \leq \ell \leq \sum_{k \leq i} d_k \leq c  \leq \sum_{k \leq i+1} d_k.$$  
 
 If $c = w(a)$ is in a lower block than $a$, then $w^{-1}(c) = a$ is in a higher block than $c$. 
 Indeed, from our choice of $a$, we know that it is in the highest possible block that it not preserved by $w$. 
 Hence, $w^{-1}(\ell)$ cannot be in a higher block than $a$. It follows that 
 $w^{-1} \Gamma_{c\ell} w = \Gamma_{aw^{-1}(\ell)}$ is in $P_D$, and therefore, a character of 
 $\overline U_D$ that appears in the functions on $\overline U_D w P_D/P_D$ must be trivial on
 $\Gamma_{c\ell}$.  This argument is valid for any $\ell$ in the range defined by $U(i)$. 
 But a character of maximal rank of $U(i)$ must be non-trivial on one of the $\Gamma_{c\ell}$. 
 Thus, a character of maximal rank of $\overline U_D$, since its restriction to $U(i)$ must have 
 maximal rank, cannot occur in $\overline U_D w P_D/P_D$, and statement b) is established.
 
 The argument for part c) of the proposition is essentially the same as the argument for b). What
 is new about part c) is that even the identity double coset is disqualified, since the condition that 
 $D'$ dominates $D$ implies that, for some $k$, the sum $\sum_ {i \leq k} d'_i$ will be strictly larger
 than  $\sum _{I \leq k} d_i$. For this $k$, the top row of $U(k)$ will then be contained in $P_{D'}$,
 so no non-trivial character of this row can appear in the action of $U_D$
 on $U_D/ (U_D \cap P_{D'})$, which means that no maximal rank character of $U_D$ can occur.  \Qed
 
 \vskip .15 in 
 
{\bf Remark 10.5.2.} Proposition 10.5.1 implies that the constituent of $I_D =  Ind_{P_D}^{GL_n} {\bf 1}$ that contains a character of maximal rank for $\overline U_D$ will occur with multiplicity one, and that it will not occur in $I_{D'}$ for any $D'$ that strictly dominates $D$. In other words, this is the constituent  $\nu_D$ identified in Scholium 10.2.1.

\vskip .2 in
 
{\bf 11. On the agreement between the various notions of ranks for $GL_n$}

\vskip .15 in

We don't know yet to show that in general all notion of ranks agree in the range of small $U$-rank representations (Conjecture 6.5.2). However, for the group $GL_n$ we can show that this is the case at least for representations of $U$-rank $k<\frac{n}{4}$. As a result, we are able to conclude (see Theorem 11.4) that, for $GL_n$ , asymptotic rank and tensor rank agree for all representations.    

We will use here the intrinsic criterion given in Proposition 9.1.3 for determining tensor rank.

We begin with a somewhat curious weaker result.

\vskip .1 in  
 
{\bf Proposition 11.1.} A representation of $GL_n$ of $U$-rank $k < \lfloor{n \over 2}\rfloor$ has tensor rank at most $2k$.

\vskip .1 in

{\it Proof.} This argument involves working with a certain collection of standard unipotent subgroups of $GL_n$. We will introduce 
the actors here before getting down to the details.

In this paragraph, ``matrix" means ``$n \times n$ matrix". As in \S10, $E_{ij}$ indicates the matrix unit with an entry of $1$ in the $i$-th row and $j$-th column, and zeroes everywhere else. Each of the groups we want to define here consists of matrices gotten by adding the identity matrix (denoted $I$) to a linear combination of certain matrix units.

We start with the group ${\mathcal  U}_1$ given by 
$$
{\mathcal  U}_1 = I + span\, \{E_{ij}: \, j \leq k < i \leq n \}.
$$
The group ${\mathcal  U}_1$ is the unipotent radical of the standard parabolic $P_{k, n-k}$ with diagonal blocks of size $k$ and $n-k$.

\vskip .05 in

Next, we define the group 
$$
{\mathcal  U}_2= I + span\, \{E_{ij}: \, j \leq k;\,\,\,  n - k < i \leq n \}.
$$
This is the standard subgroup  of $GL_n$ with matrix units in the upper right $k\times k$ corner. In the terminology of \S3.1, ${\mathcal  U}_2$  has size $(k,k)$. 

\vskip .05 in

Now we look at the group 
$$
{\mathcal  U}_3= I + span\,  \{E_{ij}: \,j \leq  k < i \leq n-k\},
$$
and we note that 
${\mathcal  U}_3$  and ${\mathcal  U}_2$ are complementary subgroups of ${\mathcal  U}_1$, so that ${\mathcal  U}_1 \simeq{\mathcal  U}_2 \times{\mathcal  U}_3$. This is a decomposition of ${\mathcal  U}_1$ of the type considered in Formula (2.1.3).

\vskip .05 in

We also define the group 
$$
{\mathcal U}_4 = I + span\{E_{n-k,k+1}\}
$$
which is a one-parameter group with non-zero entry in the upper right corner of the middle $(n - 2k)\times (n - 2k)$ diagonal block of matrices.

\vskip .05 in

In addition, we consider the group
$$
{\mathcal U}_5= I + span \{E_{ij}:\,  j \leq k+1;\,\,\, n-k  \leq i \},
$$
namely, the standard subgroup in the upper right $(k+1) \times (k+1)$ corner. It contains both ${\mathcal U}_2$ and ${\mathcal U}_4$, and it is of size $(k+1, k+1)$.

\vskip .05 in

Finally, we denote by $\widetilde H_{2k}$ the group of matrices $T$ such that $T - I$ has its first $k$ columns and last $k$ columns equal to 0. Note that  $\widetilde H_{2k}$
\vskip .03 in 
$-$ is conjugate to $H_{2k}$;
\vskip .03 in 
$-$ normalizes ${\mathcal U}_1$;
\vskip .03 in
$-$ satisfies ${\mathcal U}_1 \cap \widetilde H_{2k} = {\mathcal U}_3$.

\vskip .07 in 

Note that we can find a character $\psi$ of ${\mathcal U}_1$ that has rank $k$, but is trivial on ${\mathcal U}_3$. Thus, the restriction $\psi_{|{\mathcal U}_2}$ will also have rank $k$. Observe also  that, since $\psi$ is trivial on ${\mathcal U}_3$ it  will be invariant under conjugation by $\widetilde H_{2k}$.

Now, take an irreducible representation $\rho$ of $GL_n$ which has $U$-rank $k$. Then the $\psi$-eigenspace of $\rho_{|{\mathcal U}_1}$ 
will be non-trivial. Denote this space by $\rho_{\psi}$. Since $\psi$ is invariant under conjugation by $\widetilde H_{2k}$, the space $\rho_{\psi}$
will be invariant under $\rho_{|\widetilde H_{2k}}$. Thus, it defines a representation of
$\widetilde H_{2k}$. We call this representation $\sigma$.

\vskip .05 in 

{\bf Claim.} The representation $\sigma$ must be trivial on the commutator subgroup of $\widetilde H_{2k}$.

\vskip .05 in

Note that the above claim implies that $\sigma$ is a sum of characters, which in its turn implies that 
 $\rho$ is of tensor rank at most $2k$.
 
 \vskip .03 in
 
 Let us verify the claim. Indeed, if $\sigma$ is not trivial on the commutator subgroup of $\widetilde H_{2k}$, then in particular,
 it must be not trivial on the one-parameter subgroup ${\mathcal U}_4$. This is because the normal 
 subgroup of $\widetilde H_{2k}$ generated by ${\mathcal U}_4$ is the full commutator subgroup 
 of $\widetilde H_{2k}$.  Since ${\mathcal U}_4$ is commutative, we can decompose the space $\rho_{\psi}$
 into eigenspaces for ${\mathcal U}_4$. Let $\rho^1_{\psi}$ be one of these eigenspaces, corresponding
 to a non-trivial character $\lambda$ of ${\mathcal U}_4$.
  
  The group ${\mathcal U}_5$ contains ${\mathcal U}_2$ and ${\mathcal U}_4$, and is commutative, so it will preserve 
  the space $\rho^1_{\psi}$. So we can decompose $\rho^1_{\psi}$ into eigenspaces for 
  ${\mathcal U}_5$. The characters of ${\mathcal U}_5$ that give rise to non-zero eigenspaces in
  $\rho^1_{\psi}$ will restrict to $\psi$ on ${\mathcal U}_5 \cap {\mathcal U}_1$, and to the the character $\lambda$
  on ${\mathcal U}_4$. Such a character of ${\mathcal U}_5$ would have rank $k+1$, contradicting the 
  assumption that $\rho$ has $U$-rank only $k$. Thus, the claim is proven, and 
  $\rho$ must have tensor rank at most $2k$. \Qed
 
 \vskip .1 in

Although Proposition 11.1 seems to give a rather loose bound on tensor rank in terms of $U$-rank, it can be used to verify the following agreement:

\vskip .1 in

{\bf Proposition 11.2.} Let $\rho$ in $\widehat{GL}_n$ be a representation of $U$-rank
$k$ with  $k < {n \over 4}$. Then the tensor rank of $\rho$ is equal to $k$.

\vskip .1 in
 
{\it Proof.} We know from Proposition 11.1 that the tensor rank of $\rho$ can be at most $2k$.
Since $k < {n \over 4}$, it follows that the tensor rank $\ell$ of $\rho$ is less than ${n \over 2}$, 
and therefore, it appears in the representation associated to the dual pair $(GL_n, GL_{\ell})$,
and not in the representations for dual pairs $(GL_n, GL_m)$ for $m < \ell$. But then it follows from Part c) of Theorem 5.3.1, that the $U$-rank of $\rho$ is also $\ell$. Thus $\ell = k$, as desired. \Qed

\vskip .1 in

Our next goal is Theorem 11.4. We will use the following Lemma:

\vskip .07 in

{\bf Lemma 11.3.} Suppose $n<m$. Then,
\vskip .05 in
{\bf a)} Representations of $GL_m$ of tensor rank less than or equal to $k$ restrict to representations of $GL_n$ of tensor rank less than or equal to $k$. 
\vskip .05 in
{\bf b)} All representations of $GL_n$ of tensor rank at most $k$ appear in 
restrictions of representations of $GL_m$ of tensor rank at most $k$.

\vskip .1 in

{\it Proof.} Part a) follows from the characterization of tensor rank in terms of eigenvectors for the subgroups $H_k$. Part b) follows directly from the definition of tensor rank\Qed

\vskip .1 in

Now, using Proposition 11.2 we may conclude the following:

\vskip .1 in

{\bf Theorem 11.4 (Agreement between asymptotic and tensor ranks).} For $GL_n$, asymptotic rank and tensor rank agree for all representations. 

\vskip .1 in
 
{\it Proof.} First recall from \S 6.2 that the asymptotic rank $k$ of a representation $\rho$ of $GL_n$ is the smallest number such that, for all sufficiently large $m$, 
there is a representation $\rho_m$ of $GL_m$, of $U$-rank $k$, and such that the restriction to $GL_n$ of $\rho_m$ contains $\rho$. If $m > 4k$, Proposition 11.2 
tells us that $\rho_m$ would also be of tensor rank $k$. Hence, by Part a) of Lemma 11.2, $\rho$ has tensor rank at most $k$.

By the minimality condition in the definition of asymptotic rank, we know that representations of $GL_m$ of $U$-rank less than $k$ will not contain $\rho$ when restricted to $GL_n$. Again by the condition $m > 4k$, these representations are all the tensor rank less than $k$ representations of $GL_m$, and, by Part b) of Lemma 11.2, their restrictions to $GL_n$ will contain all representations of tensor ranks less than $k$. We conclude that $k$ is also the tensor rank of $k$, as desired. \Qed

  \vskip .3 in 

  \centerline {\bf References}
  
  \vskip .1 in

[Adams-Moy93] Adams J. and Moy A., Unipotent representations and reductive dual pairs over finite fields. {\it Trans. AMS 340 (1993), 309-321}. 
  
  \vskip .05 in

[Artin57] Artin E., Geometric Algebra. {\it Interscience Publisher (1957).} 
 
  \vskip .05 in 
  
[Aubert-Kraskiewicz-Przebinda16] Aubert A.M., Kraskiewicz W., and Przebinda T., Howe correspondence and Springer correspondence for dual pairs over a finite field. {\it Lie Algebras, Lie Superalgebras, Vertex Algebras and Related Topics, Proceedings of Symposia in Pure Mathematics 92,  AMS (2016)}.
  
  \vskip .05 in 
  
[Aubert-Michel-Rouquier96] Aubert A.M., Michel J., and Rouquier R., Correspondance de Howe pour les groupesr\'eductifs sur les corps finis.  {\it Duke Math. J.  83 (1996), 353-397}.
 
  \vskip .05 in 
  
[Auslander-Kostant71]  Auslander L., and Kostant B., Polarization and unitary representations of solvable Lie groups. {\it Invent. Math. 14 (1971), 255-354}. 
  
  \vskip .05 in
  
[Bernat-Conze-Duflo-L\'evy-Nahas-Rais-Renouard-Vergne72] Bernat P., Conze N., Duflo M., L\'evy-Nahas M, Rais M., Renouard P., Vergne M., {\it Repr\'esentations des groupes de Lie r\'esolubles. Monographies de la Soci\'et\'e  Math\'ematique de France,  No. 4., Paris (1972)}. 
  
   \vskip .05 in 
  
[Borel91] Borel A.. Linear algebraic groups. {\it Springer-Verlag, New York (1991)}.
  
\vskip .05 in

[Brauer-Weyl35] Brauer R. and Weyl H., Spinors in n dimensions. {\it American J. of Math, Vol. 57 - 2, 425--449 (1935)}. 
  
 \vskip .05 in 
  
[Carter93] Carter R., Finite Groups of Lie Type: Conjugacy Classes and Complex Characters. {\it Wiley Classics Library (1993)}.
  
  \vskip .05 in
  
[Casselman-Hecht-Milicic98] Casselman W., Hecht H., and Milicic D., Bruhat filtrations and Whittaker vectors for real groups. {\it The mathematical legacy of Harish-Chandra (Baltimore, MD, 1998), 151-190, Proc. Sympos. Pure Math. 68, AMS, Providence (2000).}
  
  \vskip .05 in
  
[Ceccherini-Silberstein-Scarabotti-Tolli10] Ceccherini-Silberstein T., Scarabotti F., and Tolli F., Representation Theory of the Symmetric Groups: The Okounkov-Vershik Approach, Character Formulas, and Partition Algebras. Cambridge Studies in Advanced Mathematics 121 (2010).
 
\vskip .05 in
  
[Curtis-Reiner62] Curtis C.W. and Reiner I., Representation theory of finite groups and associative algebras. {\it Interscience, New York, 1962.}
  
  \vskip .05 in 
  
[Diaconis88] Diaconis P., Group representations in probability and statistics. {\it IMS LNM. Ser.  vol. 11 (1988)}.
  
  \vskip .05 in 
  
[Diaconis96] Diaconis P., The cutoff phenomenon in finite Markov chains. {\it PNAS, vol. 93, 1659-1664 (1996)}.
  
  \vskip .05 in 
     
[Diaconis-Shahshahani81] Diaconis P. and Shahshahani M., Generating a random permutation with random transpositions. {\it Z. Wahrsch. 57 (1981), 159-179}.

  \vskip .05 in
  
[Duke-Howe-Li92] Duke W., Howe R., and Li J.S.,  Estimating Hecke eigenvalues of Siegel modular forms. {\it Duke Math. J. 67 (1992) no. 1, 219-240}.
  
  \vskip .05 in

[Feit82] Feit W., The representation theory of finite groups. {\it North-Holland Mathematical Library 25, Amsterdam-New York (1982)}.

\vskip .05 in

[Frobenius1896] Frobenius F.G., Über Gruppencharaktere. {\it Sitzber. Preuss. Akad. Wiss. (1896) 985-1021.}

  \vskip .05 in 
    
[Fulton97] Fulton W., Young Table. Lon. Mat. Soc. Stu. Tex. 35, Cambridge University Press (1997).
  
  \vskip .05 in   
  
[Gan-Takeda16] Gan W. T. and Takeda S., A proof of the Howe duality conjecture. {\it JAMS 29 (2016), 473-493}.  
  
  \vskip .05 in 

[Gelbart75] Gelbart S., Automorphic forms on ad\`ele groups. {\it Annals of Mathematics Studies, No. 83. Princeton University Press (1975).}

  \vskip .05 in 

[Gelbart77] Gelbart S., Examples of dual reductive pairs. {\it Automorphic forms, representations and L-functions, Proc. Sympos. Pure Math., Oregon State Univ., Corvallis, Ore., (1977), Part 1, pp. 287-296}.

  \vskip .05 in 
  
[Gelfand70] Gelfand S.I., Representations of the full linear group over a finite f]ield. {\it Math. USSR-Sb., 12 (1970) 13-39}.

  \vskip .05 in 
  
[G\'erardin77] G\'erardin P., Weil representations associated to finite fields. {\it J. Alg. 46 (1977), 54-101}.
  
  \vskip .05 in
  
[Gluck97] Gluck D., Characters and Random Walks on Finite Classical Groups. {\it Adv. in Math. 129 (1997), 46-72}.

  \vskip .05 in

[Green55] Green J.A., The characters of the finite general linear groups. {\it TAMS 80 (1955) 402-447}.

  \vskip .05 in
  
[Guillemin-Sternberg79] Guillemin V.W. and Sternberg S., Some problems in integral geometry and some related problems in micro-local analysis. {\it Amer. J. Math. 101 (1979), 915-955}.  

\vskip .05 in

[Guralnick-Larsen-Tiep17] Guralnick R., Larsen M., and Tiep P.H., Character Levels and Character Bounds. {\it https://arxiv.org/abs/1708.03844 (2017).} 
 
    \vskip .05 in
  
[Guranick-Malle14] Guralnick R. and Malle G., Rational rigidity for $E_8(p)$. {\it Comp. Math., v.150 (2014)}.

  \vskip .05 in 
   
[Gurevich-Hadani07] Gurevich S. and Hadani R., The geometric Weil representation. {\it Selecta Math. 13 (2007), 465-481}.
  
  \vskip .05 in
  
[Gurevich-Hadani09] Gurevich S. and Hadani R., Quantization of Symplectic Vector Spaces over Finite Fields. {\it J. of Symp. Geom 7 (2009), 475-502.}
  
  \vskip .05 in

[Gurevich-Howe15] Gurevich S. and Howe R., Small Representations of Finite Classical Groups. {\it Representation Theory, Number Theory, 
  and Invariant Theory: In Honor of Roger Howe on the Occasion of his 70th Birthday (2015), Progress in Mathematics (2017)}.    

 \vskip .05 in

[Gurevich-Howe17] Gurevich S. and Howe R., Harmonic Analysis on $GL_n$ over Finite Fields. {\it In Preparation (2017)}.
 
\vskip .05 in
    
[Gurevich-Howe18] Gurevich S. and Howe R., A Look at Representations of $SL(2,q)$ through the Lens of Size. {\it Joe Wolf's 80s vol. (2018)}.  
   \vskip .05 in
   
[Gurevich17] Gurevich S., Small Representations of finite Classical Groups. {\it Special Lecture at IAS, March 8, 2017: https://video.ias.edu/specialrep/2017/0308-ShamgarGurevich}.  
  
\vskip .05 in

[Harish-Chandra70] Harish-Chandra., Eisenstein series over finite fields. {\it Functional analysis and related fields, Springer (1970) 76-88}. 

\vskip .05 in

[Hildebrand92] Hildebrand M., Generating random elements in $SL(n,q)$ by random transvections. J. Alg. Comb. 1 (1992),133-150.  
  
 \vskip .05 in 
 
[Howe73-1] Howe R.,  On the character of Weil's representation. {\it Tran. AMS. 177 (1973), 287-298.}
 
   \vskip .05 in 
  
[Howe73-2] Howe R., Invariant theory and duality for classical groups over finite fields with applications to their singular representation theory.{\it Preprint, Yale University (1973).}
 
   \vskip .05 in 
 
[Howe77] Howe R.,  $\theta$-series and invariant theory. {\it Automorphic forms, representations and L-functions, Proc. Sympos. Pure Math., Oregon State Univ., Corvallis, Ore. (1977), Part 1, pp. 275-285}.
  
  \vskip .05 in
  
[Howe80] Howe R., On the role of the Heisenberg group in harmonic analysis. {\it Bul. Amer. Math. Soc., Vol. 3, Num. 2 (1980) 821-843.}
   
  \vskip .05 in   
  
[Howe82] Howe R., On a notion of rank for unitary representations of the classical groups. {\it Harmonic analysis and group representations, 223-331, Liguori, Naples, (1982).}
   
 \vskip .05 in    

[Howe85] Howe R., Dual pairs in physics: Harmonic oscillators, photons, electrons, and singletons. {\it Lectures in Applied Mathematics, Vol. 21, AMS, (1985) 179-207}.

 \vskip .05 in

[Howe87] Howe R.. The oscillator semigroup. {\it The mathematical heritage of Hermann Weyl. Proc. Sympo. Pure Math 48, AMS (1988) 61-132}. 

 \vskip .05 in   
   
[Howe17-1] Howe R., On the role of rank in representation theory of the classical groups. {\it Special Lecture at IAS, March 8, 2017: https://video.ias.edu/reptheory/2017/0308-RogerHowe}.    
   
   \vskip .05 in 
   
[Howe17-2] Howe R., Rank and Duality in Representation Theory. {\it Video of the 19th Takagi lectures, RIMS - Kyoto University, Japan, July 8-9 (2017):}  
   
$http://www.ms.u-tokyo.ac.jp/~toshi/jjm/JJM\_HP/contents/takagi/19th/takagi19-video.htm$.
  
\vskip .05 in

[Howe-Moy86] Howe R. and Moy A., Harish-Chandra homomorphisms for p-adic groups. {\it CBMS Regional Conference Series in Mathematics 59 (1986)}.
 
\vskip .05 in 
  
[Isaacs76] Isaacs I.M., Character theory of finite groups. {\it Pure and Applied Mathematics, No. 69. Academic Press, New York-London (1976).}
  
  \vskip .05 in 
  
[Jacobson62] Jacobson N., Lie Algebras. {\it Dover Publications (1979)}.
  
  \vskip .05 in 
  
[Jacquet67] Jacquet H., Fonctions de Whittaker associées aux groupes de Chevalley. {\it Bulletin de la Société Mathématique de France, 95 (1967) 243–309}.

\vskip .05 in 
  
[Kirillov62] Kirillov A.A.,  Unitary representations of nilpotent Lie groups. {\it Uspehi Mat. Nauk 17 (1962), no. 4, 57-110}.
  
  \vskip .05 in 
  
[Kirillov04] Kirillov A.A., Lectures on the orbit method. {\it GSM 64, AMS, Providence (2004).}
  
  \vskip .05 in 
  
[Knapp86] Knapp A., Representation theory of semisimple groups. An overview based on examples. 
  {\it Princeton Mathematical Series 36, Princeton University Press (1986).}
  
  \vskip .05 in 
  
[Kostant78] Kostant B., On Whittaker Vectors and Representation Theory. {\it Inventiones mathematicae 48 (1978) 101-184}.
  
  \vskip .05 in  
  
[Kudla86] Kudla S., On the local theta-correspondence. {\it Invent. Math. 83 (1986), no. 2, 229-255.}
  
  \vskip .05 in
  
[Kudla-Milson86] Kudla S. and Millson J., The theta correspondence and harmonic forms. {\it Math. Ann. 274 (1986), no. 3, 353-378.}
  
  \vskip .05 in
 
[Lam73] Lam T.Y., The algebraic theory of quadratic forms. {\it MLNS. W. A. Benjamin, Mass. (1973).}

  \vskip .05 in 
  
[Lang02] Lang S, Algebra, {\it GTM 211. Springer-Verlag, New York,  (2002).}
  
  \vskip .05 in

[Larsen-Shalev-Tiep11] Larsen M, Shalev A., and Tiep P.H., The Waring problem for finite simple groups. {\it Annals of Math. 174 (2011), 1885-1950}.

\vskip .05 in

[Li89-1] Li J.S.,  Singular unitary representations of classical groups. {\it Invent. Math. 97 (1989), no. 2, 237-255}.
  
  \vskip .05 in
  
[Li89-2] Li J.S., On the classification of irreducible low rank unitary representations of classical groups. {\it Compositio Math. 71 (1989), no. 1, 29-48.}
  
  \vskip .05 in 
  
[Li92] Li J.S.,  Nonexistence of singular cusp forms. {\it Compositio Math. 83 (1992), no. 1, 43-51.}
   
 \vskip .05 in
 
[Liebeck17] Liebeck M.W., Character ratios for finite groups of Lie type, and applications. {\it Contemporary Mathematics 694 (2017), American Math. Soc., 193-208.} 
 
 \vskip .05 in

[Liebeck-O'Brien-Shalev-Tiep10] Liebeck M.W., O'Brien E.A., Shalev A., and Tiep P.H., The Ore conjecture. {\it JEMS 12 (2010), 939-1008}.
 
 \vskip .05 in
 
[Liebeck-Shalev05] Liebeck M.W. and Shalev A., Character degrees and random walks in finite groups of Lie type. {\it Proc. LMS. 90 (2005), 61-86}.

\vskip .05 in  
  
[Lusztig84] Lusztig G., Characters of reductive groups over a finite field. {\it Annals of Mathematics Studies 107. Princeton University Press (1984)}.
  
  \vskip .05 in 
  
[Mackey49] Mackey G., A theorem of Stone and von Neumann. {\it Duke Math. J. 16, (1949). 313-326}. 

  \vskip .05 in 
  
[Malle14] Malle G., The proof of Ore's conjecture. {\it S\'em. Bourbaki, Ast\'erisque 361 (2014), 325-348}.
 
  \vskip .05 in   
  
[Matsumoto69] Matsumoto H., Sur les sous-groupes arithm\'etiques des groupes semi-simples d\'eploy\'es. {\it Ann. Sci. Ec. Norm. Sup 2 1969) 1-62.}   
   
\vskip .05 in 

[Moeglin-Vigneras-Waldspurger87] Moeglin C., Vigneras M.F., and Waldspurger J.L., Correspondence de Howe sur un corps p-adique. LNM 1291, Springer-Verlag (1987).

\vskip .05 in
  
[Nazarov95] Nazarov M., The Oscillator Semigroup over a non-archimedean field. {\it Journal of functional analysis 128, 384-438 (1995)}.      
  
\vskip .05 in   
  
[Mostow-Sampson69] Mostow G. and Sampson J., Linear algebra. {\it McGraw-Hill Book Co. (1969).}

  \vskip .05 in

[Prasad93] Prasad D., Weil representation, Howe duality, and the theta correspondence. {\it Theta functions: from the classical to the modern, 105-127, 
  CRM Proc., AMS, Providence (1993).}

  \vskip .05 in 
  
[Pukanszky67] Pukanszky L., Le\c cons sur les repr\'esentations des groupes. {\it Monographies de la Soci\'et\'e Math\'ematique de France, No. 2, Paris (1967).} 
 
  \vskip .05 in 
  
[Pukanszky78]  Pukanszky L., Unitary representations of Lie groups with co-compact radical and applications. {\it Tran. AMS. 236 (1978), 1-49.} 
  
  \vskip .05 in
 
[Sahi17] Sahi S., Private communication. {\it Weizmann Institute - Israel, June 2017.}   
 
 \vskip .05 in 
    
[Saloff-Coste04] Saloff-Coste L., Random Walks on Finite Groups. {\it vol. 110, Encyclopaedia of Math. Sciences (2004), 263-346}. 
        
  \vskip .05 in
  
 [Serre77] Serre J.P., Linear Representations of Finite Groups. {\it Springer-Verlag (1977).}
  
  \vskip .05 in
  
  [Shalev07]  Shalev A., Commutators, Words, Conjugacy Classes and Character Methods. {\it Tur. J Math 31 (2007), 131-148}. 
  
  \vskip .05 in
  
  [Shalev09] Shalev A., Word maps, conjugacy classes, and a non-commutative Waring-type theorem. {\it Annals of Math., 170 (2009), 1383-1416}. 
  
  \vskip .05 in 
  
  [Shalev15] Shalev A., Conjugacy classes, growth and complexity. {\it Finite Simple Groups: Thirty Years of the Atlas and Beyond (2015), Cont. Math. 694 (2017) 209-221.} 

  \vskip .05 in    
  
[Shalika74] Shalika J.A., The Multiplicity One Theorem for GLn. {\it Annals of Math 100 no. 2, (1974) 171-93}. 

  \vskip .05 in 

 [Srinivasan79] Srinivasan B., Weil representations of classical groups. {\it Invent. Math. 51 (1979), 143-153.}    
    
  \vskip .05 in    
  
[Steinberg62] Steinberg R., G\'en\'erateurs, relations et rev\^etements de groupes alg\'ebriques. {\it Colloque sur la Th\'eoriedes groupes Alg\'ebriques, Bruxelles (1962) 113-127.}  
  
  \vskip .05 in   
  
[Strang76] Strang G. , Linear algebra and its applications. {\it Academic Press, New York-London (1976).} 
  
  \vskip .05 in
  
[Vergne83] Vergne M., Representations of Lie Groups and the Orbit Method. {\it Emmy Noether in Bryn Mawr. Springer (1983), 59-101.} 
 
  \vskip .05 in 
  
[Vogan80] Vogan D.A., The size of infinite-dimensional representations. {\it LNM vol 795, Springer (1980), 172-178.}  
  
  \vskip .05 in 
  
[Vogan16] Vogan D.A., The size of infinite-dimensional representations. {\it Takagi Lectures (2016), Jap. J of Math, vol 12 (2017), 175-210.} 
  
  \vskip .05 in 
  
[Wallach82] Wallach N., Real reductive groups. I. {\it Pure and Applied Mathematics, 132. Academic Press, Boston, (1988).} 

 \vskip .05 in 

[Wallach92] Wallach N., Real reductive groups. II. {\it Pure and Applied Mathematics, 132-II. Academic Press, Boston (1992).} 

\vskip .05 in 

[Weil64] Weil A., Sur certains groupes d'op\'erateurs unitaires. {\it Acta Math. 11 1964, 143-211.} 

\vskip .05 in 

[Weil74] Weil  A., Basic number theory. {\it Die Grundlehren der Mathematischen Wissenschaften, Band 144. 
Springer-Verlag (1974).}

\vskip .05 in

[Weinstein80] Weinstein A.D., The symplectic ``category". {\it Proc. Conf. Differential Geometric Methods in Mathematical Physics (1980)}.

\vskip .05 in 

[Weyl46] Weyl H., The classical groups, their invariants and representations. {\it Princeton Univ. Press (1946)}.

\vskip .05 in 

[Zelevinsky81] Zelevinsky A., Representations of finite classical groups. A Hopf algebra approach. {\it LNM 869. Springer-Verlag (1981)}.

\end{document}